\documentclass[11pt]{amsart}

\usepackage{style}
\usepackage{macros}

\author[A. Leroux-Lapierre]{Alexis Leroux-Lapierre}
\address{A. Leroux-Lapierre \\McGill University, Montréal, Québec}
\email{alexis.leroux-lapierre@mail.mcgill.ca}

\begin{document}

\begin{abstract}
This paper defines an asymptotic character map which is a morphism from the Grothendieck group of category $\mathcal{O}$ of an integral filtered quantization to rational functions on the Lie algebra of a torus. We show that the asymptotic character of a module computes the equivariant multiplicity of its characteristic cycle. We then apply this construction to truncated shifted Yangians coming from simple, simply-laced Lie algebras and draw connections with characters of modules over KLR algebras using an equivalence of categories of \cite{kamnitzer2019category}. Our main theorem shows how this new formalism gives formulas relating equivariant multiplicities of Mirkovi\'{c}--Vilonen cycles and characters of modules over cyclotomic KLR algebras. We explain how this result provides evidence that the change-of-basis between Lusztig's dual canonical basis and the Mirkovi\'{c}--Vilonen basis of $\mathbb{C}[N]$ is computed by a characteristic cycle map whose domain is category $\mathcal{O}$ for truncated shifted Yangians, implying that the coefficients are non-negative integers.
\end{abstract}

\title{Category $\cO$ and asymptotic characters}
\maketitle

{
	\setcounter{tocdepth}{1}
	\hypersetup{linkcolor=black}	
	\tableofcontents
}

\setlength{\medskipamount}{0pt}
\setlength{\smallskipamount}{0pt}

\setlength{\abovedisplayskip}{\globalspaceeq}
\setlength{\belowdisplayskip}{\globalspaceeq}
\setlength{\abovedisplayshortskip}{\globalspaceeq}
\setlength{\belowdisplayshortskip}{\globalspaceeq}

\setlength{\belowcaptionskip}{\globalspacethm} 			
\setlength{\textfloatsep}{\globalspacethm}	
\setlength{\parskip}{\globalspace}			

\section{Introduction}

\subsection{Perfect bases in representation theory}

Let $G$ be a simple, simply connected, complex algebraic group. Understanding the combinatorics of the family of irreducible representations $\{V(\lambda)\}_{\lambda\in P_+}$ (e.g.~their underlying crystal graph or tensor product multiplicity formulae) has been a focal point of many Lie theoretical problems. One strategy employed to study such combinatorial questions is to produce bases which have a constrained behaviour under the $\fg$-action. For example, a first step towards this goal is to work with weight bases. For a historical review of the different criteria on bases which proved useful, the reader is referred to \cite{kamnitzer2022perfect}. The present paper is concerned with the study of \textit{perfect bases} for irreducible representations of $G$, as introduced in \cite{berenstein2006geometric}. By definition, a perfect basis has a well-behaved action of the Chevalley generators of $\fg=\Lie(G)$ which is controlled by the crystal of the underlying representation. Following \cite{baumann2021mirkovic}, we focus on \textit{biperfect bases} of the coordinate ring $\C[N]$, where $N$ is the unipotent radical of a chosen Borel subgroup. Through embeddings, biperfect bases give perfect bases of all irreducible representations at once. As of the time of writing this article, there are three known construction of biperfect bases of $\C[N]$: 
\begin{enumerate}[label={\color{colorOfLinks}\bf b\arabic*}]
\item\label{b:DCbasis} The dual canonical basis (or upper global basis), due to \cite{lusztig1990canonicalI,lusztig1990canonicalII} and \cite{kashiwara1993global}. Work of \cite{varagnolo2011canonical} demonstrates that in the simply-laced case, this basis can be reinterpreted using representation theory of KLR algebras\footnote{Note that in general, representation theory of KLR algebras over a field of positive characteristic provides other examples of biperfect bases, see \cite[Section~12]{mcnamara2019folding}} (sometimes called quiver Hecke algebras) of \cite{khovanov2009diagrammatic,rouquier20082,khovanov2011diagrammatic}.
\item\label{b:MVbasis} The Mirkovi\'{c}--Vilonen basis, due to \cite{mirkovic2007geometric} together with \cite{ginzburg1995perverse} and \cite{baumann2021mirkovic}. This basis is defined using the geometry of affine Grassmannian slices. 
\item\label{b:SDCbasis} The dual semicanonical basis, due to \cite{lusztig2000semicanonical}. This basis is defined using representation theory of preprojective algebras. 
\end{enumerate}
Work of \cite{baumann2021mirkovic} shows that all three bases are distinct for general $G$, hence posing the obvious question of comparing them. \par 

It has been highlighted in many recent papers that the above defined biperfect bases and their interactions encode rich combinatorics. We name a few notable such papers here. Work of \cite{casbi2021equivariant} draws connections between the cluster algebra structure on $\C[N]$ defined in \cite{berensetein2005cluster}, combinatorics of good Lyndon words (via basis \ref{b:DCbasis} and \cite{kleshchev2011representations,mcnamara2015finite}) and smoothness of Mirkovi\'{c}--Vilonen (MV) cycles (via basis \ref{b:MVbasis}). Subsequent work of \cite{casbi2023quantum,casbi2024auslander} investigates further links between MV cycles and the combinatorics of $q$-characters based on tools from \cite{hernandez2016cluster}. In a different direction, \cite{baumann2024bases} studies multiplicative properties of the basis \ref{b:MVbasis} and its relationship with cluster monomials of $\C[N]$. Furthermore, as all three bases give a model for the crystal $\cB(\infty)$, they encode combinatorics of MV polytopes of \cite{anderson2003polytope}, a family of polytopes which was explicitly described in \cite{kamnitzer2010mirkovic}. Finally, we remark that the study of biperfect bases of $\C[N]$ is related to the trichotomy of bases for cluster algebras discussed in \cite{qin2023cluster}. Under some technical assumptions (which the cluster structure on $\C[N]$ satisfies), cluster algebras admit three remarkable bases: the common triangular basis (due to \cite{qin2017triangular}), the theta basis (due to \cite{gross2018canonical}) and the generic basis (due to \cite{dupont2011generic}). In the case of $\C[N]$, \cite{qin2024bases} shows that basis \ref{b:DCbasis} coincides with the common triangular basis while \cite{geiss2012generic} shows that basis \ref{b:SDCbasis} coincides with the generic basis. It is conjectured that basis \ref{b:MVbasis} coincides with the theta basis, see \cite[Question~7.1]{qin2023cluster} and \cite[Conjecture~6.9]{kamnitzer2022perfect}.\par 

This article tackles the question of comparing the biperfect bases \ref{b:DCbasis} and \ref{b:MVbasis}. To establish our results, we use category $\cO$ for truncated shifted Yangians. \par 

\subsection{Shifted Yangians, the affine Grassmannian, and KLR algebras}

Shifted Yangians are a family of infinite dimensional shifted quantum groups associated to $\fg$ which depend on a choice of coweight $\mu$. They appeared first in type $A$ (for dominant $\mu$) in the work of \cite{brundan2006shifted} and were generalized to arbitrary type (for dominant $\mu$) in \cite{kamnitzer2014yangians}, where they were identified as quantizations of coordinate rings of subvarieties of the affine Grassmannian. Crucially, quotients of the shifted Yangians (which further depend on a dominant coweight $\lambda$) --- the truncated shifted Yangians --- were shown to quantize coordinate rings of affine Grassmannian slices. This result placed truncated shifted Yangians in the paradigm of quantizations of conical symplectic resolutions \cite{braden2012quantizations,braden2014quantizations} when $\lambda$ is a sum of dominant minuscule coweights. The definition of these algebras was later extended to arbitrary coweights in \cite{braverman2016coulomb}, where they were also identified as quantizations of Coulomb branches of $3d$ $\mathcal{N}=4$ quiver gauge theories in the simply-laced case. The non-simply laced case was carried out in \cite{nakajima2023coulomb}, where it was shown that they are quantizations of Coulomb branches of quiver gauge theories with symmetrizers.

Since then, the representation theory of shifted Yangians and their truncations has developed rapidly, particularly the study of representations in category $\cO$ (the analogue of the classical BGG category $\cO$).
A systematic study of category $\cO$ for shifted Yangians was initiated in \cite{kamnitzer2019highest} and in \cite{hernandez2021shifted}. 
Subsequently, it was proven in \cite{kamnitzer2019category,kamnitzer2022lie} that, given a truncated shifted Yangian, its category $\cO$ is equivalent to a full subcategory of the category of modules over a KLRW algebra of \cite{webster2017knot}. 
This endowed each of those category $\cO$ with a categorical $\fg$-action via transport de structure.

As the previous paragraphs have highlighted, truncated shifted Yangians are related both to the geometry of affine Grassmannian slices (via quantization) and to representation theory of KLR algebras (via equivalences of categories). The aim of this paper is to explain how the quantization process encodes discrepancies between bases \ref{b:DCbasis} and \ref{b:MVbasis}. To carry this out, we use a tool called the \textit{characteristic cycle} which is an algebraic analogue of the characteristic cycle of \cite{kashiwara2010deformation}. For a given module over a truncated shifted Yangian, its characteristic cycle is the scheme-theoretical support of the associated graded module. We prove that when the module lies in category $\cO$, then the top-dimensional characteristic cycle can be expressed as a $\Z_{\geq 0}$-linear combination of its MV cycles. Hence, this construction provides a good candidate for the change-of-between bases \ref{b:DCbasis} and \ref{b:MVbasis}. While we are unable to prove this conjecture, we provide evidence using the $\Dbar$ map of \cite{baumann2021mirkovic}.

\subsection{Comparison in the limit}

In \cite{baumann2021mirkovic}, a morphism of algebras $\Dbar:\C[N]\to \C(\fh)$ was introduced. There, it was shown that applying $\Dbar$ to bases \ref{b:MVbasis} and \ref{b:SDCbasis} computes a quantity of interest. For basis \ref{b:MVbasis}, \cite[Theorem~1.4]{baumann2021mirkovic} shows that the image of a basis vector under $\Dbar$ is the equivariant multiplicity of the associated MV cycle at its bottom fixed point. This allowed for a proof of the main conjecture of \cite{muthiah2021weyl}. For basis \ref{b:SDCbasis}, \cite[Theorem~11.4]{baumann2021mirkovic} shows that the image of a basis vector under $\Dbar$ computes the asymptotic behavior of flags of submodules of the associated generic preprojective algebra module. For basis \ref{b:DCbasis}, as discussed in \cite{casbi2021equivariant}, the image of a basis vector under $\Dbar$ can be computed from the character of the associated KLR module.

Inspired by the definition of the $\Dbar$ map as a coefficient of the Fourier transform of the Duistermaat-Heckmann measure, we develop a new notion of asymptotic characters, which we denote $\achi$, for modules over truncated shifted Yangians. As the name suggests, asymptotic characters can be computed by taking appropriate limits of characters of modules. On one hand, we show that for a given module of category $\cO$ which has maximal Gelfand--Kirillov dimension, $\achi$ computes the equivariant multiplicity of its characteristic cycle. On the another hand, we prove the equivalence of categories of \cite{kamnitzer2019category} mentioned previously intertwines $\achi$ and $\Dbar$, hence showing that for a simple module of $\cO$ which has maximal Gelfand--Kirillov dimension, its asymptotic character computes $\Dbar$ of the corresponding basis vector of basis \ref{b:DCbasis}. Let us note that the above conjectural change-of-basis from basis \ref{b:DCbasis} to \ref{b:MVbasis} matches with conjectural properties of the change-of-basis from the common triangular basis to the theta basis in the setting of cluster algebras. As explained in \cite[Section~7]{qin2023cluster}, results of \cite{mandel2017theta} provide evidence that common triangular basis elements are non-negative integral linear combinations of theta basis elements.

As the formalism needed for the construction of $\achi$ is completely independent of truncated shifted Yangians and only uses properties of ``well-behaved'' filtered quantizations, theorems of Section \ref{sec:asympchar} are stated in the abstract framework of filtered quantization where some strong integrality and finiteness properties hold. The necessary hypotheses for the asymptotic characters to converge and give rational functions which encode equivariant invariants of their characteristic cycles are clearly identified.

\subsection{Outline of the paper}
{\it Section \ref{sec:generalities}} introduces the notation needed throughout the paper, reviews the definition of a biperfect basis of $\C[N]$ together with the $\Dbar$ map and states the main result without introducing truncated shifted Yangians (see Theorem \ref{thm:equalityDbar}). 
{\it Section \ref{sec:asympchar}} develops tools referred to as algebraic category $\cO$ for integral quantizations and, as a consequence, proves the first part of the main theorem by defining asymptotic characters. This is carried out by showing that asymptotic characters compute the equivariant multiplicity of their support (see Theorem \ref{thm:achiandequivmult}). 
{\it Section \ref{sec:Yangian}} is concerned with shifted Yangians and their truncations. The main point of this section is carefully check that the necessary hypotheses are satisfied so that tools of the previous section can be applied. The main result of the section appears in Theorem \ref{thm:epsilonTversusachi}. A discussion on injectivity of the characteristic cycle map is added. 
{\it Section \ref{sec:KLR}} introduces KLR(W) algebras and relates them to truncated shifted Yangians according to an equivalence of categories of \cite{kamnitzer2019category}. This equivalence is then used to prove Theorem \ref{thm:DbarMvsachi} which shows that asymptotic characters can be computed using the the bar-character map of KLR algebra modules. 
Finally, {\it Section \ref{sec:stitching}} assembles results of the previous sections into Theorem \ref{thm:theMainTheorem}. Natural conjectures related to the main results are highlighted. A brief discussion on the Kac-Moody case is included and Corollary \ref{cor:conjectureNakajima} provides evidence for a conjecture from \cite{braverman2016coulomb}. 

To help the reader, we give the main result of the paper here:  
\begin{theorem*}[Theorem \ref{thm:theMainTheorem}]
There is a commutative diagram
\begin{equation*}\adjustbox{scale=1,center}{
\begin{tikzcd}[column sep=1em,row sep=0.75em]
\HH_{\tp}\big((\cWbar{}^\lambda_\mu)_-\big)\arrow[d]  & K_0(\cO_\mu^\lambda(\bR)) \arrow[l,"(1)"'] \arrow[r,"(2)"]\arrow[dd,"(3)"'] & K_0(R_{\lambda-\mu}^\lambda\modu)\arrow[d]\\
\C[N]_{-(\lambda-\mu)}\arrow[dr,"\DbarM"'] & & \C[N]_{-(\lambda-\mu)}\arrow[dl,"\DbarM"]\\
& \C(\fh) & 
\end{tikzcd}}
\end{equation*}
where arrow $(1)$ is the characteristic cycle map from Definition \ref{def:charcycle}, arrow $(2)$ is a consequence of Theorem \ref{thm:equivTheta} and arrow $(3)$ is the asymptotic character map defined in \ref{def:asymptoticCharacterDef}. 
\end{theorem*}

\subsection{Acknowledgements}

The author is especially grateful to Joel Kamnitzer for gifting him this project and offering support, advices and key ideas throughout.

The author would also like to thank: B.~Webster for suggesting that equivariant multiplicities could be computed using KLR algebras, A.~Dranowski who was part of this project at its early stages, A.~Weekes for the enlightening discussions, for pointing out Example \ref{ex:D4monomialCrystal} and for suggesting Lemma \ref{lem:restrictionPreservesGKdim} together with Remark \ref{rem:smoothFixedPoint}, P.~Baumann for carefully explaining the details of the construction of the $\Dbar$ map appearing in Section \ref{subsec:Dbarmap} and for reading this paper and providing essential feedback, É.~Casbi for explaining the combinatorics appearing in Example \ref{ex:Nakada} and for numerous discussions on the $\Dbar$ map, D.~Gaitsgory for suggesting to consider the reduced repelling set appearing in Section \ref{subsec:fixedpoints}, F.~Qin for pointing the connections with properties of bases of cluster algebras, P.~McNamara for explaning that representation theory of KLR algebras over fields of positive characteristic yield many examples of biperfect bases, T.~Pinet and A.~Kalmykov.

The author is also in debt to the anonymous referee who pointed out many typos and provided useful comments.

This work was supported by the \textit{Fonds de recherche du Québec -- Nature et Technologies} and \textit{Institut des Sciences Mathématiques}.

\subsection{Lie-theoretic notation}\label{sec:notation}

Let $B$ be a choice of Borel subgroup of $G$, together with $T\subset B$ a maximal torus. Let $N$ denote the unipotent radical of $B$. Let $B_-$ be the opposite Borel and $N_-$ its unipotent radical. The complex Lie algebras of the groups $T,N,B,G$ are denoted respectively by $\fh,\fn,\fb$ and $\fg$. \par
 
The above choices determine a Cartan datum, which consists of 
\begin{itemize}
\item the weight lattice of $\fg$, denoted by $P$,
\item an indexing set for the simple roots of $\fg$, denoted by $I$, 
\item the set of simple roots, denoted by $\{\alpha_i\}_{i\in I}\subset\fh^\ast$, 
\item the set of simple coroots, denoted by $\{\alpha_i^\vee\}_{i\in I}\subset\fh$ and
\item a non-degenerate symmetric bilinear form on $\fh^\ast$, denoted by $(.,.):\fh^\ast\times\fh^\ast\to \C$.
\end{itemize}
We normalize the bilinear form $(.,.)$ such that $d_i=\tfrac{1}{2}(\alpha_i,\alpha_i)\in\Z_{>0}$ are coprime integers. The paring between $\fh^\ast$ and $\fh$ is denoted by $\langle.,.\rangle:\fh^\ast\times\fh\to \C$. It satisfies 
\begin{equation}\label{eq:cartanPairing}
\langle \beta,\alpha_i^\vee\rangle=\dfrac{2(\alpha_i,\beta)}{(\alpha_i,\alpha_i)}
\end{equation}
for all $\beta\in P$ and $i\in I$. The corresponding root system is denoted by $\Phi$ and positive roots (respectively negative roots) are denoted by $\Phi_+$ (respectively $\Phi_-$). Dominant weights are denoted by $P_+$. The root lattice of $\fg$ is denoted by $Q:=\bigoplus_{i\in I}\Z\alpha_i$ and the positive root cone is denoted $Q_+:=\bigoplus_{i\in I}\Z_{\geq 0}\alpha_i$. The Dynkin diagram of $\fg$ is denoted by $\Gamma$. The entries of the Cartan matrix of $\fg$ are labeled by $A=\big(a_{ij}=\langle\alpha_j,\alpha_i^\vee\rangle\big)_{i,j\in I}$. The fundamental weights are denoted by $\{\varpi_i\}_{i\in I}\subset\fh^\ast$. The rank of $\fg$ is denoted by $r=|I|$.\par 

The above data also determines a root system dual to that of $\fg$. Let $G^\vee$ denote the corresponding Langlands dual group. We use the typical $(\trou)^\vee$ notation to label the data associated to $G^\vee$.

For $\lambda\in P_+$, let $V(\lambda)$ denote the (unique) finite dimensional, irreducible representation of $G$ of highest weight $\lambda$. Fix once and for all a highest weight vector $v_\lambda\in V(\lambda)$. The linear form $v_{\lambda}^\ast\in V(\lambda)^\ast$ will denote the projection on the $1$-dimensional space $\spa_\C\{v_\lambda\}\simeq \C$ satisfying $v_\lambda^\ast(v_\lambda)=1$. For $V$ any $G$-representation and $\mu \in P$, the $\mu$-weight space of $V$ is denoted by $V_\mu$. 

Given a $\fg$-crystal $\cB$, the set of maps associated it will be denoted by $(\wt,\epsilon_i,\phi_i,\tilde{e}_i,\tilde{f}_i)_{i\in I}$, where $\wt:\cB\to P$ is the weight map and for all $i\in I$, one has maps $\epsilon_i,\phi_i:\cB\to \Z\sqcup\{-\infty\}$ as well as crystal operators $\tilde{e}_i,\tilde{f}_i:\cB\to\cB\sqcup\{0\}$ satisfying the usual axioms (see for example \cite[Definition~2.13]{bump2017crystal}).

Given a $\C$-algebra $A$, $A\Modu$ denotes the category of left $A$-modules, $A\modu$ denotes the category of finitely generated left $A$-modules and $A\fmodu$ denotes the category of finite dimensional left $A$-modules. Unless otherwise specified, modules are left modules and algebras are assumed unital.

Throughout, we use parity convention of \cite[Section~2.1]{kamnitzer2019category}, namely:

\begin{convention}\label{con:parity}
Suppose that $\Gamma$ is simply-laced. For $i,j\in I$, write $i\sim j$ if $i$ is connected to $j$. Fix a decomposition $I=I_{0}\sqcup I_{1}$ into even and odd vertices which makes $\Gamma$ a bipartite graph (i.e. no edge has vertices of the same parity). Endow $\Gamma$ with an orientation such that even arrows are sources and odd arrows are targets. Now, fix a total order on $I$ such that even vertices are smaller than odd vertices, say $I=\{i_1<i_2<\dots<i_r\}$.
\end{convention}

\section{Perfect bases and the $\Dbar$ map}\label{sec:generalities}

In this section, we state some well-known definitions and results concerning (bi)perfect bases to give the reader a precise statement of the main result of the paper (Theorem \ref{thm:equalityDbar}). Elementary examples which we will refer to throughout the following sections are included.

\subsection{Biperfect bases}

For each $i\in I$, fix $h_i:=\alpha_i^\vee$ and choose Chevalley generators, i.e.~weight vectors $e_i,f_i\in \fg$ of weight $\alpha_i,-\alpha_i$ respectively satisfying $[e_i,f_i]=h_i$.

\begin{definition}[{\cite[Definition~2.5]{baumann2021mirkovic}}]
A basis $B_\lambda$ of $V(\lambda)$ is called a perfect basis if $B_\lambda$ is equipped with an upper semi-normal crystal structure $(\wt,\epsilon_i,\phi_i,\tilde{e}_i,\tilde{f}_i)$ such that 
\begin{enumerate}
\item $v_\lambda\in B_\lambda$, 
\item $B_\lambda$ is a weight basis and
\item for each $i\in I$ and $b\in B_\lambda$, one has an explicit description of the action of each $e_i\in\fg$ as
\begin{equation}\label{eq:actionPerfBasis}
e_i b=\epsilon_i(b) \tilde{e}_i(b) +\sum_{\substack{b'\in B_\lambda \\ \epsilon_i(b')<\epsilon_i(b)-1}} a_{b'} b'
\end{equation}
for some $a_{b'}\in \C$.
\end{enumerate}
\end{definition}

Note that the choice of Chevalley generators $e_i\in \fn$ impacts the basis $B_\lambda$.

Recall that the coordinate ring $\C[N]$ admits a $Q_+$-grading via conjugation by $T$, i.e.~for $f\in \C[N]$, let $(t\cdot f)(n)=f(tnt^{-1})$ for all $t\in T$ and $n\in N$.
Moreover, $\fn$ acts both on the right and on the left of $\C[N]$ by differential operators.
This action provides a pairing
\begin{equation}\label{eq:paringUnCN}
\begin{aligned}
\langle.,.\rangle:U(\fn)\times \C[N]&\to \C\\
(x,f)&\mapsto(x\cdot f)(1_N)
\end{aligned}
\end{equation}
which is a $Q_+$-graded perfect pairing.

One way of producing perfect bases simultaneously for all $V(\lambda)$'s is to produce a \textit{biperfect basis} for the coordinate ring $\C[N]$. Let $\psi_\lambda:V(\lambda)\to \C[N]$ be the map defined by 
\begin{equation}\label{eq:defpsilambda}
\big(\psi_\lambda(v)\big)(n)=v_\lambda^\ast(nv).
\end{equation}
It is routine to show that $\psi_\lambda$ is injective, $\fn$-equivariant and that it sends vectors of weight $\mu$ to functions of weight $\lambda-\mu$.

\begin{definition}[{\cite[Definition~2.2]{baumann2021mirkovic}}]\label{def:biperfectBasis}
A basis $B$ of $\C[N]$ is called biperfect if it is a perfect basis and if $B$ is equipped with a second upper semi-normal crystal structure $(\wt,\epsilon_i^\ast,\phi_i^\ast,\tilde{e}_i^\ast,\tilde{f}_i^\ast)$ sharing the same weight map as the first and such that 
\begin{equation*}
b e_i=\epsilon_i^\ast(b) \tilde{e}_i^\ast(b) +\sum_{\substack{b'\in B_\lambda \\ \epsilon_i^\ast(b')<\epsilon_i^\ast(b)-1}} a_{b'}^\ast b'
\end{equation*}
for some $a_{b'}^\ast\in\C$
\end{definition}

\begin{proposition}[{\cite[Proposition~2.7]{baumann2021mirkovic}}]
If $B$ is a biperfect basis of $\C[N]$, then $B_\lambda=\psi_\lambda^{-1}(B)$ is a perfect basis of $V(\lambda)$ (when endowed with the induced crystal operators which are compatible with the left action of $\mathfrak{n}$ on $\C[N]$).
\end{proposition}

The following example demonstrates how the previous proposition works. Throughout the paper, we will frequently come back to it to illustrate results.

\begin{example}\label{ex:C[N]SL3}
For $G=\SL_3(\C)$, consider the usual choices for $T,N,B$ and denote simple roots by $\alpha_1,\alpha_2$. In this case,
\begin{equation*}
N=\big\{\begin{smallpmatrix}
1 & a & c\\ 0 & 1 & b\\ 0 & 0 & 1
\end{smallpmatrix}\;;\; a,b,c\in\C \big\}.
\end{equation*}
The Chevalley generators $e_1=\begin{smallpmatrix}
0 & 1 & 0\\
0 & 0 & 0\\
0 & 0 & 0
\end{smallpmatrix},e_2=\begin{smallpmatrix}
0 & 0 & 0\\
0 & 0 & 1\\
0 & 0 & 0
\end{smallpmatrix}\in \fn$ act on $\C[N]$ on the left and on the right by the differential operators
\begin{equation*}
\Big(\, e_1\mapsto \dfrac{\partial}{\partial a},\hspace{0.5em} e_2\mapsto\dfrac{\partial}{\partial b}+a \dfrac{\partial}{\partial c}\,\Big)\hspace{1em}\text{and}\hspace{1em}\Big(e_1\mapsto \dfrac{\partial}{\partial a}+b\dfrac{\partial}{\partial c},\hspace{0.5em} e_2\mapsto\dfrac{\partial}{\partial b}\Big)
\end{equation*}
respectively. It is easy to check that the set 
\begin{equation*}
B:=\Big\{a^{n_1}c^{n_2}(ab-c)^{n_3}\;;\; n_i\in\Z_{\geq 0}\Big\}\cup \Big\{ b^{n_1}c^{n_2}(ab-c)^{n_3}\;;\; n_i\in\Z_{\geq 0}\Big\}
\end{equation*}
is a biperfect basis. Now, for the adjoint representation $V(\varpi_1+\varpi_2)\simeq \sl_3(\C)$, one can pick the highest weight vector $\begin{smallpmatrix}
0 & 0 & 1 \\
0 & 0 & 0 \\
0 & 0 & 0
\end{smallpmatrix}$.
Hence, for $n\in N$ and $X\in \sl_3(\C)$, an easy computation gives
\begin{equation*}
\scalebox{0.9}{$
\big(\psi_{\varpi_1+\varpi_2}(X)\big)(n)=X_{13}+aX_{23}+(ab-c)(X_{11}+aX_{21}+cX_{31})-b(X_{12}+aX_{22}+cX_{32})+c X_{33}
$}
\end{equation*}
which shows that 
\begin{equation*}
\Ima \psi_{\varpi_1+\varpi_2}\cap B=\{1,a,b,(ab-c),c,a(ab-c),bc,c(ab-c)\}.
\end{equation*}
From this, one can compute that the corresponding perfect basis $B_{\varpi_1+\varpi_2}:=\psi_{\varpi_1+\varpi_2}^{-1}(B)$ is 
\begin{equation*}
\scalebox{0.9}{$
\big\{
\begin{smallpmatrix}
0 & 0 & 1\\ 0 & 0 & 0\\ 0 & 0 & 0 
\end{smallpmatrix},
\begin{smallpmatrix}
0 & 0 & 0\\ 0 & 0 & 1\\ 0 & 0 & 0 
\end{smallpmatrix},
\begin{smallpmatrix}
0 & -1 & 0\\ 0 & 0 & 0\\ 0 & 0 & 0 
\end{smallpmatrix},
\begin{smallpmatrix}
2/3 & 0 & 0\\ 0 & -1/3 & 0\\ 0 & 0 & -1/3 
\end{smallpmatrix},
\begin{smallpmatrix}
-1/3 & 0 & 0\\ 0 & -1/3 & 0\\ 0 & 0 & 2/3
\end{smallpmatrix},
\begin{smallpmatrix}
0 & 0 & 0\\ 1 & 0 & 0\\ 0 & 0 & 0 
\end{smallpmatrix},
\begin{smallpmatrix}
0 & 0 & 0\\ 0 & 0 & 0\\ 0 & -1 & 0 
\end{smallpmatrix},
\begin{smallpmatrix}
0 & 0 & 0\\ 0 & 0 & 0\\ 1 & 0 & 0 
\end{smallpmatrix}
\big\}.
$}
\end{equation*}
\end{example}

\subsection{The \texorpdfstring{$\Dbar$}{D-bar} map}\label{subsec:Dbarmap}

We introduce two algebra homomorphisms
\begin{equation*}
\DbarPM:\C[N]\to \C(\ft)
\end{equation*}
first defined in \cite{baumann2021mirkovic} and which are known to distinguish biperfect bases (see \cite[Theorem~A.13]{baumann2021mirkovic}). They will play a key role in this paper, so we recall their definitions here. This section follows closely \cite[Section~8.1]{baumann2021mirkovic} (which extensively revolves around definitions and lemmas of \cite{bourbaki1968groupes}).

An element $x\in \fg$ is called a \textit{regular} element if
\begin{equation*}
\dim Z_\fg(x):=\dim \Ker \ad_x=\rk\fg
\end{equation*}
and a nilpotent element $x\in \fg$ is called a \textit{principal} nilpotent element if it is also regular.
Let $\fh^\reg$ denote the set of regular elements of $\fh$. Choose $e\in \fn$ a principal nilpotent element. For each $h\in \fh^\reg$, consider the subset $h+\fn\subset\fg$. The group $N$ acts on $h+\fn$ via the adjoint action and it is well-known that it acts simply transitively. Consequently, there exists a unique $n\in N$ such that $\Ad(n)h=h+e$. 

\begin{example}
Let $G=\SL_n$.The regular elements of $\fh$ are traceless diagonal matrices having distinct entries. Choose the principal nilpotent element $e$ having $1$s in entries $(i,i+1)$ for $i=1,\dots, n-1$ and zeros elsewhere. One can compute that the unique element $n\in N$ satisfying the above property has matrix entries $n_{ii}=1$, $n_{ij}=0$ for $i>j$ and 
\begin{equation*}
\textstyle n_{ij}=\prod_{k=i}^{j-1} \tfrac{1}{h_j-h_k}=\prod_{k=i}^{j-1} \tfrac{1}{(-\sum_{\ell=k}^{j-1} \alpha_\ell)(h)}
\end{equation*}
for $i<j$.
\end{example}

For any choice of principal nilpotent $e\in \fg$, the associations $h\mapsto n_h$ and $h\mapsto n_h^{-1}$ where $\Ad(n_h)h=h+e$ define two algebraic maps $\fh^\reg\to N$. Let 
\begin{align*}
\DbarP_e:\C[N]&\to \C[\fh^\reg] & \DbarM_e:\C[N]&\to \C[\fh^\reg] \\
f&\mapsto (h\mapsto f(n_h^{-1})) & f&\mapsto (h\mapsto f(n_h))
\end{align*}
be the pullback by the corresponding algebraic maps.\par 

The Weyl vector $2\rho^\vee=\sum_{\alpha>0}\alpha^\vee$ satisfies $2\rho^\vee(\alpha_i)=2$ for all $i\in I$. Hence, it determines up to scalar a principal nilpotent element $e=\sum_{i\in I} e_i$ which can be completed into an $\sl_2$ triple $(e,2\rho^\vee,f)$, where $e_i\in \fn$ is non-zero and of weight $\alpha_i$. For this fixed choice of principal nilpotent $e$, let $\DbarPM=\DbarPM_e$.

There is a useful rewriting of the maps $\DbarPM$ in terms of sequences of roots. For a sequence $\bfi=(i_1,\dots, i_d)$ of nodes of the Dynkin diagram, let $e_{\bfi}=e_{i_1}\dots e_{i_d}\in U(\fn)$ where $e_{i_j}$ is the normalization of the Chevalley generator fixed by the choice of principal nilpotent. Moreover, define $\beta_{k}^{\bfi}=\alpha_{i_1}+\dots+\alpha_{i_k}$ and put 
\begin{equation*}
\textstyle\DbarP_{\bfi}=\prod_{k=1}^{d}\tfrac{1}{\beta_{k}^{\bfi}}\hspace{2em}\text{and}\hspace{2em}\DbarM_{\bfi}=\prod_{k=0}^{d-1}\tfrac{1}{\beta_{k}^{\bfi}-\beta_{d}^{\bfi}}.
\end{equation*}	
As a consequence, one has:
\begin{proposition}[{\cite[Proposition~8.4]{baumann2021mirkovic}}]
There is an equality
\begin{equation}\label{eq:DbarVSDbarbfi}
\textstyle\DbarPM(f)=\sum_{\bfi} \langle e_{\bfi},f \rangle \,\DbarPM_{\bfi}
\end{equation}
as rational functions on $\fh$. 
\end{proposition}
 
\subsection{A weak change-of-basis}

Originally, the $\Dbar$ map appeared in the study of equivariant multiplicities of MV cycles which are irreducible components of subvarieties of the affine Grassmannian for the Langlands dual group (a precise definition of these varieties will be given in Section \ref{subsec:affineGrassSlices}). Under the geometric Satake isomorphism, (stable) MV cycles produce the basis \ref{b:MVbasis}. For $Z$ a stable MV cycle of weight $\nu\in Q_+$, denote by $b_{Z}\in \C[N]_{-\nu}$ the corresponding biperfect basis vector. As it will be explained in Section \ref{subsec:affineGrassSlices}, the association $Z\mapsto b_Z$ requires a choice of principal nilpotent $e\in \fn$ together with an isomorphism $\iota:\fh^\ast\to \fh$. Recall that $Z$ has two distinguished $T$-fixed points $L_0,L_{-\nu}\in \Gr$ which are respectively called the top and bottom fixed points. 

\begin{theorem}[{\cite[Theorem~1.4]{baumann2021mirkovic}}]
Fix an isomorphism $\iota:\fh^\ast\to \fh$ together with a choice of principal nilpotent $e\in \fn$.
Then, there is an equality
\begin{equation*}
\DbarM(b_Z)=\epsilon^T_{L_{-\nu}}([Z])
\end{equation*}
where the right-hand side denotes the equivariant multiplicity of $Z$ at its bottom fixed point. 
\end{theorem}

On the flip side, when $\fg$ is simply-laced, one can use finite-dimensional simple modules over KLR algebras to produce basis \ref{b:DCbasis}. Recall that the KLR algebra $R=\bigoplus_{\nu\in Q_+} R_{\nu}$ associated to $\fg$ admits a distinguished set of idempotents $\{e(\bfi)\}_{\bfi}$, where $\bfi$ runs over finite sequences of vertices of the Dynkin diagram. For $L$ a simple finite-dimensional $R_\nu$-module, denote by $d_L\in \C[N]_{-\nu}$ the corresponding biperfect basis vector. Modules over KLR algebras admit a \textit{character} which takes the form 
\begin{equation*}
\ch(L)=\textstyle \sum_{\bfi} \dim \big(e(\bfi)L\big) \; \bfi^{\rev}
\end{equation*}
seen as a formal linear combination of sequences. Here, $\bfi^{\rev}$ denotes the sequence obtained by reading $\bfi$ from right to left. By definition, $d_L\in \C[N]_{-\nu}$ is the unique function satisfying
\begin{equation*}
\langle e_{\bfi^\rev},d_L\rangle=\langle e_{i_d}\dots e_{i_1},d_L\rangle=\dim \big(e(\bfi)L\big)
\end{equation*}
for all sequences $\bfi=(i_1,\dots, i_d)$ such that $\sum \alpha_{i_j}=\nu$. \par 
To each simple finite-dimensional module $L$ over a cyclotomic KLR algebra, results of \cite{kamnitzer2019category} show how to produce a module $M$ of maximal Gelfand--Kirillov dimension which lies in category $\cO$ of a truncated shifted Yangian. In Section \ref{sec:asympchar}, we explain how one can make sense of the support of such modules and show that the top-dimensional support can be expressed (in the Borel--Moore homology) as a sum of MV cycles. For a fixed $L$, the above procedure gives, for each MV cycle $Z$, a non-negative integer $n_{L,Z}$ which corresponds to the multiplicity of $Z$ in the scheme-theoretic support of associated graded module of $M$. We are now ready to state one of the main results of this paper.

\begin{theorem}\label{thm:equalityDbar}
There is an equality
\begin{equation}\label{eq:DbardLequalsDbarbZ}
\textstyle \DbarM(d_L)=\sum_{Z} n_{L,Z} \DbarM(b_Z)
\end{equation}
of rational functions.
\end{theorem}
This theorem should be treated as evidence for the following:

\begin{conjecture}
The change-of-basis from basis \ref{b:DCbasis} to basis \ref{b:MVbasis} is given by the non-negative integers $n_{L,Z}$.
\end{conjecture}

The next example shows that Theorem \ref{thm:equalityDbar} encodes non-trivial combinatorics. It is greatly inspired by the exposition given in \cite[Section~3.2]{casbi2021equivariant}.

\begin{example}\label{ex:Nakada}

Fix a dominant weight $\lambda\in P_+$. For $\mu\in W\lambda$, let $w\in W$ be such that $w\lambda=\mu$. Furthermore, set $\nu=\lambda-\mu$, $d=\rht\nu$. Since $\mu$ is an extremal weight of $V(\lambda)$, the weight space $V(\lambda)_\mu$ has dimension $1$. Assume that $w$ is $\lambda$-minuscule in the sense of \cite{stembridge2001minuscule}.

On the one hand, it is shown in \cite{kleshchev2010homogeneous} that there is a unique finite-dimensional simple $R^{\lambda}_{\nu}$-module $L$ up to isomorphism, which is called a \textit{strongly homogeneous} module. Its character is given by the formula 
\begin{equation*}
\textstyle \ch(L)=\sum_{\bfi\in \redexp(w)} \bfi
\end{equation*}
where the sum runs over reduced expressions for $w$. 

On the other hand, as $\mu$ is an extremal weight, the MV cycle associated to the unique basis vector of $V(\lambda)_\mu$ is isomorphic to the affine space $\mathbb{A}^{d}$. By \cite[Section~4.1]{krylov2021almost}, the weights of the torus action at the fixed point $0\in \mathbb{A}^{d}$ are given by 
\begin{equation*}
\Phi_+^w:=\{\alpha\in \Phi_+\;;\; w(\alpha)<0\}.
\end{equation*}
Recall that when a variety is smooth at a $T$-fixed point, the equivariant multiplicity at that point is the inverse of the product of the weights of the torus action.

Following closely the arguments of Section \ref{subsec:attractingandrepelling}, it is possible to show that the coefficient $n_{L,Z}$ appearing in \eqref{eq:DbardLequalsDbarbZ} is equal to $1$. Applying Theorem \ref{thm:equalityDbar} yields
\begin{equation*}
\prod_{\alpha\in \Phi_+^w}\dfrac{1}{\alpha}=\sum_{(i_1,\dots, i_d)\in \redexp(w)} \dfrac{1}{\alpha_{i_1}(\alpha_{i_1}+\alpha_{i_2})\dots (\alpha_{i_1}+\dots+\alpha_{i_d})}
\end{equation*}
which is an equality known as \textit{Nakada's colored hook formula} \cite{nakada2008colored}. Evaluating the previous rational function at $\rho^\vee$ gives the equality
\begin{equation*}
\prod_{\alpha\in \Phi_+^w}\dfrac{1}{\rht(\alpha)}=\dfrac{|\redexp(w)|}{d!}
\end{equation*}
which is the classical \textit{Peterson--Proctor colored hook formula}. In type $A$, $|\redexp(w)|$ equals the number of standard Young tableaux of shape given by the heap of $w$ while the product $\prod_{\alpha\in \Phi_+^w} \rht(\alpha)$ is the product of the hook lengths.
\end{example}

\section{Asymptotic characters for filtered quantizations}\label{sec:asympchar}

The objective of this section is to develop a notion of \textit{asymptotic characters}. This will be done by building a framework of \textit{algebraic category $\cO$ for integral quantizations} which encompasses three notable families of examples, namely
\begin{enumerate}
\item truncated shifted Yangians of \cite{kamnitzer2014yangians},
\item quantizations of conical symplectic resolutions of \cite{braden2012quantizations,braden2014quantizations} and
\item Coulomb branch algebras of \cite{braverman2016coulomb}.
\end{enumerate}
This development is motivated by the fact that formalism of \cite{braden2012quantizations,braden2014quantizations} cannot be applied to truncated shifted Yangians whose shift is non-dominant or whose truncation coweight is not a sum of dominant minuscule coweights. However, as the results of the following sections will highlight, those cases are of crucial importance to study the interactions between bases \ref{b:DCbasis} and \ref{b:MVbasis}. One should note that the treatment presented here relies on some strong integrality and finiteness criteria (see hypotheses \ref{H:Bfinite} and \ref{H:Bintegral} for a more precise statement). Those criteria have an analogous statement in the classical BGG category $\cO$ setting: they amount to choosing a central character whose linkage class consists of integral weights. \par

The main result of this section is Theorem \ref{thm:achiandequivmult} which relates characters of modules over the quantization with the underlying geometry. Much of this section is only a careful recollection of known results. Applications to the above mentioned cases will be carried through in the following sections. The notation of this section will be slightly different from what is announced in Section \ref{sec:notation} and $\ft$ will denote the Lie algebra of a torus $T$.

\subsection{Filtered quantizations and the associated graded algebra}\label{subsec:filteredquant}

For a $\C$-algebra $A$, consider the following properties:
\begin{enumeratehypo}
\item\label{H:Afilt} The algebra $A$ is endowed with a separated, increasing $\Z$-filtration $F^\bullet$ which is compatible with the algebra structure.
\item\label{H:findimgenA} The algebra $A$ is finitely generated. 
\item\label{H:commR} The associated graded algebra 
\begin{equation*}
R:=\gr_F A=\textstyle\bigoplus_{n\in\Z} F^n A/F^{n-1}A
\end{equation*}
is commutative and finitely generated.
\item\label{H:filtExp} The filtration $F^\bullet$ admits an expansion, i.e.~there is an isomorphism of filtered $\C$-vector spaces $A\simeq \gr_F A$.
\end{enumeratehypo}
Under the hypotheses \ref{H:Afilt}-\ref{H:findimgenA}, the Gelfand--Kirillov (GK) dimension of $A$ and of finitely generated $A$-modules are defined. The reader is referred to \cite{krause2000growth} for definitions and properties of GK dimension. For $d\in \R_{\geq 0}$, one can consider the category $A\modu_{\leq d}$ which we call the \textit{GK subcategory of order $d$} and which is defined as the full subcategory of $A\modu$ whose objects are $A$-modules $M$ satisfying $\gkdim (M)\leq d$. Given $\scrC$ a full subcategory of $A\modu$, the notation $\scrC_{\leq d}$ will be used to denote the full subcategory of $\scrC$ whose objects lie in $A\modu_{\leq d}$.

The following result is an immediate consequence of \cite[Proposition~5.1]{krause2000growth}.
 
\begin{lemma}
If $\scrC$ is abelian, so is $\scrC_{\leq d}$.
\end{lemma}

For $M$ a finitely generated $A$-module, using the filtration $F^\bullet$ on $A$, one can define filtrations on $M$ which are \textit{good filtrations} (see \cite[Section~2.3.16]{chriss1997representation} for a more extensive discussion on good filtrations). For example, if $m_1,\dots, m_r\in M$ generate $M$, then the filtration $G^k(M):=\sum_{i=1}^r F^k(A) \,m_i$ where $k\in \Z$, defines a good filtration on $M$. For a chosen good filtration $G^\bullet$ on $M$, define $\gr_G M:= \bigoplus_{n\in\Z} G^n(M)/G^{n-1}(M)$. The vector space $\gr_G M$ admits a canonical structure of a finitely generated $R$-module. Though the association $M\mapsto \gr_G M$ is not functorial, the class $[\gr_G (M)]\in K_0(R\modu)$ is independent of the choice of filtration \cite[Corollary~2.3.19]{chriss1997representation}. We write the corresponding map by $[\gr(\trou)]$. Combining this result together with \cite[Lemma~6.5]{krause2000growth}, one deduces:

\begin{lemma}\label{lem:grleqd}
For $d\in \R_{\geq 0}$, the map $[\gr(\trou)]$ respects the filtrations induced by GK subcategories, i.e.~it restricts to a homomorphism of abelian groups
\begin{equation*}
[\gr(\trou)]:K_0(A\modu_{\leq d})\to K_0(R\modu_{\leq d}).
\end{equation*}
\end{lemma}

\begin{remark}\label{rem:boundedfiltGKdim}
When the filtration $F^\bullet$ on $A$ is bounded below, it is known that the GK dimension does not drop upon taking associated graded, see \cite{mcconnell1989gelfand}.
\end{remark}

Assume \ref{H:commR} holds and consider the affine scheme $X=\Spec R$. Let $V$ be an $R$-module and choose $S$ an irreducible component of $\supp V$, the support of $V$. By \cite[Lemma~5.9.2]{chriss1997representation}, there is a uniquely determined non-negative integer $n_{V,S}$, called the \textit{multiplicity of $S$ in $\supp V$}. Whenever the localization at the generic point $V_\fp$ is projective, this integer satisfies $n_{V,S}=\dim_{R_{\fp}} V_\fp$.

The \textit{support cycle} in Borel-Moore homology is the map defined by
\begin{equation}\label{eq:defsupport}
\begin{aligned}
[\supp_d(\trou)]:K_0(R\modu_{\leq d})&\to \HH_{2d}^{\BM}(X;\Z)\\
[V]&\mapsto\sum_{S\text{ $d$-dim, irred}} n_{V,S}[S]
\end{aligned}
\end{equation}
where the sum is over all $d$-dimensional irreducible components of $\supp V$. This construction is often referred to as the \textit{dévissage} principle. Here, $\HH:=\HH^{\BM}$ will always denote Borel--Moore homology. Note that $\HH^{\BM}(X,\Z)$ only depends on the $\C$-points of $X$. 

\begin{definition}[Characteristic cycle]\label{def:charcycle}
The composition 
\begin{equation*}
\CC_{\leq d}:=[\supp_d(\trou)]\circ [\gr(\trou)]:K_0(A\modu_{\leq d})\to \HH_{2d}(X;\Z)
\end{equation*}
is called the characteristic cycle map.
\end{definition}

\begin{remark}
An alternate version of Definition \ref{def:charcycle} appeared in a geometrical setting in the work of \cite{braden2012quantizations} following \cite{kashiwara2010deformation}. It is reasonable to ask whether this sheaf-theoretic characteristic cycle matches with the above defined characteristic cycle when applying global section and localization functors, see \cite[Section~2.6]{braden2014quantizations}. We do not pursue this avenue any further.
\end{remark}

Consequently, if a module $M$ has GK dimension equal to $d$, its characteristic cycle is a sum of $d$-dimensional irreducible components.

In the following sections, we will need a slight refinement of the above construction.

\begin{definition}[Support of a subcategory]\label{def:subcatSupportedOn}
Let $\scrC$ be an abelian full subcategory of $A\modu$ which is of finite length. Let $I\subset R$ be an ideal and set $Y:=\Spec R/I$. The subcategory $\scrC$ is said to be \textit{supported on $Y$} if for every simple object $L$ of $\scrC$, there exists a good filtration $G^\bullet$ of $L$ for which $I\cdot \gr_{G} L=0$.
\end{definition}

Let $\scrC$ be a subcategory which is supported on $Y=\Spec R/I$. For $M$ an object of $\scrC$ satisfying $[M]=\sum_{L} n_L[L]\in K_0(\scrC)$, let 
\begin{equation*}
\textstyle [\gr_I M]=\sum_L n_L[\gr_{G_I} L]\in K_0(R/I\modu)
\end{equation*}
where, for each $L$, the filtration $G_I^\bullet$ satisfies $I\cdot \gr_{G_I} L=0$. As previously, this linear map is well-defined (although not functorial). By definition of the support cycle, if $V$ is an $R/I$-module, then $[\supp_{ d}(V)]\in \HH_{2d}(Y;\Z)$. Define the $Y$-valued characteristic cycle by 
\begin{equation*}
\CC_{\leq d}^{Y}(M)=[\supp_d(\trou)]\circ [\gr_I(\trou)]:K_0(\scrC_{\leq d})\to \HH_{2d}(Y,\Z).
\end{equation*}

\begin{remark}
Definition \ref{def:subcatSupportedOn} has a sheaf-theoretical analogue, see \cite[Section~6.2]{braden2012quantizations}.
\end{remark}

\subsection{GK exact subcategories}\label{subsec:GKexact_subcats}

Let $A$ be a $\C$-algebra satisfying hypotheses \ref{H:Afilt} and \ref{H:findimgenA} given in Section \ref{subsec:filteredquant}.

\begin{definition}
A full subcategory $\scrC$ of $A\modu$ is \textit{GK exact} if for every short exact sequence
$\begin{tikzcd}[column sep=0.6em]
0 \arrow[r]& L \arrow[r]& M \arrow[r]& N \arrow[r]& 0
\end{tikzcd}$
in $\scrC$, one has $\gkdim(M)=\max\{\gkdim(L),\gkdim(N)\}$.
\end{definition}

\begin{remark}\label{rem:thmTauvel}
If $\bbk$ is a field of characteristic zero, $B$ is a filtered $\bbk$-algebra whose filtration is bounded below, and $\gr B$ is finitely generated and left noetherian, then the category $B\modu$ is GK exact, see \cite[Théorème~4.4]{tauvel1982dimension}.
\end{remark}

\begin{remark}
For finitely generated modules over a finitely generated commutative $\C$-algebra, GK dimension coincides with Krull dimension. Hence the category of finitely generated modules over a commutative $\C$-algebra is GK exact. 
\end{remark}

Let $\scrC$ be a GK exact full subcategory of $A\modu$. Denote by $\scrC_{<d}$ the full subcategory of objects having GK dimension strictly less than $d$. Then, $\scrC_{<d}$ is a Serre full subcategory of $\scrC_{\leq d}$ by assumption. Thus, one can consider the Serre quotient $\scrC_{=d}:=\scrC_{\leq d}/\scrC_{<d}$. By construction, it is an abelian category which satisfies $K_0(\scrC_{=d})\simeq K_0(\scrC_{\leq d})/K_0(\scrC_{<d})$.

Now, suppose furthermore that $A$ satisfies \ref{H:commR}.

\begin{lemma}\label{lem:CCpassingToTop}
Let $\scrC$ be a GK exact full subcategory of $A\modu$ supported on $Y\subset X$. The characteristic cycle passes to kernels, i.e.~there is a unique map  
\begin{equation*}
\CC_{=d}^Y:K_0(\scrC_{=d})\to \HH_{2d}(Y,\Z)
\end{equation*}
such that 
\begin{equation*}
\adjustbox{scale=0.9}{
\begin{tikzcd}[column sep=4em,row sep=-0.25em]
K_0(\scrC_{\leq d})\arrow[dd]\arrow[rd,"\CC_{\leq d}^Y"]& \\
&\HH_{2d}(Y,\Z)\\
K_0(\scrC_{=d})\arrow[ru,dashed,"\CC_{=d}^Y"']&
\end{tikzcd}}
\end{equation*}
commutes.
\end{lemma}

\begin{proof}
By Lemma \ref{lem:grleqd}, $[\gr_I(\trou)]$ induces a map on the level of Serre quotients
\begin{equation}\label{eq:commdiaggreqd}
\adjustbox{scale=0.9}{
\begin{tikzcd}[column sep=4em,row sep=1.25em]
K_0(\scrC_{\leq d})\arrow[r,"{[}\gr_I(\trou){]}"]\arrow[d]& K_0(R/I\modu_{\leq d})\arrow[d] \\
K_0(\scrC_{=d})\arrow[r,dashed,"{[}\gr_I(\trou){]}_{=d}"'] & K_0(R/I\modu_{=d})
\end{tikzcd}}
\end{equation}
which is denoted by ${[}\gr_I(\trou){]}_{=d}$. Since $K_0(R/I\modu_{<d})\subset\Ker\,[\supp_{d}(\trou)]$ by the dévissage principle, it follows that support cycles also induce a map on the level of Serre quotients 
\begin{equation*}
\adjustbox{scale=0.9}{
\begin{tikzcd}[column sep=4em,row sep=-0.25em]
K_0(R/I\modu_{\leq d})\arrow[dd]\arrow[rd,"{[}\supp_d(\trou){]}"]& \\
&\HH_{2d}(Y,\Z)\\
K_0(R/I\modu_{=d})\arrow[ru,dashed,"{[}\supp_d(\trou){]_{=d}}"']&
\end{tikzcd}}
\end{equation*}
which allows one to define $\CC_{=d}^Y$ as the composition ${[}\supp_d(\trou){]_{=d}}\circ{[}\gr_I(\trou){]}_{=d}$. 
\end{proof}

\subsection{Hamiltonian torus actions}\label{subsec:hamiltoniantorusaction}

Fix a torus $T\simeq (\C^\times)^r$ having Lie algebra $\ft=\Lie(T)$. The weight and coweight lattices of $T$ are denoted by $X^\ast(T):=\Hom_{\alg}(T,\C^\times)$ and $X_\ast(T)=\Hom_{\alg}(\C^\times,T)$ respectively.
The pairing between $X^{\ast}(T)$ and $X_\ast(T)$ is denoted by $\langle.,.\rangle$. By abuse of notation, $\langle.,.\rangle$ will also denote the pairing between $\ft^\ast$ and $\ft$.

Consider the category $\ft\Modu$ of (non-necessarily finitely generated) $\ft$-modules. An object $V$ of $\ft\Modu$ is called a \textit{weight module} if it admits a decomposition
\begin{equation*}
\textstyle V=\bigoplus_{\mu\in \ft^\ast} V(\mu)
\end{equation*}
where 
\begin{equation*}
V(\mu)=\big\{v\in V\;;\;\forall H\in \ft,\;\big(H-\langle \mu,H\rangle\big)^N v=0 \text{ for some }N\in\Z_{\geq 0}\big\}
\end{equation*}
and which satisfies $\dim V(\mu)<\infty$ for all $\mu\in \ft^\ast$. Denote by $\ft\wtModu$ the full subcategory of $\ft\Modu$ consisting of weight modules.

Let $\tstarVect$ denote the category of $\ft^\ast$-graded vector spaces which have finite-dimensional graded components. There is an obvious functor
\begin{equation*}
\ft\wtmodu\to\tstarVect
\end{equation*}
given by forgetting the $\ft$-action and keeping the $\ft^\ast$-grading coming from the $\ft$-action. The category $\tstarVect$ is endowed with a collection of shift functors $\{\Sigma_\tau\}_{\tau\in \ft^\ast}$, where each $\Sigma_\tau$ is an auto-equivalence of $\tstarVect$ defined on objects as $(\Sigma_\tau(V))(\mu)=V(\mu+\tau)$ and acting as the identity on morphisms. The weight lattice $X^\ast(T)$ gives rise to a full subcategory of $\tstarVect$ whose objects are $X^\ast(T)$-graded vector spaces which we denote by $\XastVect$. Since $T$-modules are semisimple, the data of an object of $\XastVect$ is the same as the data of a $T$-module having finite-dimensional isotypic components.

Along with the hypotheses \ref{H:Afilt} to \ref{H:filtExp} from Section \ref{subsec:filteredquant}, consider the following additional properties of the $\C$-algebra $A$: 
\begin{enumeratehypo}
\item\label{H:TactionA} There is a $T$-module structure on $A$ which is compatible with the algebra structure in the sense that for $\mu,\mu'\in X^\ast(T)$, if $a\in A(\mu)$ and $a'\in A(\mu')$, then $aa'\in A(\mu+\mu')$, where 
\begin{equation*}
A(\mu):=\{a\in A\;;\; ta=\mu(t)a\}.
\end{equation*}
\item\label{H:TactionHamiltonian} The $T$-action on $A$ is \textit{Hamiltonian}, i.e.~there is a morphism of Lie algebras $i_{\ft}:\ft\to A$ such that 
\begin{equation*}
[i_\ft(H),a]=\langle \mu, H\rangle a
\end{equation*} 
for all $a\in A(\mu)$ and $H\in \ft$, where the right-hand side is the action obtained by differentiating the action of $T$. 
\item\label{H:TpreservesF} The $T$-module structure on $A$ is compatible with the filtration $F^\bullet$ on $A$, meaning that the action of $T$ preserves the subspaces $F^nA$.
\end{enumeratehypo}

From now on, suppose that $A$ satisfies \ref{H:Afilt} to \ref{H:TpreservesF}.

Every $A$-module $M$ acquires a canonical $\ft$-module structure via pullback by the morphism $i_{\ft}:\ft\to A$. The module $M$ is called a \textit{weight module} if its pullback by $i_{\ft}$ is an object of $\ft\wtModu$. Let $A\wtmodu$ denote the full subcategory of $A\modu$ whose objects are weight modules for $\ft$ after pullback by $i_\ft$.

An $A$-module $M$ is called \textit{$T$-equivariant} if it admits a $T$-action which is compatible with the action of $A$ in the following sense. If $a\in A(\mu)$ and $m\in M(\nu)$, then $am\in M(\mu+\nu)$. The category of $T$-equivariant $A$-modules, denoted by $A\modu^T$, is the subcategory of $A\modu$ whose objects are $T$-equivariant $A$-modules and whose morphisms are also $T$-modules homomorphisms. One can define $A\wtmodu^T$ in a similar fashion.

\begin{remark}\label{rem:Tequivariantandextensions}
From an $A$-module $M$, the weights of the $\ft$-module $i_{\ft}^\ast M$ might not be integral. Furthermore, the category $A\wtmodu^T$ is not necessarily closed under extensions.
\end{remark}

Given an object $M$ of $A\wtmodu$, consider the $\ft^\ast$-graded vector space associated to $M$. If there is a weight $\tau\in \ft^\ast$ such that $\Sigma_\tau(M)$ is an object of $\XastVect$, then $\Sigma_\tau(M)$ acquires a canonical $T$-module structure providing it with the structure of a $T$-equivariant weight $A$-module.

From \ref{H:TpreservesF}, it follows that if $M\in \ob(A\wtmodu^T)$, one can choose a good filtration $G^\bullet$ for which $\gr_G M$ admits a canonical $T$-module structure.

\subsection{Category \texorpdfstring{$\cO$}{O}}

Let $A$ be a $\C$-algebra satisfying the hypotheses \ref{H:Afilt} to \ref{H:TpreservesF} and fix a cocharacter $\rho^\vee\in X_\ast(T)$. This endows $A$ with a $\Z$-grading by declaring 
\begin{equation*}
A(n)=\bigoplus_{\substack{n\in \Z\\ \langle\mu,\rho^\vee\rangle=n}} A(\mu),
\end{equation*}
where the distinguished element $\delta=(i_\ft\circ\dif\rho^\vee)(1)\in A$ makes the grading inner. Let $A_{\geq 0} =\bigoplus_{n\geq 0} A(n)$ and define $A_{>0}$, $A_{\leq 0}$, $A_{<0}$ analogously. In much the same way, the cocharacter $\rho^\vee$ makes $R$ into a $\Z$-graded algebra.\par 

\begin{remark}
Unfortunately, the notation $A(0)$ now has two different meanings. To remedy this, the trivial character for $T$ will be denote by $\bfZero\in X^\ast(T)$. 
\end{remark}

\begin{definition}[{\cite[Definition~3.10]{braden2014quantizations}}]\label{def:catO}
The category $\cO$ associated to $\rho^\vee\in X_\ast(T)$
is the full subcategory of $A\modu$ consisting of $A$-modules where $A_{\geq0}$ acts locally finitely. 
\end{definition}

As in \cite[Section~4.2]{ginzburg2003primitive}, the $B$-algebra of $A$ is defined by 
\begin{equation*}
B(A)=\dfrac{A(0)}{\sum_{n>0}A(-n)A(n)}.
\end{equation*}
The pullback by the quotient map $A_{\geq 0} \to B(A)$ gives a functor $B(A)\fmodu\to A_{\geq 0}\modu$ which allows one to define a parabolic induction functor 
\begin{equation}\label{eq:parabolicInd}
\Delta(\trou):=A\otimes_{A_{\geq 0}} (\trou):B(A)\fmodu\to A\modu.
\end{equation}
The following result is discussed in \cite[Section~4]{losev2017categories}.

\begin{lemma}\label{lem:simpleObjectsO}
For $V\in \ob(B(A)\fmodu)$, one has $A\otimes_{A_{\geq 0}} V\in\ob(\cO)$. Furthermore, if $S$ is a finite-dimensional simple $B(A)$-module, then $\Delta(S)$ has a unique maximal submodule, and consequently, a unique simple quotient. Conversely, any simple object of $\cO$ can be obtained as the unique simple quotient of a module $\Delta(S)$ for some finite-dimensional simple $B(A)$-module $S$. 
\end{lemma}

When $S$ is a simple $B(A)$-module, $\Delta(S)$ is called a \textit{generalized Verma module}.

\begin{remark}\label{rem:fingeneratedVermas}
Fix a $\C$-basis $\{s_1,\dots,s_r\}$ of $S$. The elements $\{1_A\otimes s_1,\dots, 1_A\otimes s_r\}$ form a generating set of $\Delta(S)$ as an $A_{\leq 0}$-module. Hence, generalized Verma modules are finitely generated $A_{\leq 0}$-modules. By the previous lemma, this implies that simple objects of $\cO$ are finitely generated $A_{\leq 0}$-modules.
\end{remark}

\begin{remark}\label{rem:univVerma}
When $S$ is $1$-dimensional, the generalized Verma module associated to $S$ satisfies the same universal property as the usual Lie theoretic Verma modules. Indeed, let $\lambda:B(A)\to \C$ be the algebra map making $\C\simeq \C_\lambda$ into a $1$-dimensional module isomorphic to $S$. By adjunction, the Verma module induced from $\C_\lambda$ satisfies
\begin{align*}
\Hom_{A}(A\otimes_{A_{\geq 0}} \C_\lambda,M)&\simeq \Hom_{A_{\geq 0}}(\C_\lambda,{}_{A_{\geq 0}}M)\\
&=\{m\in M\;;\; a m=0 \text{ for }a\in A_{>0}\text{ and }am=\lambda(a) m\text{ for a}\in A(0)\}
\end{align*}
as expected.
\end{remark}

\begin{remark}
If one chooses $\rho^\vee=\mathbf{0}$, then $B(A)=A$. In this case, Lemma \ref{lem:simpleObjectsO} is vacuously true since $\cO$ coincides with the category of finite-dimensional $A$-modules. 
\end{remark}

As this section will highlight, an important property which helps control category $\cO$ is the following.

\begin{definition}
The algebra $A$ is called $B$-finite if $\dim B(A)<\infty$. 
\end{definition}

For further referencing, consider the following hypothesis:
\begin{enumeratehypo}
\item\label{H:Bfinite} The algebra $A$ is $B$-finite.
\end{enumeratehypo}
In turn out that checking for $B$-finiteness can be done on the level of the associated graded algebra. In fact, hypotheses \ref{H:filtExp} and \ref{H:TpreservesF} imply that there is a vector space isomorphism $\gr B(A)\simeq B(A)$. The following lemma provides a sufficient criterion implying $B$-finiteness.

\begin{lemma}\label{lem:surjBrgA}
There is a surjection of $\C$-algebras $B(\gr A)\to \gr B(A)$. 
\end{lemma}

\begin{proof}
First, notice that the filtration $F^\bullet$ does in fact induce a filtration on $B(A)$, since for each $\mu\in X^\ast(T)$ the subspace $A(\mu)$ acquires a filtration defined by $F^k(A(\mu))=F^k A\cap A(\mu)$. Thus, if one denotes
\begin{equation*}
\textstyle M(k)=F^k A\cap \big(\sum_{n>0} A(-n)A(n)\big), 
\end{equation*}
one can write $F^k B(A)=(F^k A\cap A(0))/M(k)$. \par 

The following commutative diagram of $\C$-vector spaces can easily be deduced from the above set-up:
\begin{equation*}
\adjustbox{scale=0.8}{
\begin{tikzcd}[row sep=1em, column sep=1em]
F^{k-1} A(0) \arrow[r]\arrow[d]& F^k A(0)\arrow[r]\arrow[d] & \dfrac{F^k A(0)}{F^{k-1} A(0)}\arrow[r]\arrow[d,dashed,"\phi_k"] & 0\\
\dfrac{F^{k-1} A(0)}{M(k-1)}\arrow[r]& \dfrac{F^{k} A(0)}{M(k)}\arrow[r] & \dfrac{F^{k} A(0)}{M(k)}\big/\dfrac{F^{k-1} A(0)}{M(k-1)}\arrow[r] & 0
\end{tikzcd}}
\end{equation*}
Furthermore, the morphism $\phi_k$ is surjective since the second column morphism is. Hence, this yields a surjective morphism of $\C$-vector spaces 
\begin{equation*}
\phi:=\oplus_k \phi_k : \gr A(0)\to \gr B(A)=\bigoplus_{k\in \Z}  \dfrac{F^{k} A(0)}{M(k)}\big/\dfrac{F^{k-1} A(0)}{M(k-1)}.
\end{equation*}
One can check that $\phi$ indeed is a morphism of $\C$-algebras.\par 

Next, we show that $\sum_{n>0} (\gr A)(-n)\, (\gr A)(n)\subset \Ker \phi$. This follows the fact that if
\begin{equation*}
x\in F^k \big(\textstyle\sum_{n>0} (\gr A)(-n)\, (\gr A)(n)\big)
\end{equation*}
then, there exists $a\in F^i A\cap A(-n)$ and $b\in F^jA\cap A(n)$ such that $i+j=k$ and $x=ab \mod F^{k-1} A$. Since $ab\in M(k)$, it follows using the construction of $\phi_k$ in diagram above that $x\in \Ker \phi_k$. \par 

This concludes the proof, as $\phi$ induces a surjective morphism $B(\gr A)\to \gr B(A)$. 
\end{proof}

\begin{remark}
The above lemma is hinted at in the proof of \cite[Proposition~5.1]{braden2014quantizations} (which treats only algebras $A$ which have a non-negative filtration). Though it might be obvious to experts, we included it here for the sake of completeness.
\end{remark}

\begin{remark}
In general, the map $B(\gr A)\to \gr B(A)$ is not injective. It would be interesting to determine sufficient conditions guaranteeing that $B(\gr A)\simeq \gr B(A)$.
\end{remark}

Assuming \ref{H:Bfinite}, let $S_1,\dots, S_k$ denote a full list of isoclasses of simple $B(A)$-modules. For $i=1,\dots, k$, let $\Delta_i=\Delta(S_i)$ be the generalized Verma module associated with $S_i$ and $L_i:=\tp \Delta_i$ its maximal simple quotient.

\begin{lemma}[{\cite[Lemma~4.1]{losev2017categories}, \cite[Section~2]{ginzburg2003category}}]\label{lem:LosevCatO}
Assuming \ref{H:Bfinite}, the following properties hold: 
\begin{enumerate}
\item The simples $L_1,\dots, L_k$ form a complete set of non-isomorphic simple objects of $\cO$. 
\item Objects of $\cO$ have finite length.
\item Category $\cO$ is the full subcategory of $A\modu$ whose objects are such that 
\begin{itemize}
\item $A_{>0}$ acts locally nilpotently, 
\item $\delta$ acts locally finitely and 
\item generalized eigenspaces under $\delta$ are finite-dimensional.
\end{itemize}
\end{enumerate}
Moreover, objects of $\cO$ are weight modules.
\end{lemma}

\begin{remark}\label{rem:fingeneratedAleq}
The above result implies that objects in $\cO$ are finitely generated as $A_{\leq 0}$-modules. As discussed previously in Remark \ref{rem:fingeneratedVermas}, this is true for simple objects of $\cO$. The fact that objects of $\cO$ have finite length completes the proof.
\end{remark}

Recalling that the $\Z$-grading on $A$ is inner, one can check that $i_{\ft}(\ft)\subset A(0)$ since $\ft$ is an abelian Lie algebra. Thus, $B(A)$ inherits the structure of an algebra with a Hamiltonian torus action.

\begin{definition}
The algebra $A$ is called $B$-integral if there exists a $\tau\in \ft^\ast$ such that, for any finite-dimensional simple $B(A)$-module $S$, the vector space $\Sigma_\tau(S)$ is $X^\ast(T)$-graded. 
\end{definition}

Again, for further referencing, consider the following hypothesis:
\begin{enumeratehypo}
\item\label{H:Bintegral} The algebra $A$ is $B$-integral.
\end{enumeratehypo}
The following observation is immediate from Lemma \ref{lem:LosevCatO}.

\begin{lemma}\label{lem:TequivO}
If $A$ is satisfies \ref{H:Bfinite} and \ref{H:Bintegral}, then all generalized Verma modules can be endowed with the structure of a $T$-equivariant $A$-module and consequently, the same holds for all simple objects of $\cO$. 
\end{lemma}

The following example illustrates the criteria of $B$-finiteness and $B$-integrality.

\begin{example}\label{ex:BalgLie}
Let $\fg$ be a finite-dimensional simple Lie algebra over $\C$ with Cartan $\fh$ and Borel $\fb$. Fix $\lambda\in\fh^\ast$ and consider the $\C$-algebra 
\begin{equation*}
A:=U(\fg)/\langle z-\chi_\lambda(z)\;;\;z\in Z(\fg)\rangle
\end{equation*}
where $\chi_\lambda(z)$ is the scalar appearing in the action of $z$ on the highest weight vector of the Verma module $\Delta(\lambda)$ (often referred to as the central character associated to $\lambda$). Let $\rho^\vee=\tfrac{1}{2}\sum_{\alpha\in \Phi_+}\alpha\in\fh$ be the Weyl vector. This algebra is filtered using the \textit{Kazhdan filtration} of \cite[Section~3.2]{brundan2008highest} (i.e.~we use the good grading induced by $\rho^\vee$) and it satisfies the hypotheses \ref{H:Afilt} to \ref{H:TpreservesF}. Thus, it gives rise to a category $\cO$ in the sense of Definition \ref{def:catO}. The triangular decomposition of $U(\fg)$ induces a surjection $\Sym\fh \to B(A)$. Using the (twisted) Harish-Chandra homomorphism, we see that
\begin{equation*}
B(A)\simeq \Sym\fh/\langle \psi(f)-f(\lambda+\rho)\;;\; f\in\Sym\fh^W \rangle,
\end{equation*}
where $\psi(f)(\mu)=f(\mu+\rho)$ for all $\mu\in \fh^\ast$. The Chevalley--Shephard--Todd theorem states that $\Sym\fh$ is free of finite rank as a $\Sym\fh^W$-module. Thus, if $f_1,\dots, f_k\in \Sym \fh$ are a $\Sym\fh^W$-generating set of $\Sym\fh$, their images in $B(A)$ form a $\C$-basis. Hence, $B(A)$ is finite-dimensional and $A$ satisfies \ref{H:Bfinite}. It is a classical result of Harish-Chandra (see \cite[Theorem~1.10]{humphreys2008representations} for example) that all central characters $\chi:Z(\fg)\to \C$ are of the form $\chi=\chi_\mu$ for some $\mu\in\fh^\ast$ and that $\chi_\mu=\chi_\lambda$ implies that $\mu=w\cdot\lambda$ for some $w\in W$ (the Weyl group action here is the dot action). Hence, all simple $B(A)$-modules are morphisms of the form $w\cdot \lambda:B(A)\to \C$. This shows that the weights of the simple $B(A)$-modules are $E=\{w\cdot \lambda\;;\; w\in W\}$ (known as the linkage class of $\lambda$). Therefore, $A$ satisfies \ref{H:Bintegral} if there is a weight $\tau\in \fh^\ast$ such that $\tau+E$ consists of only integral weights. This condition implies that for all $w\in W$, the weight $\lambda-w\lambda\in \fh^\ast$ is integral. In particular, choosing $\lambda$ integral suffices. Notice that category $\cO$ for $A$ according to Definition \ref{def:catO} is not a block of the classical BGG category $\cO$ for $\fg$ as we allow for generalized weight spaces. However, for $\lambda\in\fh^\ast$ regular, \cite[Théorème~1]{soergel1986equivalences} shows that it is equivalent to the block $\cO_\lambda$ of the BGG category $\cO$.
 
Consider the above example in the case of $\fg=\sl_2$. It is well-known that $Z(\sl_2)=\C[\Omega]$, where $\Omega=\tfrac{1}{2}h^2+ef+fe$ is the Casimir element. By writting $\Omega=\tfrac{1}{2}h^2+h+2fe$, we see that 
\begin{equation*}
B(A)\simeq \C[h]/\langle (h-\lambda)(h-(-\lambda-2))\rangle
\end{equation*}
since $\chi_\lambda(\Omega)=\tfrac{1}{2}(\lambda+1)^2-\tfrac{1}{2}$. The algebra $\gr B(A)$ is isomorphic to the coinvariant algebra $\C[x]/(x^2)$. Notice that $\lambda\neq -1$ if and only if $B(A)$ is a semisimple algebra. If $\lambda=0$, its linkage class is $\{0,-2\}$ and category $\cO$ has two indecomposable projective objects, one of which can be presented as 
\begin{equation}\label{eq:presentationprojectivesl2}
\dfrac{U(\sl_2)}{U(\sl_2)\spa_\C\{e^2, he, h-fe, \Omega\}}.
\end{equation}
Note that for the above module, $h$ does not act semisimply.
\end{example}

\subsection{Attracting and repelling sets}\label{subsec:attractingandrepelling}

Let $A$ be a $\C$-algebra satisfying \ref{H:Afilt} to \ref{H:Bintegral}, together with a fixed $\rho^\vee$. The scheme $X=\Spec R$ admits a $\C^\times$-action given via $\rho^\vee$.\par 

Following \cite{drinfeld2013algebraic}, let $R_-:= R/\langle R_{>0}\rangle$ and define the scheme-theoretic repelling set of $X$ under the above $\C^\times$-action as $X_-:=\Spec R_-$. Equivalently, $R_-$ is the largest quotient of $R$ which is $\Z_{\leq 0}$-graded. The scheme-theoretic attracting set is defined analogously.\par 

\begin{example}\label{ex:repellingsl2}
Building on Example \ref{ex:BalgLie} in the case of $\fg=\sl_2$ and $\lambda=0$, it is well-known that $A$ quantizes the algebra of functions on the nilpotent cone 
\begin{equation*}
X=\{\big(\begin{smallmatrix}a & b\\ c & -a\end{smallmatrix}\big)\in\mathbf{M}_{2\times2}(\C)\;;\; a^2+bc=0\}.
\end{equation*}
Write $\alpha,\beta,\gamma\in \C[X]=R$ for the corresponding coordinate functions. The choice of Weyl vector 
\begin{equation*}
t\mapsto \begin{smallpmatrix}
t & 0 \\ 0 & t^{-1}
\end{smallpmatrix}
\end{equation*}
endows $X$ with a $\C^\times$-action via matrix conjugation. This endows $R$ with $\Z$-grading on $R$ such that $\deg \beta=-2$, $\deg \gamma=2$ and $\deg \alpha=0$. It follows that $R_-\simeq \C[\alpha,\beta]/(\alpha^2)$. 
\end{example}

Given a finite-dimensional simple $B(A)$-module $S$, let $\Delta(S)$  be the corresponding generalized Verma module. Choose $s_1,\dots, s_r\in S$ a set of generators of $S$ as a $B(A)$-module. The elements $1\otimes s_1,\dots, 1\otimes s_r\in \Delta(S)$ generate $\Delta(S)$ as an $A$-module. The \textit{highest weight filtration} associated to the choice of generators $\{s_1,\dots, s_r\}$ is the filtration $G_{\hw}^\bullet$ defined by
\begin{equation*}
G_{\hw}^k\big(\Delta(S)\big):=F^{k}A\cdot \spa_{\C} \{1\otimes s_1,\dots, 1\otimes s_r\}.
\end{equation*}
It is a good filtration on $\Delta(S)$. By Lemma \ref{lem:LosevCatO}, $G_{\hw}^\bullet$ induces a filtration on its simple top. We call any filtration which is the image of a filtration $G_{\hw}^\bullet$ a highest weight filtration. Again, Lemma \ref{lem:LosevCatO} guarantees that every simple object of $\cO$ has such a filtration. 

Recall that a filtration $G^\bullet$ on an $A$-module $M$ is a \textit{separated filtration} if $\bigcap_{n\in Z} G^nM=0$. To guarantee that objects of $\cO$ have a separated filtration, we impose a compatibility between the filtration on $A$ and the cocharacter $\rho^\vee$, namely:

\begin{enumeratehypo}
\item\label{H:triangular} If $\mu$ is a weight of $F^k A$ with $k<0$, then $\langle\mu,\rho^\vee\rangle>0$
\end{enumeratehypo}
Notice that such a condition is trivially satisfied when $F^\bullet$ is bounded below. This hypothesis is also closely related to the notion of a generalized \textit{triangular decomposition} appearing in \cite{ginzburg2003primitive}. From now on, suppose that $A$ satisfies \ref{H:triangular}.

\begin{lemma}
Any highest weight filtration on an object of $\cO$ is separated.
\end{lemma}

\begin{proof}
It suffices to prove the statement for generalized Verma modules. Let $M$ be such an object and fix a basis $\{s_1,\dots,s_r\}$ of the corresponding finite-dimensional simple $B(A)$-module. Let $G_{\hw}^\bullet$ be the corresponding highest weight filtration. If $k\in \Z_{<0}$, one has 
\begin{equation*}
G_{\hw}^k(M)=F^k(A)\cdot \spa_\C\{1\otimes s_1,\dots,1\otimes s_r\}.
\end{equation*}
By \ref{H:triangular}, every weight vector of $F^k(A)$ belongs to $A_{>0}$. By construction, $A_{>0}$ annihilates elements of the form $1\otimes s_i$ for $i=1,\dots, r$ and thus, $G_{\hw}^k(M)=0$. The statement follows.
\end{proof}

\begin{example}
Consider the $\C$-algebra $A=\C[t,t^{-1}]$ filtered according to the degree. 
Consider the (trivial) Hamiltonian $\C^\times$-action, where $t^n$ has weight $1$ for all $n\in \Z$.  
Let $M=\C$ be the $A$-module on which the variable $t$ acts trivially. 
For every cocharacter $\rho^\vee\in X_\ast(\C^\times)\simeq \Z$, this object belongs to the corresponding category $\cO$.
The module $M$ has a good filtration $G^\bullet$ defined by $G^kM:=F^kA\cdot 1\simeq \C$. This filtration is not a separated filtration.
This example shows that \ref{H:triangular} is necessary in order to ensure that objects of $\cO$ have a separated filtration. 
\end{example}

It would be interesting to determine if techniques similar to \cite[Section~5.1]{finkelberg2018comultiplication} can always be applied to modify the filtration on $A$ so that \ref{H:triangular} is automatically satisfied.

\begin{lemma}\label{lem:repellerAnnihilator}
If $L$ is a simple object of $\cO$, then for any highest weight filtration $G_{\hw}^\bullet$ on $L$, $\langle R_{>0} \rangle\subset\ann(\gr_{G_{\hw}} L)$ and $\gr_{G_{\hw}} L$ is a module over $R_-$. 
\end{lemma}

\begin{proof}

By Lemma \ref{lem:LosevCatO}, let $S$ be a finite-dimensional simple $B(A)$-module such that $L=\tp \Delta(S)$. Let $G^\bullet=G_{\hw}^\bullet$ be a highest weight filtration on $\Delta(S)$. The result will follow if we show that $R_{>0}$ acts by zero on $\gr_{G} \Delta(S)$.  

Choose a set of generators $s_1,\dots, s_r\in S$ for the $B(A)$-module $S$ giving rise to $G^\bullet$. Choose $r\in R_{>0}$ with $r\in F^{m_1} A/F^{m_1-1} A$ for some $m_1\in\Z$ and let $a\in (F^{m_1} A)\cap(A_{>0})$ be such that $r=a\mod F^{m_1-1} A$. In a similar way, choose $x\in \gr_G \Delta(S)$ with $x\in G^{m_2}(\Delta(S))/G^{m_2-1}(\Delta(S))$ for some $m_2\in \Z$ and let $v\in G^{m_2}(\Delta(S))$ be such that $x=v\mod G^{m_2-1}(\Delta(S))$. As mentioned previously, there exists $a_1,\dots, a_r\in A$ such that $v=\sum_{\ell=1}^r a_\ell \otimes s_\ell$. Notice that for each $\ell$, $a_\ell\in F^{m_2}A$ which implies that $a\,a_\ell-a_\ell\, a\in F^{m_1+m_2-1}A$ by \ref{H:commR} and thus, 
\begin{equation*}
(a\,a_\ell-a_\ell\, a)\cdot 1\otimes s_\ell\in G^{m_1+m_2-1}\Delta(S).
\end{equation*}
It follows that
\begin{align*}
a \cdot a_\ell \otimes s_\ell \mod G^{m_1+m_2-1}(\Delta(S))&=(a \, a_\ell) \cdot 1\otimes s_\ell \mod G^{m_1+m_2-1}(\Delta(S))\\
&=(a_\ell \, a) \cdot 1\otimes s_\ell \mod G^{m_1+m_2-1}(\Delta(S))\\
&=a_\ell \cdot a\otimes s_\ell \mod G^{m_1+m_2-1}(\Delta(S))\\
&=0\mod G^{m_1+m_2-1}(\Delta(S))
\end{align*}
as $a\in A_{>0}$ implies $a\otimes s_\ell=1\otimes a s_\ell=0$. This proves that $rx=0$, which shows that $\langle R_{>0}\rangle$ acts as zero on $\gr_G \Delta(S)$ as required.
\end{proof}

\begin{example}
In Example \ref{ex:repellingsl2}, both $\Delta(0)$ and $\Delta(-2)$ can be presented as
\begin{equation*}
\Delta(\lambda)=\dfrac{U(\sl_2)}{U(\sl_2)\spa_{\C}\{h-\lambda,e\}}
\end{equation*}
and upon passing to associated graded (using the highest weight filtration coming from the generator $1\in U(\sl_2)$), both yield the $R$-module $\C[\beta]$ but with two different $T$-actions. Clearly, $R_{>0}$ acts by zero on $\C[\beta]$. Notice that for the projective module $P(-2)$ (see \eqref{eq:presentationprojectivesl2}), when passing to associated graded (using the good filtration again coming from $1\in U(\sl_2)$), $R_{>0}$ does not act by zero on it. 
\end{example}

\begin{lemma}\label{lem:dimmugr}
Assuming \ref{H:triangular} holds, if $M\in \ob(\cO)$ and $G^\bullet$ is a highest weight filtration, then $\dim M(\mu)=\dim (\gr_G M)(\mu)$. 
\end{lemma}

\begin{proof}
For $i\in \Z$, consider the following short exact sequence of vector spaces: 
\begin{equation*}
\begin{tikzcd}[column sep =1em]
0 \arrow[r] & G^{i-1}(M)\arrow[r] &  G^{i}(M)\arrow[r] & G^{i}(M)/G^{i-1}(M)\arrow[r] & 0
\end{tikzcd}
\end{equation*}
Since taking the $\mu$ (generalized) weight-space is exact, we have a short exact sequence
\begin{equation*}
\begin{tikzcd}[column sep =1em]
0 \arrow[r] & G^{i-1}(M)(\mu)\arrow[r] & G^{i}(M)(\mu)\arrow[r] & (G^{i}(M)/G^{i-1}(M))(\mu)\arrow[r] & 0
\end{tikzcd}
\end{equation*}
which implies that $\dim (G^{i}(M)/G^{i-1}(M))(\mu)=\dim G^{i}(M)(\mu)-\dim G^{i-1}(M)(\mu)$. Since the filtration on $M$ is separated by assumption and since $\dim M(\mu)<\infty$, there exists $n_1,n_2\in \Z$ with $n_1\leq n_2$ such that $\dim G^{n_1}(M)(\mu)=0$ and $\dim G^{n_2}(M)(\mu)=\dim M(\mu)$. It follows that 
\begin{align*}
\dim (\gr M)(\mu)&=\textstyle\sum_{i=n_1}^{n_2} \dim (G^{i}(M)/G^{i-1}(M))(\mu)\\
&=\textstyle\sum_{i=n_1}^{n_2}\big(\dim G^{i}(M)(\mu)-\dim G^{i-1}(M)(\mu)\big)=\dim M(\mu)
\end{align*}
where we used the fact that the sum is telescoping.
\end{proof}

\subsection{Fixed points}\label{subsec:fixedpoints}

Building on work of Section \ref{subsec:attractingandrepelling} and following again \cite{drinfeld2013algebraic}, define the scheme-theoretic fixed points of $X$ as $X_0=\Spec B(R)$. Notice that $B(R)$ is the largest graded quotient of $R$ which is concentrated in degree $0$ and $B(R)$ is identified with the $0$-graded component of $R_-$. Hence, there is an inclusion of algebras $i_-:B(R)\to R_-$. 

Consider the two final hypotheses:

\begin{enumeratehypo}
\item\label{H:BRlocal} The scheme $X_0$ has a unique closed point, or equivalently, the algebra $B(R)$ is local. 
\item\label{H:Bcomm} The $B$-algebra of $A$ is commutative.
\end{enumeratehypo}

\begin{remark}\label{rem:hyposExpPlusLocal}
From Lemma \ref{lem:surjBrgA}, it follows that \ref{H:filtExp} and \ref{H:BRlocal} imply \ref{H:Bfinite}.  
\end{remark}

\begin{remark}
When the algebra $A$ arises as the quantization of a conical symplectic resolution $Y\to X$, the fact that $B(A)$ is commutative can often be deduced from the relationship between $\dim B(A)$ and $T$-fixed points of $Y$, see \cite{losev2021localization}.
\end{remark}

\begin{remark}
When $B(A)$ is commutative, all generalized Verma modules satisfy the universal property of Remark \ref{rem:univVerma}. 
\end{remark}

Assuming \ref{H:BRlocal}, let $\fm$ denote the unique maximal ideal of $B(R)$. It follows that $B(R)$ is a finite-dimensional, $\C$-algebra with nilradical $\fm$ satisfying  $B(R)/\fm\simeq \C$. Consequently, if $\fj\subset R_-$ denotes the ideal generated by $i_-(\fm)$, the $0$-graded component of $\tilde{R}_-:=R_-/\fj$ is isomorphic to $\C$. Furthermore, since $\fm$ is nilpotent, so is $\fj$. \par 

\begin{example}
Recall the setting from Example \ref{ex:repellingsl2}. The $B$-algebra of $A$ is commutative and the $B$-algebra of $R$ is isomorphic to $\C[\alpha]/(\alpha^2)$. It is a local algebra and the algebra $\tilde{R}_-$ is isomorphic to $\C[\beta]$. 
\end{example}

The $\C$-algebra $\tilde{R}_-$ can be identified with the tensor product $B(R)/\fm\otimes_{B(R)} R_-$ via the morphism $i_-$. Hence, the \textit{reduced scheme-theoretic repelling set} $\tilde{X}_-:=\Spec \tilde{R}_-$ is the fibered product of $X_-$ and $\{\ast\}$ over $X_0$. Note that the reduced scheme-theoretic repelling set might not be a reduced scheme. Moreover, on the level of $\C$-points, $\tilde{X}_-$ and $X_-$ are isomorphic as affine algebraic varieties.

\begin{remark}\label{rem:smoothFixedPoint}
When $X$ has a unique $\C^\times$-fixed point which is smooth, both $X_0$ and $X_{-}$ are reduced by \cite[Proposition~1.4.20]{drinfeld2013algebraic}
\end{remark}

\begin{lemma}\label{lem:BalgebraAnnihilator}
If $L$ is a simple object of $\cO$, then for any highest weight filtration $G^\bullet_{\hw}$ on $L$, $\fm\subset\ann(\gr_{G_{\hw}} L)$ and $\gr_{G_{\hw}} L$ is a module over $\tilde{R}_-$. 
\end{lemma}

\begin{proof}
Recall that we assume hypothesis \ref{H:Bcomm}. Consequently, every finite-dimensional simple $B(A)$-module is $1$-dimensional. By Lemma \ref{lem:LosevCatO}, one can choose a finite-dimensional $B(A)$-module $S$ such that $L=\tp \Delta(S)$. Consider $s\in S$ such that $s\neq 0$. Then, $1\otimes s\in \Delta(S)$ generates $\Delta(S)$ as an $A$-module. Moreover, let $G^\bullet$ be the corresponding highest weight filtration on $\Delta(S)$. Choose a minimal integer $m\in \Z$ such that $1\otimes s \in G^m \Delta(S)$ and $1\otimes s \not\in G^{m-1} \Delta(S)$ and define $x:=1\otimes s \mod G^{m-1} \Delta(S)$. For any $y\in \gr\Delta(S)$, there exists $r'\in R$ such that $y=r'x$. Thus, since $R$ is commutative, the statement of the lemma is equivalent to showing that if $r\in \fm\subset B(R)$, then $r$ acts as zero on $x\in \gr_G\Delta(S)$.

Let $a\in F^n A\cap A(0)$ be such that $a\mod F^{n-1} A=r$. Then, 
\begin{align*}
r\cdot x&= 1\otimes as \mod G^{m+n-1} \Delta(S)\\
&=\lambda(a)\,1\otimes s\mod G^{m+n-1} \Delta(S)
\end{align*}
for some $\lambda(a)\in \C$. If $n>0$, then $1\otimes s\in G^{m}\Delta(S)\subset G^{m+n-1}\Delta(S)$ and it follows that $rx=0$. Assume that $\lambda(a)\neq 0$. If $n<0$, then, $\tfrac{1}{\lambda(a)} a \cdot 1\otimes s= 1\otimes s\in G^{m+n}\Delta(S)\subset G^{m-1}\Delta(S)$ which contradicts our assumption on $m$. Hence, assume $n=0$. Let $N\in \Z_{\geq 0}$ be such that $r^N=0$. Then, 
\begin{equation*}
0=r^N x=\lambda(a)^N\,1\otimes s\mod G^{m-1} \Delta(S)
\end{equation*}
which shows that $1\otimes s\in G^{m-1} \Delta(S)$, contradicting once again the hypothesis on $m$. It follows that $\lambda(a)=0$, from which the result follows.
\end{proof}

\begin{corollary}\label{cor:catOsupport}
Category $\cO$ is supported on $\tilde{X}_-$.
\end{corollary}

\begin{proof}
This follows immediately from the Lemmas \ref{lem:LosevCatO}, \ref{lem:repellerAnnihilator} and \ref{lem:BalgebraAnnihilator}.
\end{proof}

\subsection{Characters}\label{subsec:characters}

Let $\charRing$ be the $\C$-vector space spanned by infinite sums of the form $\sum_{\mu\in\ft^\ast} d_\mu e^{\mu}$ with $d_\mu\in\C$. For two elements $p=\sum_{\mu\in\ft^\ast} p_\mu e^\mu,q=\sum_{\mu\in\ft^\ast} q_\mu e^\mu\in \charRing$, if $q$ has finite support, the product of $p$ and $q$ is defined as 
\begin{equation*}
p\cdot q = \sum_{\nu\in\ft^\ast} \big(\sum_{\nu=\mu_1+\mu_2}p_{\mu_1}q_{\mu_2}\big)\,e^\nu.
\end{equation*}
Furthermore, if $p,q\in \charRing$ are both supported in a union of cones of the form $\lambda+Q_-\subset\ft^\ast$, their product is also defined by the same formula as above.  

For an object $V$ of $\tstarVect$, define the character (or the $\ft^\ast$-graded dimension) of $V$ by 
\begin{equation*}
\chi(V)=\textstyle\sum_{\mu\in \ft^\ast} \dim M(\mu)\, e^\mu\in \charRing.
\end{equation*}

The following result is well-known:

\begin{lemma}\label{lem:chigrouphom}
The character map $\chi:K_0(\tstarVect)\to \charRing$ is a group homomorphism.
\end{lemma}

It is easy to notice that shift functors of Section \ref{subsec:hamiltoniantorusaction} are well-behaved with respect to characters. In fact, one can deduce the rather trivial identity
\begin{equation}\label{eq:shifts_and_chars}
\chi(\Sigma_\tau(V))=e^\tau\cdot \chi(V).
\end{equation}
By composition with the forgetful functor, the character map extends to $K_0(A\wtmodu)$ as well as $K_0(R\wtmodu^T)$. The resulting map will still be denoted by $\chi$.

As noted previously in Section \ref{subsec:hamiltoniantorusaction}, upon passing to associated graded, $T$-equivariant $A$-modules give rise to $T$-equivariant $R$-modules. Recall that the algebra $A$ satisfies \ref{H:Bintegral} by assumption. As an immediate corolllary of Lemma \ref{lem:dimmugr}, one has:

\begin{lemma}\label{lem:grchar}
For $M\in \ob(\cO)$, one has $\chi(M)=\chi(\gr(M))$. 
\end{lemma}

On $K_0(A\wtmodu^T)$, as $\chi$ commutes with $[\gr(\trou)]$, it follows that the composition
\begin{equation*}
[\gr(\trou)]\circ [\Sigma_\tau]:K_0(\cO)\to K_0(\tilde{R}_-\wtmodu^T)
\end{equation*}
is injective whenever $\chi:K_0(\cO)\to \charRing$ is.

Consider the following two additional properties on $\cO$: 
\begin{enumeratehypoO}
\item\label{HO:GKexact} Category $\cO$ is Gelfand-Kirillov exact.
\item\label{HO:GKassograded} The GK dimension of objects of $\cO$ does not change upon taking associated graded (see Remark \ref{rem:boundedfiltGKdim}).  
\end{enumeratehypoO}

\begin{remark}
Assuming \ref{HO:GKassograded}, the GK dimension of objects of $\cO$ are integers.
\end{remark}

\begin{remark}
When $A$ arises as the quantization of a conical symplectic singularity, the filtration $F^\bullet$ on $A$ is bounded below. Thus, both properties \ref{HO:GKexact} and \ref{HO:GKassograded} hold. 
\end{remark}

\begin{corollary}\label{cor:CCinj}
If $\chi$ is injective on $K_0(\cO)$ and the two additional hypotheses \ref{HO:GKexact} and \ref{HO:GKassograded} hold, then for $d=\dim_\C \tilde{X}_-$, $\CC_{=d}^{\tilde{X}_-}$ is injective. 
\end{corollary}

\begin{proof}
If $\chi$ is injective on $K_0(\cO)$, then $[\gr(\trou)]\circ [\Sigma_\tau]$ is injective. The commutative diagram \eqref{eq:commdiaggreqd} then guarantees then injectivity of $[\gr(\trou)]_{=d}$ by \ref{HO:GKassograded}. Finally, as noted in \cite[p.290]{chriss1997representation}, the top dimensional support map
\begin{equation*}
[\supp(\trou)]_{=d}:K_0(\tilde{R}_-\modu_{=d})\to \HH_{2d}(\tilde{X}_-,\Z)=\HH_{2d}(X_-;\Z)
\end{equation*}
is an isomorphism. Consequently, $\CC_{=d}^{\tilde{X}_-}$ is the composition of two injective maps, which is injective.
\end{proof}

\subsection{Equivariant Hilbert polynomials and equivariant multiplicities}\label{subsec:equivHilbert}

Consider the following function which is defined in terms of characters of objects of $\cO$.

\begin{definition}[Asymptotic character]\label{def:asymptoticCharacterDef}
For $M\in\ob(\cO_{\leq d})$, the asymptotic character of $M$, denoted by $\achi_d(M)$, is the $\C$-valued function on $\ft$ defined pointwise as 
\begin{equation*}
\achi_d(M)(h)=\textstyle\lim_{n\to\infty}\tfrac{1}{n^{d}}\sum_{\mu} \dim M(\mu)\; e^{\langle\mu,h\rangle/n}
\end{equation*}
where $h\in \ft$. 
\end{definition}

Up to this point, it is not clear whether or not taking the asymptotic character gives a well-defined function.

Recall that the $\C$-algebra $\tilde{R}_-$ has a $T$-action, it is finitely generated and it satisfies $\tilde{R}_-(0)=\C$. Fix a set of generators $r_1,\dots,r_N\in \tilde{R}_-$ such that $r_i$ is a weight vector of weight $-\beta_i\in X^\ast(T)$ for the $T$-action. By construction, the generators are such that $\langle-\beta_i,\rho^\vee\rangle < 0$. Consider the algebra $C=\C[x_1,\dots, x_N]$, together with the $T$-action given by specifying that $x_i$ has weight $-\beta_i$. By construction, there is a surjective morphism of algebras $\pi:C\to \tilde{R}_-$ which is $T$-equivariant.

The objective of this section is to show that $\achi_d$ is indeed a rational function which computes the equivariant multiplicity of the class $\CC^{\tilde{X}_-}_{\leq d}(M)\in \HH_{2d}(\tilde{X}_-;\Z)=\HH_{2d}(X_-;\Z)$. This will be done by building upon work from previous sections. The following diagram gives an idea of the proof:
\begin{center}
{
\usetikzlibrary{decorations.pathmorphing}
\tikzcdset{arrows=rightsquigarrow}
\begin{tikzcd}
\cO_{\leq d} \arrow[r,"\Sigma_\tau"] & A\wtmodu_{\leq d}^T \arrow[r,"\gr_{G_{\hw}}(\trou)"] &  \tilde{R}_-\wtmodu^T_{\leq d}\arrow[r,"\text{pullback}"] & C\wtmodu^T_{\leq d}
\end{tikzcd}
}
\end{center}
We will show that the computation of $\achi_d$ on the rightmost node of the above diagram computes  equivariant multiplicity of the support of $T$-equivariant $C$-modules and that applying the above composition on the level of Grothentieck groups does not change the value of $\achi_d$, hence proving the result. \par

The reader should keep in mind the fact that the manipulations are happening on the level of Grothendieck groups, as the proof will rely on two ``non-functorial" associations. The first one is the fact that we will show the result on simple objects of $\cO$ and extend it linearly, thus not distinguishing an object of $\cO$ from the direct sum of its simple factors. The second one is the fact that taking associated graded is not functorial. 

The data of the $\C$-algebra $C$ together with its torus action appearing in Section \ref{subsec:fixedpoints} matches the framework needed to define \textit{equivariant Hilbert polynomials} of finitely generated $T$-equivariant $C$-modules following \cite[Section~6.6]{chriss1997representation} (which are also called \textit{multidegrees} in commutative algebra, see \cite{miller2005combinatorial}).\par 

We recall the construction of the equivariant Hilbert polynomial. By \cite[Proposition~6.6.6]{chriss1997representation}, the character of a finitely generated $T$-equivariant $C$-module $V$ can be written as 
\begin{equation*}\label{eq:characterChrissGinz}
\chi(V)=\dfrac{\sum_\lambda n_\lambda e^\lambda}{\prod_{i=1}^N (1-e^{-\beta_i})}
\end{equation*}
where the numerator has finitely many non-zero terms. Let $\chi_{C,V}:=\sum_\lambda n_\lambda e^\lambda$, seen as a function on $T$. Using the exponential map $\exp:\ft\to T$, this produces a function on $\ft$ via pullback which can be Taylor expanded around $0\in\ft$ as 
\begin{equation*}
\chi_{C,V}=\textstyle\sum_{\lambda\in X^\ast(T)}\sum_{m=0}^\infty n_\lambda\frac{\lambda^m}{m!} =\sum_{i=0}^{+\infty} p^{(i)}
\end{equation*}
where 
\begin{equation*}
p^{(i)}=\textstyle\tfrac{1}{i!}\sum_{\lambda} n_\lambda \lambda^i
\end{equation*}
is a homogeneous polynomial function on $\ft$ of degree $i$. The equivariant Hilbert polynomial of the pair $(C,V)$ is the first non-vanishing $p^{(i)}$ and it is denoted $p_{C,V}$. \par

\begin{example}
If one considers the torus $T=(\C^\times)^n$, one has $X^{\ast}(T)\simeq \Z^n$. Hence, if $t=(t_1,\dots, t_n)\in T$ and $\lambda=(\lambda_1,\dots, \lambda_n)\in X^\ast(T)$, the function $e^{\lambda}$ on $T$ is the evaluation at $\lambda$, namely
\begin{equation*}
e^{\lambda}(t)=t_1^{\lambda_1}\dots t_n^{\lambda_n}.
\end{equation*}
For $h\in\ft$, pulling back to $\ft\simeq \C^n$ by $\exp:\ft\to T$ gives 
\begin{equation*}
e^{\lambda}(h_1,\dots,h_n)=e^{\langle \lambda,h\rangle}=e^{\lambda_1h_1+\dots+\lambda_nh_n}.
\end{equation*}
\end{example}

Let $d=\dim\supp V$. It is proven in \cite[Theorem~6.6.12]{chriss1997representation} that $p_{C,V}$ is a homogeneous polynomial of degree $N-d$ and that
\begin{equation}\label{eq:equivHilbertpolymultiplicities}
p_{C,V}=\sum_{i=1}^r n_i \,p_{C,Z_i}
\end{equation}
where $Z_1,\dots, Z_r$ are the irreducible components of dimension $d$ appearing in $\supp V$ respectively with multiplicities $n_1,\dots, n_r$ and where the polynomial $p_{C,Z_i}$ is defined as $p_{C,\C[Z_i]}$

The following calculus lemma shows that one can recover $p_{C,V}$ from a well-chosen limit.

\begin{lemma}\label{lem:asympcharcomm}
Let $V$ be a finitely generated $T$-equivariant $C$-module with $d=\dim\supp V$. Then, the following equality holds pointwise
\begin{equation*}
\lim_{n\to\infty} \dfrac{1}{n^d} \sum_{\mu} \dim V(\mu)\, e^{\mu/n} =\dfrac{p_{C,V}}{\prod_{i=1}^N\beta_i},
\end{equation*}
seen as a rational function on $\ft$ (and consequently, the limit is well-defined).
\end{lemma}

\begin{proof}
Choose $h\in \ft$ such that $\langle \beta_i,h\rangle\neq 0$ for all $i$. It follows by homogeneity of the $p^{(i)}$'s that
\begin{align*}
\lim_{n\to\infty} n^{N-d}\cdot\chi_{C,V}(h/n)&=\lim_{n\to\infty} n^{N-d} \textstyle\sum_{i=0}^{+\infty} p^{(i+N-d)}(h/n)\\
&=\lim_{n\to\infty} \textstyle\sum_{i=0}^{+\infty} \dfrac{1}{n^i} p^{(i+N-d)}(h)=p_{C,V}(h).
\end{align*}
Using similar calculus manipulations, one has
\begin{align*}
\lim_{n\to\infty} n^N\cdot \prod_{i=1}^N(1-e^{-\langle\beta_i,h\rangle/n})=\lim_{n\to\infty}n^N\cdot \prod_{i=1}^N\Big(\dfrac{\langle\beta_i,h\rangle}{n}+O(\tfrac{1}{n^2})\Big)=\prod_{i=1}^N \langle\beta_i,h\rangle.
\end{align*}
Combining those two elementary manipulations with \eqref{eq:characterChrissGinz}, it follows that 
\begin{equation*}
\lim_{n\to\infty} \dfrac{1}{n^d} \sum_{\mu} \dim V(\mu)\, e^{\langle\mu,h\rangle/n} =\lim_{n\to\infty} \dfrac{n^{N-d}}{n^N} \dfrac{\chi_{C,V}(h/n)}{\prod_{i=1}^N(1-e^{-\langle\beta_i,h\rangle/n})}=\dfrac{p_{C,V}(h)}{\prod_{i=1}^N\langle\beta_i,h\rangle}
\end{equation*}
which completes the proof.
\end{proof}

This proves the following result on the non-commutative level.

\begin{proposition}\label{prop:kernelOfAsymptoticChar}
Let $M\in \ob(\cO_{\leq d})$. The expression $\achi_d(M)$ is a $\C$-valued rational function on $\ft$. It is a homomorphism of abelian groups
\begin{equation*}
\achi_d:K_0(\cO_{\leq d})\to \C(\ft)
\end{equation*}
whose kernel contains $K_0(\cO_{<d})$.
\end{proposition}

\begin{proof}
The fact that $\achi_d$ is a morphism comes from Lemma \ref{lem:chigrouphom} and from linearity of limits. Hence, we only need to prove that $\achi_d$ is a rational function for simple objects of $\cO_{\leq d}$. 

Let $L\in \ob(\cO_{\leq d})$ be a simple object. Then, by assumption, there is a $\tau\in \ft^\ast$ such that $L_\tau:=\Sigma_\tau(L)$ is a $T$-equivariant $A$-module. Consequently, upon choosing a highest weight filtration, Corollary \ref{cor:catOsupport} shows that $V:=\gr_{G_{\hw}} L_\tau$ is an object of $\tilde{R}_-\wtmodu^T_{\leq d}$ which is finitely generated. Via pullback, $V$ becomes a $C$-module. It is easy to see that pullback to $C$ preserves characters, as the map $C\to \tilde{R}_-$ is $T$-equivariant. Applying Lemma \ref{lem:asympcharcomm} shows that 
\begin{equation}\label{eq:asymptoticCharValue}
\lim_{n\to\infty} \tfrac{1}{n^d} \textstyle\sum_{\mu} \dim V(\mu)\, e^{\mu/n}=\begin{cases}
p_{C,V}/(\prod_{i=1}^N \beta_i) & \text{if }\dim \supp V=d\\
0 & \text{if }\dim \supp V<d
\end{cases}.
\end{equation} 
As mentioned in Section \ref{subsec:characters}, there is an equality of characters $\chi(L)=e^{-\tau}\cdot \chi(L_\tau)$ and applying Lemma \ref{lem:grchar}, one has $\chi(L)=e^{-\tau} \cdot\chi(V)$. Since $e^{-\tau/n}\to 1$ pointwise when $n\to \infty$, it follows that 
\begin{equation*}
\achi_d(L)=\lim_{n\to\infty} \tfrac{1}{n^d} \textstyle \sum_{\mu} \dim V(\mu)\, e^{\mu/n}
\end{equation*}
which, combined with \eqref{eq:asymptoticCharValue} shows the desired result. 
\end{proof}

Following \cite[Section~9.2]{baumann2021mirkovic} as well as \cite[Chapter~8]{miller2005combinatorial} and \cite{brion1997equivariant}, we explain how equivariant Hilbert polynomials are related to equivariant multiplicities. For $C$ as above, consider the affine space $W=\Spec C$. For $V$ a $T$-equivariant finitely generated $C$-module, as in \eqref{eq:defsupport}, there is a well-defined class $[\supp V]\in \HH_{2d}^T(W)$ where $d=\dim \supp V$. Moreover, if $Z_1,\dots, Z_r$ denotes the set of $d$-dimensional irreducible components of $\supp V$, we have
\begin{equation*}
[\supp V]=\textstyle\sum_{i=1}^r n_i [Z_i]
\end{equation*}
in $\HH_{2d}^T(W)$, where $n_i$ is the multiplicity of $\supp V$ along $Z_i$. The ring $\HH^\bullet_T(\mathrm{pt})\simeq \Sym\ft^\ast$ acts on $\HH_\bullet^T(W)$ in a natural way. Let $k=\Frac(\HH^\bullet_T(\mathrm{pt}))$ denote the fraction field of $\HH^\bullet_T(\mathrm{pt})$. By \cite[Corollary~4.2]{brion1997equivariant}, $k\otimes_{\HH^\bullet_T(\mathrm{pt})}\HH_\bullet^T(W)$ is a $1$-dimensional $k$-vector space. Thus, for $x\in\HH_\bullet^T(W)$, there exists a unique $\epsilon^T(x)\in k$ such that
\begin{equation*}
\epsilon^T(x)\cdot [\{0\}]=[x]. 
\end{equation*}
The linear map $\epsilon^T:\HH_\bullet^T(W)\to k$ is called the equivariant multiplicity. For $V$ a $T$-equivariant $C$-module, recall that $p_{C,V}$ defined above is a polynomial function on $\ft$ of degree $N-d$, hence $p_{C,V}\in \Sym^{N-d}\ft^\ast\simeq\HH^{2N-2d}_T(\mathrm{pt})$.

\begin{lemma}[{\cite[Section~5]{brion1997equivariant}}]\label{lem:actionHilbertPoly}
The equality
\begin{equation}\label{eq:actionPCM}
[\supp V]=p_{C,V}\cdot [W]
\end{equation}
holds in $k\otimes_{\HH^\bullet_T(\mathrm{pt})}\HH_{2d}^T(W)$.
\end{lemma}

As a special case of Lemma \ref{lem:actionHilbertPoly}, if one picks the $1$-dimensional $C$-module $\C$ with a trivial $T$-action, its support is $\{0\}\subset W$. This module has trivial character, which shows that its equivariant Hilbert polynomial is $\prod_{i=1}^N \beta_i$. Applying \eqref{eq:actionPCM} to this module yields 
\begin{equation*}
[\{0\}]=\textstyle\big(\prod_{i=1}^N \beta_i\big)[W].
\end{equation*}
Thus, by definition, one has 
\begin{equation*}
[\supp V]=p_{C,V}\,[W]=\dfrac{p_{C,V}}{\big(\prod_{i=1}^N \beta_i\big)} [\{0\}]
\end{equation*}
which shows that the equivariant multiplicity of $\supp V$ is given by
\begin{equation*}
\epsilon^T(\supp V)=\dfrac{p_{C,V}}{\big(\prod_{i=1}^N \beta_i\big)}.
\end{equation*}
Note that when $U\subset W$ is a smooth variety, \cite[Theorem~4.2(iii)]{brion1997equivariant} shows that $1/\epsilon^T([U])$ is equal to the product of the weights of the $T$-action on the tangent space of $U$ at $0$.

Coming back to our original setting one has: 
\begin{theorem}\label{thm:achiandequivmult}
Let $Z_1,\dots,Z_r$ be the $d$-dimensional irreducible components of $\tilde{X}_-$. For $M\in\ob(\cO_{\leq d})$ such that $\CC_{\leq d}^{\tilde{X}_-}([M])=\sum n_i[Z_i]\in \HH_{2d}(\tilde{X}_-;\Z)$, one has
\begin{equation*}
\achi_{d}(M)=\textstyle\sum_{i=1}^r n_i \epsilon^T(Z_i)=\epsilon^T\big(\CC_{\leq d}([M])\big).
\end{equation*}
\end{theorem}

\begin{proof}
From \eqref{eq:equivHilbertpolymultiplicities} and \eqref{eq:asymptoticCharValue}, one deduces that 
\begin{equation*}
\achi_d(M)=\textstyle\sum_{i=1}^r n_i \,\dfrac{p_{C,\C[Z_i]}}{\prod_{i=1}^N \beta_i}
\end{equation*}
which is precisely the definition of $\sum_{i=1}^r n_i \epsilon^T(Z_i)$. The result follows by $\Z$-linearity.
\end{proof}

\begin{example}
Again in the case of Example \ref{ex:repellingsl2}, the character of the projective module $P(-2)$ is 
\begin{equation*}
\chi(P(-2))=1+2e^{-2}+2e^{-4}+2e^{-6}+\dots =1+2e^{-2}(1-e^{-2})^{-1}.
\end{equation*}
The support of its associated graded module contains two copies of the line defined by $Z=\{a=c=0\}\subset X$ since $P(-2)$ has two composition factors isomorphic to $\Delta(-2)$. As predicted by Theorem \ref{thm:achiandequivmult}, one has
\begin{equation*}
\achi_{1}(P(-2))=\lim_{n\to \infty} \tfrac{1}{n}\cdot \Big(1+\tfrac{2e^{-2/n}}{1-e^{-2/n}}\Big)= 2\cdot \tfrac{1}{2}=2\cdot \epsilon^T(Z).
\end{equation*}
\end{example}

\begin{remark}\label{rem:is_Bintegral_necessary}
As the referee pointed out, it is not necessary to assume \ref{H:Bintegral} in order to prove Theorem \ref{thm:achiandequivmult}. By the assumptions \ref{H:TactionA} and \ref{H:TactionHamiltonian}, for every simple object $L$ of $\cO$, there exists a weight $\tau\in \ft^\ast$ such that $\Sigma_\tau(L)$ is a $T$-equivariant $A$-module. While the choice of $\tau$ may depend on $L$, taking the asymptotic character is insensible to this dependency (see \eqref{eq:shifts_and_chars}). However, for results such as Corollary \ref{cor:CCinj}, \ref{H:Bintegral} is necessary (in Example \ref{ex:BalgLie}, taking $\fg=\sl_3$ along with a non-integral central character exhibits the necessity). 
\end{remark}

\section{Truncated shifted Yangians}\label{sec:Yangian}

This section defines \textit{shifted Yangians} of \cite{kamnitzer2014yangians} as well as some remarkable quotients called \textit{truncated shifted Yangians}. A brief review of the geometry of the affine Grassmannian together with the affine schemes quantized by the previously mentioned algebras is given in Section \ref{subsec:affineGrassSlices}. Doing so, we show how the formalism developed in Section \ref{sec:asympchar} can be applied to study the representation theory of truncated shifted Yangians, proving that asymptotic characters of modules in category $\cO$ of truncated shifted Yangians compute the equivariant multiplicity of the sum of MV cycles of their characteristic cycles.

The main results appear in Section \ref{subsec:asymptoticcharYangians}. Many definitions and lemmas make use of sets of parameters which are defined in Appendix \ref{sec:setofparameters}. 

\subsection{Definitions and properties}
The following definition is taken from \cite[Definition~3.1]{kamnitzer2022hamiltonian}. It first appeared in \cite{kamnitzer2014yangians}, generalizing \cite{brundan2006shifted}.

\begin{definition}[{\cite[Definition~3.1]{kamnitzer2022hamiltonian}}]\label{def:shiftedYangian}
Given a coweight $\mu\in P^\vee$, the shifted Yangian $Y_\mu(\fg)$ is the $\C$-algebra with generators $e_{i,r}, f_{i,r}, h_{i,p}$ for $ i\in I$, $r\in \mathbb{Z}_{>0} $, $p\in\mathbb{Z}$, and subject to the following relations: 
\begin{align}
[h_{i,p},h_{j,q}]&=0\label{eq:commHH}\\
[e_{i,r},f_{j,s}]&=2\delta_{ij}h_{i,r+s-1}\label{eq:EFH}\\
[h_{i,p+1},e_{j,r}]-[h_{i,p},e_{j,r+1}]&=		(\alpha_i,\alpha_j)		(h_{i,p}e_{j,r}+e_{j,r}h_{i,p})\label{eq:HEroot}\\
[h_{i,p+1},f_{j,r}]-[h_{i,p},f_{j,r+1}]&=		-(\alpha_i,\alpha_j)	(h_{i,p}f_{j,r}+f_{j,r}h_{i,p})\label{eq:HFroot}\\
[e_{i,r+1},e_{j,s}]-[e_{i,r},e_{j,r+1}]&=		(\alpha_i,\alpha_j)		(e_{i,r}e_{j,s}+e_{j,s}e_{i,r})\\
[f_{i,r+1},f_{j,s}]-[f_{i,r},f_{j,r+1}]&=		-(\alpha_i,\alpha_j)	(f_{i,r}f_{j,s}+f_{j,s}f_{i,r})\\
\textstyle\sum_{\sigma\in \fS_m} [e_{i,r_{\sigma(1)}},[e_{i,r_{\sigma(2)}},[\dots[e_{i,r_{\sigma(m)}},e_{j,s}]\dots ]]]&=0\label{eq:symE}\\
\textstyle\sum_{\sigma\in \fS_m} [f_{i,r_{\sigma(1)}},[f_{i,r_{\sigma(2)}},[\dots[f_{i,r_{\sigma(m)}},f_{j,s}]\dots ]]]&=0\label{eq:symF}\\
h_{i,p}&=0 \text{ when }p<-\langle\alpha_i,\mu\rangle\\
h_{i,-\langle\alpha_i,\mu\rangle}&=1\label{eq:Hmueq1}
\end{align}
where, in \eqref{eq:symE} and \eqref{eq:symF}, $i\neq j$, $m=1-a_{ij}$ and $r_1,\dots,r_m\in\Z_{\geq0}$. 
\end{definition}

\begin{remark}
The above normalization of the relations corresponds to setting $\hbar=2$ in the original definition of \cite{kamnitzer2014yangians}. It is the same convention as \cite{kamnitzer2019category}. 
\end{remark}

\begin{remark}
Setting $\mu=0$, one recovers the usual Yangian $Y(\fg)$ presented according to Drinfeld's second realization (with $\hbar=2$).
\end{remark}

The current notation will be used to bookkeep generators of $Y_\mu$. We write
\begin{align*}
h_i(u)&=\textstyle\sum_{p\in\mathbb{Z}} h_{i,p}u^{-p} = u^{\langle\alpha_i,\mu\rangle} + h_{i,-\langle\alpha_i,\mu\rangle+1} u^{\langle\alpha_i,\mu\rangle-1}+h_{i,-\langle\alpha_i,\mu\rangle+2} u^{\langle\alpha_i,\mu\rangle-2}+\cdots\\
e_i(u)&=\textstyle\sum_{r\in\Z_{>0}} e_{i,r} u^{-r}=e_{i,1}u^{-1}+e_{i,2}u^{-2}+\dots\\
f_i(u)&=\textstyle\sum_{r\in\Z_{>0}} f_{i,r} u^{-r}=f_{i,1}u^{-1}+f_{i,2}u^{-2}+\dots
\end{align*}
seen as elements of $Y_\mu\parparl u^{-1}\parparr$, the algebra of formal Laurent series in the variable $u^{-1}$ with coefficients in $Y_\mu$.

One can check that $Y_\mu$ is a finitely generated algebra. It is generated by the finite-dimensional subspace
\begin{equation}\label{eq:gensetYangian}
\scalebox{0.9}{$
\spa_\C\Big\{e_{i,1},f_{i,1},h_{i,p}\;;\;i\in I, p\in\{-\langle\alpha_i,\mu\rangle+1,-\langle\alpha_i,\mu\rangle+2\}\cup\{q\;;\;	-\langle\alpha_i,\mu\rangle+3\leq q\leq -1\}\Big\}
$}
\end{equation}
since, if one defines
\begin{equation*}
s_{i,1}=h_{i,-\langle\alpha_i,\mu\rangle+2}-\tfrac{1}{2} (h_{i,-\langle\alpha_i,\mu\rangle+1})^2,
\end{equation*}
then
\begin{equation*}
[\tfrac{1}{4d_i} s_{i,1},e_{i,r}]=e_{i,r+1}\hspace{1em}\text{and}\hspace{1em}[\tfrac{1}{4d_i} s_{i,1},f_{i,r}]=-f_{i,r+1}.
\end{equation*}
All other generators can then be obtained by applying \eqref{eq:EFH}.\par

Fix a coweight $\mu\in P^\vee$. Let $\lambda=\sum_{i\in I} \lambda_i\varpi_i^\vee\in P^\vee_+$ be such that $\lambda-\mu=\sum_{i\in I} m_i\alpha_i^\vee\in Q_+^\vee$. Furthermore, let $\bR=(R_i)_{i\in I}$ be a set of parameters of level $\lambda$ (see Appendix \ref{sec:setofparameters} for the definition of sets of parameters\footnote{While the appendix requires $\fg$ to be simply-laced, the definitions of sets of parameters and GT-weights also make sense for non-simply-laced types.}). Following \cite[Lemma~2.1]{gerasimov2005class}, define the currents $a_i(u)\in Y_\mu\parparl u^{-1}\parparr$ by the implicit equation
\begin{equation*}
h_i(u)= p_i(u)\cdot \dfrac{\prod_{j\neq i}\prod_{r=1}^{-a_{ji}} (u-d_ia_{ij}-2rd_j)^{m_j}}{u^{m_i}(u-2d_i)^{m_i}}\cdot\dfrac{\prod_{j\neq i}\prod_{r=1}^{-a_{ji}}  a_j(u-d_{i} a_{ij}-2rd_j)}{a_i(u)a_i(u-2d_i)}
\end{equation*}
where $p_i(u)=\prod_{c\in R_i}(u-c)$. These currents are well-defined and they can be written as 
\begin{equation*}
\textstyle a_i(u)=\sum_{p\in\Z} a_{i,p} u^{-p} = 1+a_{i,1} u^{-1}+a_{i,2}u^{-2}+\dots
\end{equation*}
and from \cite[Equation~(2.8)]{gerasimov2005class}, one can check that 
\begin{equation}\label{eq:GKLOmodes}
\textstyle h_{i,p}=\sum_{j\in I} -a_{ji} a_{j,p+\langle\alpha_j,\mu\rangle}+\Big(\text{terms involving $a_{j,q}$ for }0\leq q\leq p+\langle\alpha_j,\mu\rangle-1\Big).
\end{equation}
Since the Cartan matrix of $\fg$ is invertible, the previous equation shows that the subalgebra generated by the modes $\{h_{i,p}\}_{p\in \Z, i\in I}$ is equal to the subalgebra generated by the modes $\{a_{i,p}\}_{p\in \Z, i\in I}$. For $p=-\langle\alpha_i,\mu\rangle+1$, \eqref{eq:GKLOmodes} reduces to 
\begin{equation}\label{eq:actionHvsA}
h_{i,-\langle\alpha_j,\mu\rangle+1}=\textstyle\sum_{j\in I} -a_{ji} a_{j,1}+\sum_{j\in I} d_ja_{ji} m_j -e_1(R_i).
\end{equation}
The following definition is \cite[Definition~3.13]{kamnitzer2022hamiltonian} and it first appeared in \cite{kamnitzer2014yangians}.

\begin{definition}[Truncated shifted Yangian]\label{def:truncatedshiftedYangians}
The truncated shifted Yangian $Y_\mu^\lambda(\bR)$ is the quotient 
\begin{equation*}
Y_\mu^{\lambda}(\bR):=Y_\mu/I^\lambda_\mu(\bR)
\end{equation*}
by an ideal $I^\lambda_\mu(\bR)$ (whose definition can be found in \cite[Theorem~3.12]{kamnitzer2022hamiltonian}). The corresponding quotient map is denoted by $\pi_\mu^\lambda(\bR):Y_\mu\to Y_\mu^\lambda(\bR)$. 
\end{definition}

The only information about $I^\lambda_\mu(\bR)$ which will be needed here is that it contains the two-sided ideal $\langle a_{i,p} \;;\; p > m_i\rangle$. This follows immediately from the definition of $I^\lambda_\mu(\bR)$.

\begin{remark}
In type A, the ideal $I^\lambda_\mu(\bR)$ coincides with the two-sided ideal $\langle a_{i,p} \;;\; p > m_i\rangle$, see \cite[Theorem~1.6]{kamnitzer2014yangians} and \cite[Theorem~A.5]{kamnitzer2022hamiltonian}.
\end{remark}

In the quotient $Y_\mu^\lambda(\bR)$, the currents $a_i(u)$ are polynomials in $u^{-1}$ of degree $m_i$. For a $Y_\mu^\lambda(\bR)$-module $M$ and for $\bS=(S_i)_{i\in I}$ a set of parameters of height $\lambda-\mu$, define
\begin{equation}\label{eq:bSweightspace}
\scalebox{0.925}{$
W_\bS(M):=\big\{m\in M \;;\;\exists N>0, (a_{i,s}-(-1)^{s}e_{s}(S_i))^Nm=0 \text{ for all $i\in I$ and $1\leq s\leq m_i$}\big\},
$}
\end{equation}
i.e.~the subspace of generalized eigenvectors for the modes $a_{i,s}$ upon which they act with eigenvalue equal to $(-1)^{s}$ times the $s^\text{th}$ symmetric function of the multiset $S_i$.

\begin{remark}\label{rem:erreurWbS}
The reader should note that the definition of $W_{\bS}$ above does not match with the definition given in \cite[Section~5]{kamnitzer2019category} as there is a missing sign. However, \eqref{eq:bSweightspace} matches with the original definition of \cite[Section~3.6]{kamnitzer2019highest}.
\end{remark}

\subsection{The PBW basis and filtrations}\label{subsec:PBWbasisYmu}

The algebra $Y_\mu$ admits a PBW basis in the sense of \cite[Section~3.12]{finkelberg2018comultiplication}. The PBW basis elements are denoted by $e_{\beta,r},f_{\beta,r},h_{i,p}$ where $\beta\in \Delta_+$, $r\in \Z_{>0}$, $i\in I$ and $p\in \Z$, $p\geq  -\langle\alpha_i,\mu\rangle+1$. Moreover, for each splitting $\mu=\mu_1+\mu_2$, there is a filtration on $Y_\mu$ denoted by $F_{\mu_1,\mu_2}^\bullet$. It is defined by setting $F_{\mu_1,\mu_2}^k Y_\mu$ to be the $\C$-span of PBW variables of degree $\leq k$, where 
\begin{equation}\label{eq:splitmu}
\deg e_{\beta,r}=\langle\beta,\mu_1\rangle+r,\hspace{1em}\deg f_{\beta,r}=\langle\beta,\mu_2\rangle+r,\hspace{1em}\deg h_{i,p}=\langle\alpha_i,\mu\rangle+p.
\end{equation}
Properties of this filtration are recalled in \cite[Proposition~4.2]{kamnitzer2022hamiltonian}. This yields a filtration $F_{\mu_1,\mu_2,\bR}^\bullet$ on $Y_\mu^\lambda(\bR)$ by declaring $F_{\mu_1,\mu_2,\bR}^k Y_\mu^\lambda(\bR):=\pi_{\bR}(F_{\mu_1,\mu_2}^k Y_\mu)$
This quotient filtration is both separated and exhaustive. \par 

Using the PBW basis defined above, one sees that $Y_\mu$ admits a \textit{triangular decomposition}, i.e.~there is an isomorphism of vector spaces
\begin{equation}\label{eq:triangDecompYangian}
Y_\mu\simeq Y_\mu^-\otimes_\C Y_\mu^{=}\otimes_\C Y_\mu^+
\end{equation}
where $Y_\mu^-$ is the subalgebra of $Y_\mu$ generated by the modes $f_{i,r}$, $Y_\mu^=$ is the subalgebra generated by the modes $h_{i,p}$ and $Y_\mu^+$ is the subalgebra generated by the modes $e_{i,r}$.
 
It is shown in \cite[Theorem~4.7]{kamnitzer2022hamiltonian} and \cite[Theorem~5.15]{finkelberg2018comultiplication} that the algebras $\gr_{F_{\mu_1,\mu_2}}Y_\mu$ and $\gr_{F_{\mu_1,\mu_2,\bR}} Y_\mu^\lambda(\bR)$ are commutative and that they do not depend on the choice of splitting $\mu=\mu_1+\mu_2$. Moreover, the filtration $F_{\mu_1,\mu_2,\bR}^\bullet$ admits an expansion.

As a corollary of the previous discussion, one has : 

\begin{corollary}\label{cor:filtHypoYangian}
The hypotheses \ref{H:Afilt}, \ref{H:findimgenA}, \ref{H:commR} and \ref{H:filtExp} hold for the algebra $Y_\mu^\lambda(\bR)$. 
\end{corollary}

\subsection{Generalized affine Grassmannian slices}\label{subsec:affineGrassSlices}

Generalized affine Grassmannian slices are affine varieties which were first defined in \cite[Section~2(ii)]{braverman2016coulomb}. A definition, together with some of the main properties of those varieties as well as some notation which will extensively be used here can be found in \cite[Section~2]{kamnitzer2022hamiltonian}. We also introduce the affine Grassmannian of $G$, the spherical Schubert varieties of \cite{lusztig1983singularities} and state theorems relating them to representations of the Langlands dual group which are due to \cite{ginzburg1995perverse,beilinson1991quantization,mirkovic2007geometric}. We use terminology and notation of \cite{baumann2021mirkovic}. 

We briefly recall definitions concerning the geometry of $\Gr_G:=G\parparl t\parparr/G\brabral t\brabrar$, the affine Grassmannian of $G$. It is a reduced projective ind-scheme over $\C$ and the reader is referred to \cite{zhu2017introduction} for the details. For each $\mu\in P^\vee$, let $L_\mu$ denote the image of $t^\mu$ in $\Gr$. For a dominant coweight $\lambda\in P_+^\vee$, write $\Gr^\lambda:=G\brabral t\brabrar L_\lambda$. The Cartan decomposition of $G$ guarantees that 
\begin{equation*}
\Gr=\textstyle\bigsqcup_{\lambda\in P_+^\vee} \Gr^\lambda
\end{equation*}
and the \textit{spherical Schubert variety} associated to $\lambda\in P_+^\vee$ is defined as $\overline{\Gr^\lambda}$. Furthermore, for each $\mu \in P^\vee$, consider the \textit{semi-infinite orbit} $S_{\pm}^\mu:=N_{\pm}\parparl t\parparr L_\mu$. The intersection $\Gr^\lambda\cap S_-^\mu$ is non-empty precisely when $\mu\leq \lambda$. In this case, $\Gr^\lambda\cap S_-^\mu$ has pure dimension $\rht(\lambda-\mu)$ and the irreducible components of $\overline{\Gr^\lambda\cap S_-^\mu}$ are called \textit{MV cycles} of type $\lambda$ and weight $\mu$. For $\nu\in Q_+^\vee$, the set $S_+^0\cap S_-^{-\nu}$ is non-empty and it has pure dimension $\rht(\nu)$. The irreducible components of $\overline{S_+^0\cap S_-^{-\nu}}$ are called \textit{stable MV cycles} of weight $\nu$. Each stable MV cycle of weight $\nu$ has two distinguished $T$-fixed points $L_0,L_{-\nu}\in \Gr$ which are respectively called the top and bottom fixed points. 

One of the consequences of the celebrated \textit{geometric Satake isomorphism} is that there is an isomorphism of vector spaces 
\begin{equation}\label{eq:geomSatake}
V(\lambda)_\mu\simeq \HH_{\tp}\big(\overline{\Gr^\lambda}\cap S_-^\mu\big)
\end{equation}
which depends on the choice of principal nilpotent element $e\in \fn$ and on an isomorphism $\iota:\fh^\ast\to \fh$. This implies that MV cycles of type $\lambda$ and weight $\mu$ give a natural basis of the $\mu$-weight space of the representation of the Langlands dual group $V(\lambda)$. As per \cite[Proposition~3]{anderson2003polytope}, multiplication on the left by $t^{\lambda}$ provides a bijection 
\begin{equation*}
\scalebox{0.9}{$
\{\text{stable MV cycles $Z$ of weight $\lambda-\mu$ such that }t^\lambda Z\subset\overline{\Gr^\lambda}\}\xrightarrow[]{\sim}\{\text{MV cycles of type $\lambda$ and weight $\mu$}\}.
$}
\end{equation*}
Recall the family of maps $\{\psi_\lambda:V(\lambda)\to \C[N^\vee]\}_{\lambda\in P_+^\vee}$ defined in \eqref{eq:defpsilambda}. By \cite[Proposition~6.1]{baumann2021mirkovic}, for each stable MV cycle $Z$, there is a well-defined function $b_Z\in \C[N^\vee]$ such that if $t^\lambda Z\subset\overline{\Gr^\lambda}$, then $b_Z=\psi_\lambda([t^\lambda Z])$. Here, $[t^\lambda Z]$ is viewed as an element of $V(\lambda)_\mu$ using \eqref{eq:geomSatake}. This is the definition of basis \ref{b:MVbasis}. 

For $\mu\in P^\vee$, let $\cW_\mu$ be the subscheme of $G\parparl t^{-1}\parparr$ given by 
\begin{equation}\label{def:sliceWmu}
\cW_\mu:=N_1\brabral t^{-1}\brabrar T_1\brabral t^{-1}\brabrar t^{\mu} N_{-,1}\brabral t^{-1}\brabrar
\end{equation}
and for $\lambda\in P_+^\vee$, define the \textit{generalized affine Grassmannian slice} as 
\begin{equation*}
\cWbar{}^\lambda_\mu:=\cW_\mu\cap \overline{G[t]t^\lambda G[t]}.
\end{equation*}
This subscheme is non-empty precisely when $\lambda\geq \mu$. In this case, it is an affine variety of dimension $2\rht(\lambda-\mu)$ which admits an action of the torus $T$ together with a $\C^\times$-action by \textit{loop rotation}. Let $2\rho^\vee:\C^\times \to T$ be twice the Weyl vector. Under the $\C^\times$-action coming from $2\rho^\vee$, \cite[Lemma~2.8]{krylov2018integrable} shows that the set of fixed points satisfies 
\begin{equation}\label{eq:fixedPointsSlices}
(\cWbar{}^\lambda_\mu)^{\C^\times}=\begin{cases}
\{t^\mu\} & \mu \in \wt(V(\lambda))\\
\emptyset & \text{otherwise}
\end{cases}.
\end{equation}
Consequently, this result can be restated as saying that when $\mu\in \wt(V(\lambda))$, the $B$-algebra of $\C[\cWbar{}^\lambda_\mu]$ is local. The generalized affine Grassmannian slices are related to the geometry of the affine Grassmannian by \cite[Theorem~3.1]{krylov2018integrable} which states that there is a $T$-equivariant isomorphism
\begin{equation}\label{eq:isoVasya}
(\cWbar{}^\lambda_\mu)_- \simeq \overline{\Gr^\lambda}\cap S_-^\mu
\end{equation}
where the left-hand side of the above equation denotes the repelling set under the $\C^\times$-action. Consequently, top-dimensional irreducible components of $(\cWbar{}^\lambda_\mu)_-$ are in bijection with MV cycles of type $\lambda$ and weight $\mu$. and the unique $\C^\times$-fixed point of the irreducible components of $(\cWbar{}^\lambda_\mu)_- $ is sent to the bottom fixed point of the corresponding MV cycles.

Finally, recall from Section \ref{subsec:Dbarmap} that there are algebra maps
\begin{equation}\label{eq:Dbarmap}
\DbarPM:\C[N^\vee]\to \C(\fh^\vee)
\end{equation}
which, by \cite[Theorem~1.4]{baumann2021mirkovic}, satisfy (after choosing a principal nilpotent element and identifying $\fh^\ast\simeq \fh$) the equality 
\begin{equation*}
\DbarM(b_Z)=\epsilon^T_{L_{-\nu}}(Z)
\end{equation*}
where $Z$ is a stable MV cycle of weight $\nu$. Here, since $Z\subset\Gr$ has multiple $T$-fixed points, we add a subscript to the equivariant multiplicity to indicate it is taken with respect to which fixed point. Notice that on the left-hand side of the equality, both $\DbarM$ and $b_Z$ depend on the choice of principal nilpotent element, but the right-hand side does not. Using isomorphism \eqref{eq:isoVasya}, the equivariant multiplicity of an irreducible component of $(\cWbar{}^\lambda_\mu)_-$ at its unique fixed point is equal to the equivariant multiplicity of the corresponding stable MV cycle at its bottom fixed point.

\subsection{Torus actions}
	
The algebra $Y_\mu$ admits a Hamiltonian torus action, where the Cartan subalgebra $\fh$ of $\fg$ embeds in $Y_\mu$ as
\begin{align*}
i_\fh:\fh&\to Y_\mu\\
h_i&\mapsto \tfrac{1}{2}h_{i,-\langle\alpha_i,\mu\rangle+1}.
\end{align*}

\begin{remark}
The factor of $\tfrac{1}{2}$ is to account for the choice of $\hbar=2$. At $\mu=0$, the above embedding coincides with the restriction to $\fh$ of the embedding $U(\fg)\to Y(\fg)$.  
\end{remark}

Using the relations \eqref{eq:commHH},\eqref{eq:HEroot} and \eqref{eq:HFroot}, it is easy to check that 
\begin{equation*}
[i_\fh(h_i),h_{j,p}]=0,\hspace{1em}[i_\fh(h_i),e_{j,r}]=(\alpha_i,\alpha_j)e_{j,r},\hspace{1em}[i_\fh(h_i),e_{j,r}]=-(\alpha_i,\alpha_j)f_{j,r}
\end{equation*}
which shows that the torus $T\simeq (\C^\times)^r$ acts on $Y_\mu$, where $r$ is the rank of $\fg$. It is easy to see that $T$ preserves the filtration $F_{\mu_1,\mu_2}^\bullet$. Finally, the ideal $I_\mu^\lambda(\bR)$ is a $T$-submodule since it is the kernel of the $T$-equivariant map of \cite[Theorem~3.12]{kamnitzer2022hamiltonian}. Consequently, the algebra $Y_\mu^\lambda(\bR)$ is endowed with a Hamiltonian torus action, which shows: 

\begin{corollary}\label{cor:torusHypoYangian}
The hypotheses \ref{H:TactionA}, \ref{H:TactionHamiltonian} and \ref{H:TpreservesF} hold for the algebra $Y_\mu^\lambda(\bR)$.
\end{corollary}

As mentioned in Section \ref{subsec:PBWbasisYmu}, the commutative algebras $\gr_{F_{\mu_1,\mu_2}} Y_\mu$ and $\gr_{F_{\mu_1,\mu_2}} Y^\lambda_\mu(\bR)$ are independent of the choice of splitting $\mu=\mu_1+\mu_2$. In fact, \cite[Theorem~5.15]{finkelberg2018comultiplication} and \cite[Theorem~4.7]{kamnitzer2022hamiltonian} show that there are  isomorphisms of graded algebras 
\begin{equation}\label{eq:isoPoissonAlg}
\gr_{F_{\mu_1,\mu_2}} Y_\mu\simeq  \C[\cW_\mu]\hspace{2em}\text{and}\hspace{2em}\gr_{F_{\mu_1,\mu_2}} Y^\lambda_\mu(\bR)\simeq  \C[\cWbar{}^\lambda_\mu]
\end{equation}
which are in fact isomorphisms of Poisson algebras. 

\begin{remark}\label{rem:antiTequiv}
We warn the reader that the morphisms of \eqref{eq:isoPoissonAlg} are unfortunately anti-$T$-equivariant morphisms. This comes from the fact that if a set $X$ has a an action by a group $G$, the action on functions $f:X\to \C$ is given by $(gf)(x)=f(g^{-1}x)$. To remedy this, we will have to carry a sign through the upcoming computations.
\end{remark}

\subsection{The \texorpdfstring{$B$}{B}-algebra}

The cocharacter $2\rho^\vee$ endows $Y_\mu^\lambda(\bR)$ (and $Y_\mu$) with a $\Z$-grading. For this cocharacter, the $B$-algebra of $Y_\mu^\lambda(\bR)$ has been systematically studied in \cite{weekes2016highest,kamnitzer2019highest} where it was shown that is finite-dimensional when $\mu$ is dominant, see \cite[Proposition~3.12]{kamnitzer2019highest}. As one can expect, this extends to the general case.

\begin{lemma}
The algebra $Y_\mu^\lambda(\bR)$ is $B$-finite, i.e.~\ref{H:Bfinite} holds for $Y_\mu^\lambda(\bR)$. 
\end{lemma}

\begin{proof}
By Corollary \ref{cor:filtHypoYangian}, \ref{H:filtExp} holds for $Y_\mu^\lambda(\bR)$ and \eqref{eq:fixedPointsSlices} implies that $\C[\cWbar{}^\lambda_\mu]$ is local, i.e.~it satisfies \ref{H:BRlocal}. By Remark \ref{rem:hyposExpPlusLocal}, the result follows. 
\end{proof}

\begin{remark}
The previous proof shows that if $\mu\not\in \wt(V(\lambda))$, then $B(Y_\mu^\lambda(\bR))\simeq \{0\}$. 
\end{remark}

\begin{lemma}
The $B$-algebra of $Y_\mu^\lambda(\bR)$ is commutative, i.e.~\ref{H:Bcomm} holds for $Y_\mu^\lambda(\bR)$. 
\end{lemma}

\begin{proof}
Using the triangular decomposition \eqref{eq:triangDecompYangian}, one sees that $B(Y_\mu^\lambda(\bR))$ is a quotient of $Y_\mu^=$ which is a commutative algebra (see \eqref{eq:commHH}).
\end{proof}

In order to prove that $Y_\mu^\lambda(\bR)$ is $B$-integral under certain assumptions on $\bR$, we need following result which requires $\fg$ to be simply-laced since it revolves around technology developed in \cite{kamnitzer2019highest,kamnitzer2019category}. Some of these results admit generalizations to the non-simply-laced case, see \cite{varagnolo2025representations}.

\begin{theorem}[{\cite[Corollary~5.22]{kamnitzer2019category}}]\label{thm:monCrystalBalg}
Suppose $\fg$ is simply-laced and let $\bR$ be an integral set of parameters of level $\lambda$. There is a bijection between maximal ideals of $B(Y_\mu^\lambda(\bR))$ and elements of the $\mu$-weight space of the product monomial crystal $\cB(\lambda,\bR)_\mu$.
\end{theorem}

The bijection of the above theorem identifies maximal ideals of $B(Y_\mu^\lambda(\bR))$ with sets of parameters $\bS=(S_i)_{i\in I}$ of height $\lambda-\mu$ such that $y_{\bR}z_{\bS}^{-1}\in \cB(\lambda,\bR)_\mu$. For a set of parameters $\bS$, the corresponding maximal ideal is the kernel of the morphism $\phi_{\bS}$ defined by
\begin{equation}\label{eq:morphismBalgYangian}
\begin{aligned}
\phi_\bS:B(Y_\mu^\lambda(\bR))&\to \C,\\
a_{i,s}&\mapsto (-1)^{s}e_{s}(S_i).
\end{aligned}
\end{equation}
Using \eqref{eq:actionHvsA}, we compute the value of $\phi_{\bS}$ on the generators $h_i\in \fh$ via $i_\fh$. First, define
\begin{equation}\label{eq:actionAizero}
\textstyle\lambda(S_i):=\tfrac{1}{2}\Big(2e_1(S_i)-\sum_{j\sim i} e_1(S_j)+2m_i-\sum_{j\sim i} m_j-e_1(R_i)\Big).
\end{equation}
The weight of the $1$-dimensional $B(A)$-module defined by $\phi_\bS$ is $\lambda(\bS):=\sum_{i\in I} \lambda(S_i)\varpi_i$. We are now ready to prove:

\begin{lemma}\label{lem:truncShiftedBintegral}
Suppose $\fg$ is simply-laced. If $\bR$ is an integral set of parameters, the algebra $Y_\mu^\lambda(\bR)$ is $B$-integral, i.e.~\ref{H:Bintegral} holds for $Y_\mu^\lambda(\bR)$. 
\end{lemma}

\begin{proof} 
Notice that $e_1(S_j)$ is the sum of $m_j$ integers having the same parity as $j$ and $e_1(R_i)$ is the sum of $\lambda_i$ elements of the same parity as $i$. Thus, working modulo $2$ gives 
\begin{align*}
\textstyle\sum_{j\sim i} e_1(S_j) -\sum_{j\sim i} m_j -e_1(R_i)\,&\textstyle\equiv_2\, \sum_{j\sim i} j\cdot m_j-\sum_{j\sim i} m_j-i\cdot m_i\\
&\textstyle\equiv_2\,i\cdot (\lambda_i+\sum_{j\sim i} m_i)
\end{align*}
since $j\equiv_2 i+1$ when $i\sim j$. Crucially, the previous equation is independent of $\bS$ and it only depends on $\lambda$ and $\mu$. If one defines $\tau=\sum_{i\in I} \tau_i\varpi_i\in \fh^\ast$ by 
\begin{equation*}
\tau_i=\begin{cases}
1/2 & \text{if }i\cdot (\lambda_i+\sum_{j\sim i} m_i)\text{ is odd} \\
0 & \text{otherwise}
\end{cases}
\end{equation*}
then the shift functor $\Sigma_\tau$ makes every simple $B(Y_\mu^\lambda(\bR))$-module $X^\ast(T)$-graded.
\end{proof}

Combining the lemmas of the section yields: 

\begin{corollary}\label{cor:BalgebraHypoYangian}
Suppose $\fg$ is simply-laced. When $\bR$ is an integral set of parameters, the algebra $Y_\mu^\lambda(\bR)$ satisfies hypotheses \ref{H:Bfinite}-\ref{H:Bintegral}-\ref{H:BRlocal}-\ref{H:Bcomm}. 
\end{corollary}

\subsection{Category \texorpdfstring{$\cO_\mu^{\lambda}(\bR)$}{O for truncated shifted Yangians}}

Let $\bR$ be an integral set of parameters. Having checked the necessary hypotheses on $Y_\mu^\lambda(\bR)$ (see Corollaries \ref{cor:filtHypoYangian}, \ref{cor:torusHypoYangian} and \ref{cor:BalgebraHypoYangian}), we are now ready to study its category $\cO$ (see Definition \ref{def:catO}) which we will denote by $\cO_\mu^\lambda(\bR)$. By Lemma \ref{lem:LosevCatO} and Theorem \ref{thm:monCrystalBalg}, this category has $|\cB(\lambda,\bR)_\mu|$ simple objects. From \eqref{eq:morphismBalgYangian}, one can associate to each $\bS$ satisfying $y_{\bR}z^{-1}_{\bS}\in \cB(\lambda,\bR)_\mu$ a morphism $\phi_{\bS}:B\to \C$. Let $\Delta(\bS)$ be associated generalized Verma module. 

Since objects of $\cO_\mu^\lambda(\bR)$ are weight modules with finite-dimensional weight spaces, their weight spaces can be decomposed further using \eqref{eq:bSweightspace}. It follows that as vectors spaces, objects of $\cO_\mu^\lambda(\bR)$ can be decomposed as $M\simeq \oplus_{\bS} \,W_{\bS}(M)$.

\begin{definition}
For $\mu\in P^\vee$, let $\cO_\mu$ be the full subcategory of $Y_\mu\modu$ such that for all objects $M$ of $\cO_\mu$, there exists $\lambda=\sum_{i\in I} \lambda_i\varpi_i\geq \mu$ together with a set of parameters $\bR$ of level $\lambda$ such that $I_\mu^\lambda(\bR) M=0$ and $M\in\ob \cO_\mu^\lambda(\bR)$.
\end{definition}

By definition, objects of $\cO_\mu$ are $Y_\mu$-modules which are naturally modules over a truncated shifted Yangian and which lie in the corresponding category $\cO_\mu^\lambda(\bR)$. The reader should keep in mind that this definition for category $\cO$ for shifted Yangians is not equivalent to the definition of \cite[Section~3.3]{hernandez2021shifted} and that it is contains a priori more objects that category $\cO^{\text{fin}}$ of \cite[Section~9]{hernandez2021shifted}. However, as \cite[Theorem~8.4]{hernandez2021shifted} shows, both categories have the same simple objects.

By choosing $\mu=\mu+0$ in \eqref{eq:splitmu}, it follows that the subalgebra
\begin{equation*}
Y^{-,=}_\mu:=Y_{\mu}^{-}\otimes_{\C} Y_\mu^{=}\subset Y_\mu,
\end{equation*}
i.e.~the subalgebra generated by the modes $f_{i,r}$ and $h_{i,p}$, lies in $F_{\mu,0}^k$ for $k\geq 0$. Using the triangular decomposition \eqref{eq:triangDecompYangian}, one has:

\begin{corollary}
Hypothesis \ref{H:triangular} holds for $Y_\mu^\lambda(\bR)$. 
\end{corollary}

Let $\lambda\in P_+^\vee$ and choose a set of parameters $\bR$ of level $\lambda$. As previously, we write $\lambda-\mu=\sum_{i\in I} m_i\alpha_i^\vee\in Q_+^\vee$. Consider the two-sided ideal $J_\mu^\lambda(\bR)\subset Y^{-,=}_\mu$ defined by 
\begin{equation*}
J_\mu^\lambda(\bR):=\langle a_{i,p} \;;\; p>m_i\rangle
\end{equation*}
along with the quotient $Y^{-,=}_\mu(\bR):=Y^{-,=}_\mu/J_\mu^\lambda(\bR)$. Since the modes $a_{i,p}$ and $h_{i,p}$ generate the same commutative subalgebra of $Y_\mu$, we can use the same argument as in \eqref{eq:gensetYangian} in order to show that the algebra $Y^{-,=}_\mu(\bR)$ is finitely generated. The non-negative filtration $F_{\mu,0}^\bullet$ on $Y^{-,=}_\mu$ induces a non-negative filtration on $Y^{-,=}_\mu(\bR)$.

Given an object $M$ of $\cO^\lambda_\mu(\bR)$, its pullback to $Y_\mu$ is still finitely generated. Moreover, the triangular decomposition \eqref{eq:triangDecompYangian} guarantees that it is finitely generated as a $Y^{-,=}_\mu$-module. As it was indicated below Definition \ref{def:truncatedshiftedYangians}, the defining ideal of the projection $Y_\mu\twoheadrightarrow Y_\mu^\lambda(\bR)$ contains the modes $a_{i,p}$ for $i\in I$ and $p>m_i$. Consequently, $M$ is naturally a module over the algebra $Y^{-,=}_\mu(\bR)$. 

We are now ready to prove:

\begin{lemma}\label{lem:restrictionPreservesGKdim}
The restriction functor 
\begin{equation*}
\cO_\mu^\lambda(\bR)\to Y_\mu^{-,=}(\bR)\modu
\end{equation*}
preserves GK dimension. 
\end{lemma}

\begin{proof}
Let $M$ be an object of $\cO_\mu^\lambda(\bR))$ and choose a finite-dimensional subspace $M_0\subset M$ which generates $M$ under the action of $Y_\mu^{-,=}$. Let $A_0$ and $B_0$ finite-dimensional subspaces of $Y_\mu^{-,=}(\bR)$ and $Y_\mu^\lambda(\bR)$ respectively which are spanned by the generators described in \eqref{eq:gensetYangian}. By construction, we have 
\begin{equation*}
(A_0)^n M_0= (B_0)^n M_0 
\end{equation*}
for all $n\in \Z_{\geq 0}$ since we can order elements according to the PBW order. This immediately implies the statement.
\end{proof}

Using Remarks \ref{rem:boundedfiltGKdim} and \ref{rem:thmTauvel}, it follows that the category of finitely generated $Y_\mu^{-,=}(\bR)$-module satisfies \ref{HO:GKexact} and \ref{HO:GKassograded}. The previous lemma implies that $\cO_\mu^{\lambda}(\bR)$ also satisfies both of these hypotheses.

\subsection{Asymptotic characters for truncated shifted Yangians}\label{subsec:asymptoticcharYangians}

We now apply the main results of Section \ref{sec:asympchar}. Again, we fix $\lambda\in P_+^\vee$, $\mu\in P^\vee$ with $\mu\leq \lambda$,  $d=\rht(\lambda-\mu)$, and $\bR$ an integral set of parameters of level $\lambda$.

Recall that the variety $\cWbar{}^\lambda_\mu$ has complex dimension $2d$ and, by \eqref{eq:isoVasya}, its repelling set $(\cWbar{}^\lambda_\mu)_-$ has complex dimension $d$. Identify irreducible components of $(\cWbar{}^\lambda_\mu)_-$ with MV cycles of type $\lambda$ and weight $\mu$ using the isomorphism \eqref{eq:isoVasya}.

We will need a bound on the GK dimension of objects in $\cO_{\mu}^\lambda(\bR)$. Given a simple object $L$ of $\cO_{\mu}^\lambda(\bR)$, its GK dimension can be computed using the product monomial crystal $\cB(\lambda,\bR)$ and Theorem \ref{thm:monCrystalBalg}. If $L$ is associated to a monomial $y\in\cB(\lambda,\bR)$, \cite[Proposition~9.18]{kamnitzer2022lie} shows that $\gkdim L$ is equal to $\rht(\lambda'-\mu)$, where $y$ lies in the normal subcrystal $\cB(\lambda')\subset \cB(\lambda,\bR)$. Since $\cO_{\mu}^\lambda(\bR)$ is GK exact (see the previous section), its objects have GK dimension less than or equal to $d$. Thus, we write $\cO_{\mu}^\lambda(\bR)_{\tp}=\cO_{\mu}^\lambda(\bR)_{=d}$. 

As in Section \ref{subsec:GKexact_subcats}, consider the two homomorphisms of abelian groups 
\begin{align*}
\CC_{\leq \tp}:=\CC_{\leq d}:K_0(\cO_\mu^{\lambda}(\bR))&\to \HH_{\tp}\big((\cWbar{}^\lambda_\mu)_-\big)\\
\CC_{\tp}:=\CC_{=d}:K_0(\cO_\mu^{\lambda}(\bR)_{\tp})&\to \HH_{\tp}\big((\cWbar{}^\lambda_\mu)_-\big)
\end{align*}
where $X=\cWbar{}^\lambda_\mu$. Above, we suppressed the superscript corresponding to the reduced scheme-theoretic repelling set, but the reader should keep in mind that we consider the characteristic cycle supported on $\widetilde{(\cWbar{}^\lambda_\mu)}_-$.

As mentioned in Remark \ref{rem:antiTequiv}, the identification $\gr Y_\mu^\lambda(\bR)\simeq \C[\cWbar{}^\lambda_\mu]$ is anti-$T$-equivariant. This implies that, when passing to associated graded, the character of a $T$-equivariant $Y_\mu^\lambda(\bR)$-modules is related to the character of its associated graded module by the map $e^\lambda\mapsto e^{-\lambda}$. Let $\varsigma:\C(\fh)\to \C(\fh)$ be the linear map which sends a rational function $f$ to the rational function $\varsigma(f)=(h\mapsto f(-h))$. This slight modification could have been avoided by working with the Weyl vector $-2\rho^\vee$ (giving rise to the category $\cO^+$ of \cite{kamnitzer2019category}). 

The following result is a corollary of Theorem \ref{thm:achiandequivmult} together with \cite[Theorem~1.4]{baumann2021mirkovic}.

\begin{theorem}\label{thm:epsilonTversusachi}
Let $M$ be an object of $\cO_\mu^{\lambda}(\bR)$. Then, 
\begin{equation*}
\varsigma(\achi_d(M))=\textstyle\sum n_Z\,\epsilon^T_{L_{\mu}}(Z)=\textstyle\sum n_Z\, \DbarM(b_{t^{-\lambda}Z})
\end{equation*}
where the sum is over MV cycles $Z$ appearing with multiplicity $n_Z$ in the characteristic cycle $\CC_{\leq \tp}(M)$. Equivalently, the diagram 
\begin{equation*}
\begin{tikzcd}[column sep=1.2em,row sep=0.75em]
\HH_{\tp}\big((\cWbar{}^\lambda_\mu)_-\big)\arrow[d,"\psi_\lambda"'] & K_0(\cO_\mu^\lambda(\bR))\arrow[l,"\CC_{\leq \tp}"']\arrow[dd,"\varsigma\circ\achi_d"]\\
\C[N^\vee]_{-(\lambda-\mu)} \arrow[rd,"\DbarM"'] & \\
& \C(\fh) 
\end{tikzcd}
\end{equation*}
is commutative.
\end{theorem}

One can also apply Corollary \ref{cor:CCinj} and deduce that :

\begin{lemma}
If $\chi$ is injective on $K_0(\cO_\mu^{\lambda}(\bR))$, then $\CC_{\tp}$ is injective.
\end{lemma}

As mentioned previously, the geometric Satake isomorphism \eqref{eq:geomSatake} implies that the rank of the top-dimensional Borel--Moore homology $\HH_{\tp}\big((\cWbar{}^\lambda_\mu)_-\big)$ is precisely $\dim V(\lambda)_\mu$. From \cite[Corollary~9.19]{kamnitzer2022lie}, it follows that the rank of $K_0\big((\cO_\mu^{\lambda}(\bR))_{\tp}\big)$ is also $\dim V(\lambda)_\mu$. Hence, one has:

\begin{theorem}
If $\chi$ is injective on $K_0(\cO_\mu^{\lambda}(\bR))$, then the complexified homomorphism
\begin{equation*}
\mathrm{id}_\C\otimes \CC_{\tp}:\C\otimes K_0(\cO_\mu^{\lambda}(\bR)_{\tp})\to \C\otimes \HH_{\tp}((\cWbar{}^\lambda_\mu)_-)
\end{equation*}
is an isomorphism.
\end{theorem}

More generally, we expect that the top-dimensional characteristic cycle $\CC_{\tp}$ is an isomorphism for any $\lambda$, $\mu$ and $\bR$. Note that here, the hypothesis that $\bR$ is integral is crucial (see Remark \ref{rem:is_Bintegral_necessary}).

Notice that if the highest weights of the generalized Verma modules of $\cO_\mu^{\lambda}(\bR)$ are distinct, then $\chi$ is injective since Verma modules form a $\Z$-basis of $K_0(\cO_\mu^{\lambda}(\bR))$. This (sufficient but not necessary) criteria can be checked using the product monomial crystal.

Let $\cB$ be the monomial crystal for $\fg^\vee$. Consider the map 
\begin{equation*}
\varpi^\vee:\cB\to \fh^\ast
\end{equation*}
given by $\varpi^\vee(y_{i,k}):=\tfrac{k}{2}\varpi_i$ and which satisfies $\varpi^\vee(m\cdot m')=\varpi^\vee(m)+\varpi^\vee(m')$. From the definition, one has 
\begin{equation}\label{eq:varpiveez}
\varpi^\vee(z_{i,k})=(\tfrac{k}{2}+\tfrac{1}{2})\alpha_i.
\end{equation}

By \cite[Corollary~5.22]{kamnitzer2019category}, if $m\in \cB(\lambda,\bR)_\mu$, then $-\varpi^\vee(m)$ is the highest weight of the Verma module associated with this monomial under the embedding $i_\fh$. Using notation of \eqref{eq:actionAizero}, it is not hard to check that if $y_{\bR}z_{\bS}^{-1}\in \cB(\lambda,\bR)$, then $\lambda(\bS)=-\varpi^\vee(y_{\bR}z_{\bS}^{-1})$.

\begin{definition}
A set of parameters $\bR$ is called $\varpi^\vee$-parted if $|\varpi^\vee(\cB(\lambda,\bR))|=|\cB(\lambda,\bR)|$. 
\end{definition}

Hence, a set of parameters $\bR$ is $\varpi^\vee$-parted if $\varpi^\vee$ sends distinct monomials in $\cB(\lambda,\bR)$ to distinct points in $\fh^\ast$.

\begin{example}\label{ex:D4monomialCrystal}
Let $\fg^\vee=\mathfrak{so}_{8}$ (type $D_4$) and let $2$ be the trivalent node of the Dynkin diagram of $\fg$. Then, 
\begin{equation*}
\textstyle\cB(\varpi_2^\vee,c)_{\mu=0}=\big\{\,\tfrac{y_{1,c-1}}{y_{1,c-5}}\,,\,\tfrac{y_{2,c-2}}{y_{2,c-4}}\,,\,\tfrac{y_{3,c-1}}{y_{3,c-5}}\,,\,\tfrac{y_{4,c-1}}{y_{4,c-5}}\, \big\}
\end{equation*}
which shows that 
\begin{equation*}
\varpi^\vee(\cB(\varpi_2^\vee,c)_{\mu=0})=\big\{\,2\varpi_1,\varpi_2,2\varpi_3,2\varpi_4\,\big\}.
\end{equation*}
In the product monomial crystal 
\begin{equation*}
\cB(2\varpi_2,(R_1=\emptyset,R_2=\{c_1,c_2\},R_3=\emptyset,R_4=\emptyset))=\cB(\varpi_2^\vee,c_1)\cB(\varpi_2^\vee,c_2),
\end{equation*}
the two monomials 
\begin{equation*}
\textstyle\frac{y_{1,c_1-1}y_{4,c_2-1}}{y_{1,c_1-5}y_{4,c_2-5}}\hspace{1em}\text{and}\hspace{1em} \frac{y_{1,c_2-1}y_{4,c_1-1}}{y_{1,c_2-5}y_{4,c_1-5}}
\end{equation*}
are both sent to $2(\varpi_1+\varpi_4)$ under the map $\varpi^\vee$. It follows that for all $c_1,c_2\in \C$, the set of parameters $\bR=(\emptyset,\{c_1,c_2\},\emptyset,\emptyset)$ is never $\varpi^\vee$-parted.
\end{example}

In the case when $\varpi_i^\vee$ is minuscule, there is an explicit description of the monomials appearing in the crystal $\cB(\varpi_i^\vee,c)$ given in \cite[Proposition~2.10]{kamnitzer2019highest}. This results states that for $\mu\in W\varpi_i^\vee$, the (unique) monomial of coweight $\mu$ has the form 
\begin{equation*}
y_{\mu,c}:=\textstyle\prod_{j\in I} y^{\langle\alpha_j,\mu\rangle}_{j,c-k_{\mu,j}}
\end{equation*}
for some integer $k_{\mu,j}$ of the same parity as $j$ which is independent of $c$. It follows that 
\begin{equation}\label{eq:varpiveeYmuc}
\textstyle\varpi^\vee(y_{\mu,c})=\sum_{j\in I}\tfrac{1}{2} \langle\alpha_j,\mu\rangle(c-k_{\mu,j}) \cdot\varpi_j=\tfrac{1}{2} c\mu^\vee -\tfrac{1}{2}\sum_{j\in I} \langle\alpha_j,\mu\rangle k_{\mu,j} \varpi_j
\end{equation}
which, as a function of $c$, is non-constant. As Example \ref{ex:D4monomialCrystal} showed, this is not necessarily the case if $\varpi_i^\vee$ is not minuscule. The following lemma shows that when $\lambda$ is a sum of minuscule coweights, being $\varpi^\vee$-parted is a ``generic'' condition.  

\begin{lemma}\label{lem:hyperplanesAndVarpiParted}
Suppose that $\lambda=\sum_{k=1}^N \varpi_{i_k}^\vee\in P^\vee_+$ is a sum of minuscule coweights and let $\bR$ be a set of parameters of level $\lambda$. There exists finitely many affine hyperplanes in $\C^N$ such that if the elements of $\bR$ are chosen outside of those hyperplanes, then $\cB(\lambda,\bR)$ is $\varpi^\vee$-parted. 
\end{lemma}

\begin{proof}
Let $r:=|\textstyle\bigotimes_{j=1}^N\cB(\varpi_{i_j}^\vee)|$. Consider the map 
\begin{equation*}
\varpi^\vee_{i_1,\dots,i_N}:\C^N\to \mathcal{P}(\fh^\ast)\hspace{2em}(\mathcal{P}\text{ denotes the powerset})
\end{equation*}
defined by 
\begin{equation*}
\textstyle\varpi^\vee_{i_1,\dots,i_N}(c_1,\dots,c_N)=\varpi^\vee\Big(\prod_{j=1}^N \cB(\varpi_{i_j}^\vee,c_j)\Big).
\end{equation*}
If one chooses $c_1,\dots,c_N\in \C$ for which $c_{j_1}+\dots+c_{j_m}\not\in\tfrac{1}{2}\Z$ for all $m$ and for all $1\leq j_1<\dots <j_m\leq N$, then, by using \eqref{eq:varpiveeYmuc}, one has $|\varpi^\vee_{i_1,\dots,i_N}(c_1,\dots,c_N)|=r$. Consequently, the subset 
\begin{equation*}
\{(c_1,\dots,c_N)\in \C^N \;;\; |\varpi^\vee_{i_1,\dots,i_N}(c_1,\dots,c_N)|=r\}
\end{equation*}
is defined by finitely many (non-constant) linear polynomials. Thus, it is the union of finitely many affine hyperplanes. 
\end{proof}

Specializing to integral sets of parameters, one has:

\begin{corollary}
If $\lambda$ is a sum of minuscule coweights and if the elements of $\bR$ lie outside a finite collection of affine hyperplanes, then $\CC_{\tp}$ is injective.
\end{corollary}

\section{KLR algebras}\label{sec:KLR}

This section defines \textit{KLR algebras}, their \textit{cyclotomic quotients}, as well as \textit{KLRW algebras}. We explain how they are related to truncated shifted Yangians via \textit{parity KLRW algebras}. The definition of KLR and KLRW algebras are given for general $\fg$, but for the parity KLRW algebra, we require $\fg$ to be simply-laced. The main results appear in Section \ref{subsec:equivOfCats}. 

\subsection{Definitions and first properties}

We recall here the construction of two families of diagrammatic $\C$-algebras often referred to as KLR algebras (due independently to \cite{khovanov2009diagrammatic,khovanov2011diagrammatic} and \cite{rouquier20082}) and KLRW algebras (due to \cite{webster2017knot}). We also recall the construction of cyclotomic quotients of KLR algebras (also due to \cite{khovanov2009diagrammatic,khovanov2011diagrammatic}).

For both these families, elements of the algebras are finite $\C$-linear combinations of braid-like, decorated, planar diagrams which are considered up to isotopy. Relations in the algebras are local relations (which are given below in \eqref{eq:KLRrel1}-\eqref{eq:KLRWrel4}).  Each diagram is allowed to contain finitely many strands which are either black or red. Black strands are labelled by vertices of the Dynkin diagram of $\fg$. Red strands are labeled by dominant integral weights. Black strands can be decorated with finitely many dots. When a dot is labelled by a positive integer, say $n$, this is equivalent to labelling the given strand with $n$ dots. The diagrams are required to locally look like the following:
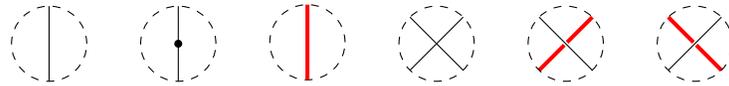
\begin{equation*}
\begin{tikzpicture}
\draw[dashed] (0,0) circle (0.5);
\draw (0,0.5) -- (0,-0.5);
\end{tikzpicture}
\hspace{2em}
\begin{tikzpicture}
\draw[dashed] (0,0) circle (0.5);
\draw (0,0.5) -- (0,-0.5);
\fill (0,0) circle (\circwidth);
\end{tikzpicture}
\hspace{2em}
\begin{tikzpicture}
\draw[dashed] (0,0) circle (0.5);
\draw[thick,red] (0,0.5) -- (0,-0.5);
\end{tikzpicture}
\hspace{2em}
\begin{tikzpicture}
\draw[dashed] (0,0) circle (0.5);
\draw (-0.3536,0.3536) -- (0.3536,-0.3536);
\draw (-0.3536,-0.3536) -- (0.3536,0.3536);
\end{tikzpicture}
\hspace{2em}
\begin{tikzpicture}
\draw[red,thick] (-0.3536,-0.3536) -- (0.3536,0.3536);
\draw[thick,white] (-0.3536,0.3536) -- (0.3536,-0.3536);
\draw (-0.3536,0.3536) -- (0.3536,-0.3536);
\draw[dashed] (0,0) circle (0.5);
\end{tikzpicture}
\hspace{2em}
\begin{tikzpicture}
\draw[red,thick] (-0.3536,0.3536) -- (0.3536,-0.3536);
\draw[thick,white] (-0.3536,-0.3536) -- (0.3536,0.3536);
\draw (-0.3536,-0.3536) -- (0.3536,0.3536);
\draw[dashed] (0,0) circle (0.5);
\end{tikzpicture}
\end{equation*}
Note that this implies that red strands are not allowed to cross. A diagram having only straight lines is called a straight-line diagram. Multiplication of two diagrams is defined as vertical concatenation: 
\begin{equation*}
\begin{tikzpicture}
\filldraw[rounded corners,fill opacity=0.1] (0,0) rectangle (1.5,0.75) node[opacity=1,black] at (0.75,0.375) {diag. 1};
\end{tikzpicture}\hspace{0.5em}
\raisebox{0.75em}{$\cdot$}\hspace{0.5em}
\begin{tikzpicture}
\filldraw[rounded corners,fill opacity=0.1] (0,0) rectangle (1.5,0.75) node[opacity=1,black] at (0.75,0.375) {diag. 2};
\end{tikzpicture}\hspace{0.5em}
\raisebox{0.75em}{$=$}\hspace{0.5em}
\raisebox{-1.1em}{\begin{tikzpicture}
\filldraw[rounded corners,fill opacity=0.1] (0,0) rectangle (1.5,0.75) node[opacity=1,black] at (0.75,0.375) {diag. 1};
\filldraw[rounded corners,fill opacity=0.1] (0,0) rectangle (1.5,-0.75) node[opacity=1,black] at (0.75,-0.375) {diag. 2};
\end{tikzpicture}}
\end{equation*}
In the above example, if the labeling on the bottom of diagram $1$ does not match the labeling on the top of diagram $2$, the result of the concatenation is declared to be zero. Otherwise, the resulting diagram is defined in the obvious way.

The relations \eqref{eq:KLRrel1}-\eqref{eq:KLRrel5} will be referred to as the KLR relations:

{\mathdiagramsetup
\begin{equation}\label{eq:KLRrel1}
\begin{tikzpicture}
\draw (1,1) -- (0,0) node[below] {$i$};
\draw (0,1) -- (1,0) node[below] {$i$};
\fill (0.25,0.75) circle (\circwidth);
\end{tikzpicture}
-
\begin{tikzpicture}
\draw (1,1) -- (0,0) node[below] {$i$};
\draw (0,1) -- (1,0) node[below] {$i$};
\fill (0.75,0.25) circle (\circwidth);
\end{tikzpicture}
=
\begin{tikzpicture}
\draw (0,1) -- (0,0) node[below] {$i$};
\draw (0.5,1) -- (0.5,0) node[below] {$i$};
\end{tikzpicture}
\hspace{3em}
\begin{tikzpicture}
\draw (1,1) -- (0,0) node[below] {$i$};
\draw (0,1) -- (1,0) node[below] {$i$};
\fill (0.25,0.25) circle (\circwidth);
\end{tikzpicture}
-
\begin{tikzpicture}
\draw (1,1) -- (0,0) node[below] {$i$};
\draw (0,1) -- (1,0) node[below] {$i$};
\fill (0.75,0.75) circle (\circwidth);
\end{tikzpicture}
=
\begin{tikzpicture}
\draw (0,1) -- (0,0) node[below] {$i$};
\draw (0.5,1) -- (0.5,0) node[below] {$i$};
\end{tikzpicture}
\end{equation}
\begin{equation}
\begin{tikzpicture}
\draw (1,1) -- (0,0) node[below] {$i$};
\draw (0,1) -- (1,0) node[below] {$j$};
\fill (0.25,0.75) circle (\circwidth);
\end{tikzpicture}
=
\begin{tikzpicture}
\draw (1,1) -- (0,0) node[below] {$i$};
\draw (0,1) -- (1,0) node[below] {$j$};
\fill (0.75,0.25) circle (\circwidth);
\end{tikzpicture}
\hspace{3em}
\begin{tikzpicture}
\draw (1,1) -- (0,0) node[below] {$i$};
\draw (0,1) -- (1,0) node[below] {$j$};
\fill (0.25,0.25) circle (\circwidth);
\end{tikzpicture}
=
\begin{tikzpicture}
\draw (1,1) -- (0,0) node[below] {$i$};
\draw (0,1) -- (1,0) node[below] {$j$};
\fill (0.75,0.75) circle (\circwidth);
\end{tikzpicture}
\hspace{1em}
\text{for }i\neq j
\end{equation}
\renewcommand{\arraystretch}{2.5}
\begin{equation}\label{eq:KLRpsisquare}
\begin{tikzpicture}
\draw plot [smooth, tension=1] coordinates {(0,1) (0.4,0.5) (0,0)} node[below] {$i$};
\draw plot [smooth, tension=1] coordinates {(0.5,1) (0.1,0.5) (0.5,0)} node[below] {$j$};
\end{tikzpicture}
=
\left\{
{\begin{tabular}{cl}	
\begin{tikzpicture}
\draw[opacity=0] (0,1) -- (0,0) node[below] {$i$};
\node at (0,0.5) {$0$};
\end{tikzpicture} &  \text{if }$i=j$\\
\begin{tikzpicture}
\draw (0,1) -- (0,0) node[below] {$i$};
\draw (0.5,1) -- (0.5,0) node[below] {$j$};
\end{tikzpicture} & \text{if }$(\alpha_i,\alpha_j)=0$\\
\begin{tikzpicture}
\draw (0,1) -- (0,0) node[below] {$i$};
\draw (0.5,1) -- (0.5,0) node[below] {$j$};
\fill (0,0.5) circle (\circwidth) node[left] {$-a_{ij}$};
\end{tikzpicture}
+
\begin{tikzpicture}
\draw (0,1) -- (0,0) node[below] {$i$};
\draw (0.5,1) -- (0.5,0) node[below] {$j$};
\fill (0.5,0.5) circle (\circwidth) node[right] {$-a_{ji}$};
\end{tikzpicture} & \text{if }$i\neq j$\text{ and }$(\alpha_i,\alpha_j)\neq 0$
\end{tabular}}\right.
\end{equation}

\begin{equation}
\begin{tikzpicture}
\draw plot [smooth, tension=1] coordinates {(0,1) (1,0)} node[below] {$k$};
\draw plot [smooth, tension=1] coordinates {(0.5,1) (0.25,0.5) (0.5,0)} node[below] {$j$};
\draw plot [smooth, tension=1] coordinates {(1,1) (0,0)} node[below] {$i$};
\end{tikzpicture}
=
\begin{tikzpicture}
\draw plot [smooth, tension=1] coordinates {(0,1) (1,0)} node[below] {$k$};
\draw plot [smooth, tension=1] coordinates {(0.5,1) (0.75,0.5) (0.5,0)} node[below] {$j$};
\draw plot [smooth, tension=1] coordinates {(1,1) (0,0)} node[below] {$i$};
\end{tikzpicture}
\hspace{1em}
\text{unless }i=k\text{ and }(\alpha_i,\alpha_j)\neq 0
\end{equation}
\begin{equation}\label{eq:KLRrel5}
\begin{tikzpicture}
\draw plot [smooth, tension=1] coordinates {(0,1) (1,0)} node[below] {$i$};
\draw plot [smooth, tension=1] coordinates {(0.5,1) (0.25,0.5) (0.5,0)} node[below] {$j$};
\draw plot [smooth, tension=1] coordinates {(1,1) (0,0)} node[below] {$i$};
\end{tikzpicture}
-
\begin{tikzpicture}
\draw plot [smooth, tension=1] coordinates {(0,1) (1,0)} node[below] {$i$};
\draw plot [smooth, tension=1] coordinates {(0.5,1) (0.75,0.5) (0.5,0)} node[below] {$j$};
\draw plot [smooth, tension=1] coordinates {(1,1) (0,0)} node[below] {$i$};
\end{tikzpicture}
=
\displaystyle\sum_{s+t+1=-a_{ij}}
\begin{tikzpicture}
\draw plot [smooth, tension=1] coordinates {(0,1) (0,0)} node[below] {$i$};
\draw plot [smooth, tension=1] coordinates {(0.5,1) (0.5,0)} node[below] {$j$};
\draw plot [smooth, tension=1] coordinates {(1,1) (1,0)} node[below] {$i$};
\fill (0,0.5) circle (\circwidth) node[left] {$s$};
\fill (1,0.5) circle (\circwidth) node[right] {$t$};
\end{tikzpicture}
\end{equation}
\noindent Fix $\lambda=\sum_{i\in I}\lambda_i\varpi_i\in P_+$. The relations \eqref{eq:KLRWrel1}-\eqref{eq:KLRWrel4} will be referred to as the KLRW relations: 
\begin{equation}\label{eq:KLRWrel1}
\begin{tikzpicture}
\draw[thick,red] plot [smooth, tension=1,thick,red] coordinates {(0.5,1) (0.25,0.5) (0.5,0)} node[below,black] {$\lambda$};
\draw plot [smooth, tension=1] coordinates {(0,1) (1,0)} node[below] {$j$};
\draw plot [smooth, tension=1] coordinates {(1,1) (0,0)} node[below] {$i$};
\end{tikzpicture}
-
\begin{tikzpicture}
\draw[thick,red] plot [smooth, tension=1] coordinates {(0.5,1) (0.75,0.5) (0.5,0)} node[below,black] {$\lambda$};
\draw plot [smooth, tension=1] coordinates {(0,1) (1,0)} node[below] {$j$};
\draw plot [smooth, tension=1] coordinates {(1,1) (0,0)} node[below] {$i$};
\end{tikzpicture}
=\delta_{i,j}\displaystyle\sum_{s+t+1=\lambda_i}
\begin{tikzpicture}
\draw[thick,red] plot [smooth, tension=1] coordinates {(0.5,1) (0.5,0)} node[below,black] {$\lambda$};
\draw plot [smooth, tension=1] coordinates {(0,1) (0,0)} node[below] {$i$};
\draw plot [smooth, tension=1] coordinates {(1,1) (1,0)} node[below] {$i$};
\fill (0,0.5) circle (\circwidth) node[left] {$s$};
\fill (1,0.5) circle (\circwidth) node[right] {$t$};
\end{tikzpicture}
\end{equation}

\begin{equation}
\begin{tikzpicture}
\draw[red,thick] (1,1) -- (0,0) node[below,black] {$\lambda$};
\draw (0,1) -- (1,0) node[below] {$i$};
\fill (0.25,0.75) circle (\circwidth);
\end{tikzpicture}
=
\begin{tikzpicture}
\draw[red,thick] (1,1) -- (0,0) node[below,black] {$\lambda$};
\draw (0,1) -- (1,0) node[below] {$i$};
\fill (0.75,0.25) circle (\circwidth);
\end{tikzpicture}
\hspace{3em}
\begin{tikzpicture}
\draw[red, thick] (0,1) -- (1,0) node[below,black] {$\lambda$};
\draw (1,1) -- (0,0) node[below] {$i$};
\fill (0.25,0.25) circle (\circwidth);
\end{tikzpicture}
=
\begin{tikzpicture}
\draw[red, thick] (0,1) -- (1,0) node[below,black] {$\lambda$};
\draw (1,1) -- (0,0) node[below] {$i$};
\fill (0.75,0.75) circle (\circwidth);
\end{tikzpicture}
\end{equation}

\begin{equation}
\begin{tikzpicture}
\draw[red, thick] plot [smooth, tension=1] coordinates {(0,1) (1,0)} node[below,black] {$\lambda$};
\draw plot [smooth, tension=1] coordinates {(0.5,1) (0.25,0.5) (0.5,0)} node[below] {$j$};
\draw plot [smooth, tension=1] coordinates {(1,1) (0,0)} node[below] {$i$};
\end{tikzpicture}
=
\begin{tikzpicture}
\draw[red, thick] plot [smooth, tension=1] coordinates {(0,1) (1,0)} node[below,black] {$\lambda$};
\draw plot [smooth, tension=1] coordinates {(0.5,1) (0.75,0.5) (0.5,0)} node[below] {$j$};
\draw plot [smooth, tension=1] coordinates {(1,1) (0,0)} node[below] {$i$};
\end{tikzpicture}
\hspace{3em}
\begin{tikzpicture}
\draw[red,thick] plot [smooth, tension=1] coordinates {(1,1) (0,0)} node[below,black] {$\lambda$};
\draw plot [smooth, tension=1] coordinates {(0,1) (1,0)} node[below] {$j$};
\draw plot [smooth, tension=1] coordinates {(0.5,1) (0.25,0.5) (0.5,0)} node[below] {$i$};
\end{tikzpicture}
=
\begin{tikzpicture}
\draw[red,thick] plot [smooth, tension=1] coordinates {(1,1) (0,0)} node[below,black] {$\lambda$};
\draw plot [smooth, tension=1] coordinates {(0,1) (1,0)} node[below] {$j$};
\draw plot [smooth, tension=1] coordinates {(0.5,1) (0.75,0.5) (0.5,0)} node[below] {$i$};
\end{tikzpicture}
\end{equation}

\begin{equation}\label{eq:KLRWrel4}
\begin{tikzpicture}
\draw[red,thick] plot [smooth, tension=1] coordinates {(0.5,1) (0.1,0.5) (0.5,0)} node[below,black] {$\lambda$};
\draw plot [smooth, tension=1] coordinates {(0,1) (0.4,0.5) (0,0)} node[below] {$i$};
\end{tikzpicture}
=
\begin{tikzpicture}
\draw[red,thick] (0.5,1) -- (0.5,0) node[below,black] {$\lambda$};
\draw (0,1) -- (0,0) node[below] {$i$};
\fill (0,0.5) circle (\circwidth) node[left] {$\lambda_i$};
\end{tikzpicture}
\hspace{3em}
\begin{tikzpicture}
\draw[red,thick] plot [smooth, tension=1] coordinates {(0,1) (0.4,0.5) (0,0)} node[below,black] {$\lambda$};
\draw plot [smooth, tension=1] coordinates {(0.5,1) (0.1,0.5) (0.5,0)} node[below] {$i$};
\end{tikzpicture}
=
\begin{tikzpicture}
\draw (0.5,1) -- (0.5,0) node[below] {$i$};
\draw[red,thick]  (0,1) -- (0,0) node[below,black] {$\lambda$};
\fill (0.5,0.5) circle (\circwidth) node[right] {$\lambda_i$};
\end{tikzpicture}
\end{equation}
}

\begin{definition}
The \textit{KLR algebra} associated to $\fg$, denoted by $R$, is the non-unital associative $\C$-algebra whose elements are finite $\C$-linear combinations of the above-described diagrams which contain exclusively black strands and which are subject to the KLR relations. The multiplication is given by vertical concatenation of diagrams.
\end{definition}

\begin{definition}
Let $\lambda\in P_+$ and fix a decomposition $\ulambda=(\lambda^{(1)},\dots,\lambda^{(\ell)})$ with $\lambda^{(k)}\in P_+$ for all $k=1,\dots,\ell$ and such that $\sum_{k=1}^\ell \lambda^{(k)}=\lambda$. The \textit{KLRW algebra} associated to $\fg$ and to the decomposition $\ulambda$, denoted by $\tildeT^{\ulambda}$, is the non-unital associative $\C$-algebra whose elements are finite $\C$-linear combinations of the above-described diagrams which contains exactly $\ell$ red strands labeled by $\lambda^{(1)},\dots,\lambda^{(\ell)}$ from left to right (note that this is well defined since red strands are not allowed to cross). The multiplication is given by vertical concatenation of diagrams.
\end{definition}

Both algebras $R$ and $\tildeT^{\ulambda}$ are $\Z$-graded, where the grading of a diagram is determined locally by the following rules: 
{\mathdiagramsetup
\begin{gather*}
\deg\big(
\begin{tikzpicture}
\draw (0,0.5) -- (0,-0.5) node[below] {$i$};
\fill (0,0) circle (\circwidth);
\end{tikzpicture}
\big)=(\alpha_i,\alpha_i)	
\hspace{1em}
\deg\big(
\begin{tikzpicture}
\draw (-0.3536,0.3536) -- (0.3536,-0.3536) node [below] {$j$};
\draw(0.3536,0.3536) --  (-0.3536,-0.3536) node [below] {$i$};
\end{tikzpicture}
\big)=-(\alpha_i,\alpha_j)\hspace{1em}
\deg\big(
\begin{tikzpicture}
\draw[red,thick] (0.3536,0.3536) -- (-0.3536,-0.3536) node [below,black] {$\lambda$};
\draw (-0.3536,0.3536) -- (0.3536,-0.3536) node [below] {$i$};
\end{tikzpicture}
\big)=\deg\big(
\begin{tikzpicture}
\draw[red,thick] (-0.3536,0.3536) -- (0.3536,-0.3536) node [below,black] {$\lambda$};
\draw (0.3536,0.3536) -- (-0.3536,-0.3536) node [below] {$i$};
\end{tikzpicture}
\big)=(\lambda,\alpha_i)
\end{gather*}
}Using \eqref{eq:cartanPairing}, one can routinely check that the relations \eqref{eq:KLRrel1}-\eqref{eq:KLRWrel4} are in fact homogeneous with respect to the above given degrees. 

\subsubsection{Properties of \texorpdfstring{$R$}{KLR algebras}}

The algebra $R$ comes naturally equipped with a family of (non-necessarily primitive) orthogonal idempotents. Let $\Seq$ denote the set of finite sequences in the alphabet $I$. For each $\bfi=(i_1,\dots,i_m)\in \Seq$, consider the straight-line diagram
{\mathdiagramsetup
\begin{equation*}
e(\bfi)=
\begin{tikzpicture}
\draw (0,1) -- (0,0) node[below] {$i_1$};
\draw (0.5,1) -- (0.5,0) node[below] {$i_2$};
\node at (1,0.5) {$\dots$};
\draw (1.4,1) -- (1.4,0) node[below] {$i_m$};
\end{tikzpicture}.
\end{equation*}
}This element is clearly an idempotent of $R$. Moreover, these elements satisfy $e(\bfi)e(\bfj)=0$ when $\bfi\neq \bfj$. For $\nu\in Q_+$ of height $\rht(\nu)=m$, consider the set of sequences of weight $\nu$ in the alphabet $I$, namely 
\begin{equation*}
\Seq(\nu)=\{(i_1,\dots,i_m)\in I^m \;;\; \textstyle \sum_{j=1}^m\alpha_{i_j} =\nu \}.
\end{equation*}
This set is finite and one can define $e(\nu)=\sum_{\bfi\in \Seq(\nu)} e(\bfi)$. The KLR algebra of weight $\nu$ is the idempotent sandwich $R_\nu:= e(\nu) R e(\nu)$. It is a unital $\C$-algebra and its unit is $e(\nu)$. \par 

\begin{example}
For $\fg=\sl_2$ and $\nu=n\alpha_1$, the algebra $R_\nu$ is isomorphic to the nilHecke algebra on $n$ generators.
\end{example}

There is a family of non-unital homomorphisms 
\begin{equation}\label{eq:iotaKLR}
\iota_{\nu_1,\nu_2}:R_{\nu_1}\otimes R_{\nu_2}\to R_{\nu_1+\nu_2}
\end{equation}
which are defined diagrammatically as horizontal concatenation of diagrams, i.e. 
\begin{equation*}
\begin{tikzpicture}
\filldraw[rounded corners,fill opacity=0.1] (0,0) rectangle (1.5,0.75) node[opacity=1,black] at (0.75,0.375) {diag. 1};
\end{tikzpicture}\hspace{0.5em}
\raisebox{0.75em}{$\otimes$}\hspace{0.5em}
\begin{tikzpicture}
\filldraw[rounded corners,fill opacity=0.1] (0,0) rectangle (1.5,0.75) node[opacity=1,black] at (0.75,0.375) {diag. 2};
\end{tikzpicture}\hspace{0.5em}
\raisebox{0.75em}{$\mapsto$}\hspace{0.5em}
{\begin{tikzpicture}
\filldraw[rounded corners,fill opacity=0.1] (0,0) rectangle (1.5,0.75) node[opacity=1,black] at (0.75,0.375) {diag. 1};
\filldraw[rounded corners,fill opacity=0.1] (1.5,0) rectangle (3,0.75) node[opacity=1,black] at (2.25,0.375) {diag. 2};
\end{tikzpicture}}\hspace{0.25em}\raisebox{0.75em}{.}
\end{equation*}
It sends the tensor product of the units $e(\nu_1)\otimes e(\nu_2)$ to an idempotent of $R_{\nu_1+\nu_2}$. For $\nu$ such that $\nu\geq \alpha_i$, using the map $\iota_{\alpha_i,\nu-\alpha_i}$, one can define functors
\begin{equation}\label{eq:defFoncteursEi}
\begin{gathered}
\cE_i^\nu:R_{\nu}\modu\to R_{\nu-\alpha_i}\modu\\
\cE_i^\nu=\iota_{\alpha_i,\nu-\alpha_i}\big(e(\alpha_i)\otimes e(\nu-\alpha_i)\big)(\trou)
\end{gathered}
\end{equation}
where the $R_{\nu-\alpha_i}$-module structure is given via the pullback by the canonical algebra morphism $R_{\nu-\alpha_i}\to R_{\alpha_i}\otimes R_{\nu-\alpha_i}$. Being the composition of exact functors, the functors $\{\cE_i^\nu\}_{i\in I}$ are exact. When $\nu-\alpha_i$ is not a sum of positive roots, let $\cE_i^\nu=0$ and define $\cE_i:=\bigoplus_{\nu\in Q_+}\cE_i^\nu$.

For $\lambda=\sum \lambda_i\varpi_i \in P_+$, let $\cI^\lambda_\nu$ be the two-sided ideal of $R_\nu$ generated by all dotted straight-line diagrams of the form
\begin{equation*}
\raisebox{-2.5em}{\begin{tikzpicture}
\draw (0,1) -- (0,0) node[below] {\tiny$i_1$};
\draw (0.5,1) -- (0.5,0) node[below] {\tiny$i_2$};
\node at (1,0.5) {$\dots$};
\draw (1.4,1) -- (1.4,0) node[below] {\tiny$i_{m-1}$};
\draw (1.9,1) -- (1.9,0) node[below] {\tiny$i_{m}$};
\fill (1.9,0.5) circle (\circwidth) node[right] {\small$\lambda_{i_m}$};
\end{tikzpicture}}
\end{equation*}
where $(i_1,\dots,i_m)$ runs over all elements of $\Seq(\nu)$. The quotient $R^\lambda_\nu:=R_\nu/\cI^\lambda_\nu$ is called the cyclotomic quotient at weight $\lambda$ of $R_\nu$. It is a finite-dimensional algebra \cite[Proposition~2.3]{lauda2011crystals}. 

\begin{remark}
The above convention for the ideal $\cI^\lambda_\nu$ differs from \cite{lauda2011crystals}.
\end{remark}

\subsubsection{Properties of \texorpdfstring{$\tildeT^{\ulambda}$}{KLRW algebras}}
Throughout this section, fix $\lambda\in P_+$ together with a decomposition $\ulambda=(\lambda^{(1)},\dots,\lambda^{(\ell)})$ with $\lambda^{(k)}\in P_+$ for all $k=1,\dots,\ell$ and such that $\sum_{k=1}^\ell \lambda^{(k)}=\lambda$. \par 

In an analogous way, the algebra $\tildeT^{\ulambda}$ also gives rise to a family of (non-necessarily primitive) orthogonal idempotents. For $\bfi=(i_1,\dots,i_m)\in \Seq$, consider the set of integer-valued weakly increasing functions $\kappa:[1,\ell]\to [0,m]$. Such $\kappa$ is called a \textit{charge}. For a charge $\kappa$, consider the diagram associated to the idempotent $e(\bfi)$, but where a red strand labeled $\lambda^{(k)}$ is added between the $\kappa(k)$-th and $\kappa(k+1)$-th black strands. This yields a straight-line diagram of the form
\begin{equation*}
e(\bfi,\kappa)=
\raisebox{-2.5em}{\begin{tikzpicture}
\draw (0,1) -- (0,0) node[below] {\small$i_1$};
\draw (0.5,1) -- (0.5,0) node[below] {\small$i_2$};
\node at (1,0.5) {$\dots$};
\draw (1.4,1) -- (1.4,0) node[below] {\small$i_{\kappa(k)}$};
\draw[thick,red] (1.9+0.25,0) -- (1.9+0.25,1)  node[above,black] {\small$\lambda^{(k)}$};
\draw (2.4+0.5,1) -- (2.4+0.5,0) node[below] {\small$i_{\kappa(k+1)}$};
\node at (2.9+0.5,0.5) {$\dots$};
\draw (3.3+0.5,1) -- (3.3+0.5,0) node[below] {\small$i_m$};
\end{tikzpicture}}.
\end{equation*}
Thus, the data of a charge $\kappa$ is equivalent to the data of the positions of the $\ell$ red strands with respect to the $m$ black strands. Note that when $\kappa(k+1)=0$, this means that the red strand labelled $\lambda^{(k)}$ is positioned to the left of all the black strands and when $\kappa(k)=m$, this means that the red strand labeled $\lambda^{(k)}$ is positioned to the right of all the black strands. The element $e(\bfi,\kappa)$ is clearly an idempotent of $\tildeT^{\ulambda}$. Moreover, these elements satisfy $e(\bfi,\kappa_1)e(\bfj,\kappa_2)=0$ whenever $(\bfi,\kappa_1)\neq (\bfj,\kappa_2)$.\par

For $\mu\in P$ such that $\nu=\lambda-\mu\in Q_+$, let $e(\nu,\tinybullet)=\sum_{\bfi\in \Seq(\nu)}\sum_{\kappa} e(\bfi,\kappa)$. The KLRW algebra associated to the pair $\ulambda$ and $\mu\in P$, where $\nu=\lambda-\mu\in Q_+$, is defined by the idempotent sandwich $\tildeT^{\ulambda}_\mu:=e(\nu,\tinybullet)\,\tildeT^{\ulambda}\,e(\nu,\tinybullet)$. It is a unital $\C$-algebra and its unit is $e(\nu,\tinybullet)$.\par

Consider the two-sided ideal $\cJ^{\ulambda}_\mu$ of $\tildeT^{\ulambda}_\mu$ generated by the idempotents $e(\bfi,\kappa)$ such that $\kappa(k)<m$ for all $k$ (the rightmost strand of these idempotents is black). The quotient 
\begin{equation*}
T^{\ulambda}_\mu:=\tildeT^{\ulambda}_\mu/\cJ^{\ulambda}_\mu
\end{equation*}
is called the \textit{cyclotomic} KLRW algebra. It is a finite-dimensional algebra.

Consider the idempotent of $\tildeT^{\ulambda}_\mu$ defined by
\begin{equation*}
e_{\cyc}:=\sum_{\substack{\bfi\in\Seq(\lambda-\mu)\\\kappa\;;\;\kappa(1)=m}} e(\bfi,\kappa)	
\end{equation*}
which is referred to as the \textit{cyclotomic idempotent}. In the diagram of $e_{\cyc}$, all the red strands are positioned on rightmost end of the diagram. Abusing notation, let $e_{\cyc}$ be the class of the cyclotomic idempotent in the quotient $T^{\ulambda}_\mu$.

\begin{lemma}[{\cite[Proposition~5.31]{webster2017knot}}]\label{lem:isoKLRcycloKLRW}
There is an isomorphism of $\C$-algebras 
\begin{equation*}
(e_{\cyc})T^{\ulambda}_\mu (e_{\cyc})\simeq R^\lambda_{\lambda-\mu}.
\end{equation*}
\end{lemma}

\subsection{Some categorification results}\label{subsec:someCategorification}

One of the main goals which motivated the introduction of both the algebras $R$ and $\tildeT^{\ulambda}$ was to categorify $U_q(\fn_-)$ as well as $\fg$-representations. We will only need partial information concerning these categorification results.

\begin{remark}
This paper will be concerned with \textit{ungraded} results. An $R$-module $M$ is called \textit{nilpotent} if the operators given by dotted straight-line diagrams act nilpotently on $M$. The category of finite-dimensional nilpotent $R$-modules will be denoted by $R\nilmod$. The same definition applies for $\tilde{T}^{\ulambda}\,$-modules together with modules over the quotients and idempotent sandwiches defined above. As shown in \cite[Lemma~3.13]{kamnitzer2019category}, a finite-dimensional module (over $R$ or $\tilde{T}^{\ulambda}$) admits a $\Z$-grading if and only if it is nilpotent. As noted in \cite[Section~2.7]{lauda2011crystals}, when working with a finite-dimensional $\Z$-graded algebras (e.g.~$\smash{R{}^{\lambda}_{\lambda-\mu}}$ or $\smash{T{}^{\ulambda}_{\mu}}$), the notion of nilpotent module is superfluous: any finitely generated module is nilpotent and admits a unique (up to shift and isomorphism) graded lift. 
\end{remark}

\begin{example}
Let $\fg=\sl_2$ and take $\nu=\alpha$. The algebra $R_\nu\simeq \C[x]$ has infinitely-many simple modules, but there is only one simple $\C[x]$-module on which $x$ acts nilpotently.
\end{example}

Following Section \ref{sec:notation}, for $\lambda\in P_+$, let $V(\lambda)$ denote the (unique) finite-dimensional, irreducible, $\fg$-representation of highest weight $\lambda$. Its $\mu$-weight space is denoted by $V(\lambda)_\mu$. Given a list $\ulambda=(\lambda^{(1)},\dots,\lambda^{(\ell)})$ such that $\lambda^{(k)}\in P_+$ for all $k$, define the tensor product representation $V(\ulambda):=V(\lambda^{(1)})\otimes \dots \otimes V(\lambda^{(\ell)})$. Note that $V(\ulambda)$ always contains a copy of the irreducible $V(\lambda)$, where $\lambda=\sum_k \lambda^{(k)}$. \par 

\subsubsection{Categorifications arising from KLR algebras}

The following theorem is a (partial) restatement in the ungraded case of the main theorem of \cite{khovanov2009diagrammatic,khovanov2011diagrammatic} and \cite{rouquier20082}.

\begin{theorem}\label{thm:RmuisoCNmu}
Upon passing to Grothendieck groups, the functors
\begin{equation*}
\cE_i:\bigoplus_{\nu\in Q_+} R_\nu\fmodu\to \bigoplus_{\nu\in Q_+} R_\nu\fmodu,\hspace{1em}i\in I
\end{equation*}
give a left $U(\fn)$-module structure to the vector space $\bigoplus_{\nu\in Q_+} \C\otimes_{\Z}K_0(R_\nu\nilmod)$. Moreover, there is a $Q_+$-graded isomorphism of $\C$-vector spaces
\begin{equation}\label{eq:isoRmuCNmu}
\gamma:\bigoplus_{\nu\in Q_+}\C\otimes_{\Z}K_0(R_\nu\nilmod)\xrightarrow[]{\sim} \C[N]
\end{equation}
which is compatible with the $U(\fn)$-action. Under this identification, the class $[M]$ of an $R_\nu$-module is mapped to the unique vector $d_M\in \C[N]_{-\nu}$ satisfying
\begin{equation}\label{eq:pairingOnSequence}
\langle e_{i_1}\dots e_{i_m} ,d_M\rangle=\dim\big( (\cE_{i_1}\circ\dots \circ \cE_{i_m})(M) \big)=\dim e\big(\,(i_m,\dots, i_1)\,\big)M
\end{equation}
for all sequences $\mathbf{i}=(i_1,\dots,i_m)\in \Seq(\nu)$.
\end{theorem}

The following result is the main theorem of \cite{kang2012categorification}, stated in the ungraded case.

\begin{theorem}\label{thm:categorificationLV}
The category $\bigoplus_{\nu\in Q_+} R^\lambda_{\nu}\modu$ categorifies
the $\fg$-representation $V(\lambda)$, where the action of the Chevalley generators $\{e_i\}_{i\in I}$ is given by the functors $\{\cE_i\}_{i\in I}$. Moreover, the pullback functor
\begin{equation*}
R^\lambda_{\lambda-\mu}\modu\to R_{\lambda-\mu}\nilmod
\end{equation*}
categorifies the $U(\fn)$-equivariant injective map $\psi_{\lambda}:V(\lambda)_\mu\to \C[N]_{-(\lambda-\mu)}$.
\end{theorem}

\subsubsection{Categorifications arising from KLRW algebras}

It is shown in \cite[Theorem~B]{webster2017knot} that the previous statements can be extended to tensor products of representations. Using an analogous version of \eqref{eq:defFoncteursEi} for KLRW algebras, one can construct restriction functors $\{\cE_i^{\ulambda}:T^{\ulambda}_\mu\modu\to T^{\ulambda}_{\mu+\alpha_i}\modu\}_{i\in I}$, see \cite{webster2017knot}.

\begin{theorem}
The category $\bigoplus_{\mu\in P} T^{\ulambda}_\mu\modu$ categorifies the tensor product representation $V(\ulambda)$, where the action of the Chevalley generators $\{e_i\}_{i\in I}$ is given by the functors $\{\cE_i^{\ulambda}\}_{i\in I}$. 
Moreover, the exact functor 
\begin{equation*}
e_{\cyc}(\trou): T^{\ulambda}_\mu\modu\to R^{\lambda}_{\lambda-\mu}\modu
\end{equation*}
categorifies the projection map $V(\ulambda)_\mu\twoheadrightarrow V(\lambda)_\mu$.
\end{theorem}
{\mathdiagramsetup
\begin{example}\label{ex:KLRWsl3}
The following builds upon Example \ref{ex:C[N]SL3}. Let $\fg=\sl_3$ and fix $\lambda=\varpi_1+\varpi_2$, $\ulambda=(\varpi_1,\varpi_2)$. The cyclotomic quotient $T^{\ulambda}_0$ has dimension $19$ (see \cite[p.45]{webster2017knot}). The non-zero idempotents $e(\bfi,\kappa)$ of this algebra are  
\begin{equation*}
\mathsf{e}_1=
\begin{tikzpicture}
\redstrandLabelled{1}{1}
\redstrandLabelled{1.5}{2}
\blackstrandLabelled{0}{2}
\blackstrandLabelled{0.5}{1}
\end{tikzpicture}
\hspace{0.5em},\hspace{0.5em}
\mathsf{e}_2=
\begin{tikzpicture}
\redstrandLabelled{1}{1}
\redstrandLabelled{1.5}{2}
\blackstrandLabelled{0}{1}
\blackstrandLabelled{0.5}{2}
\end{tikzpicture}
\hspace{0.5em},\hspace{0.5em}
\mathsf{e}_2'=
\begin{tikzpicture}
\redstrandLabelled{0.5}{1}
\redstrandLabelled{1.5}{2}
\blackstrandLabelled{0}{1}
\blackstrandLabelled{1}{2}
\end{tikzpicture}
\hspace{0.5em},\hspace{0.5em}
\mathsf{e}_3=
\begin{tikzpicture}
\redstrandLabelled{0}{1}
\redstrandLabelled{1.5}{2}
\blackstrandLabelled{0.5}{1}
\blackstrandLabelled{1}{2}
\end{tikzpicture}
\end{equation*}
and one sees that $T^{\ulambda}_0 \mathsf{e}_2\simeq T^{\ulambda}_0 \mathsf{e}_2'$. Computing directly shows that the projective modules $P_i:=T^{\ulambda}_0 \mathsf{e}_i$ are indecomposable and non-isomorphic. One also has $e_{\cyc}=\mathsf{e}_1+\mathsf{e}_2$. Another direct computation shows that the following set is a basis of $R_{\alpha_1+\alpha_2}^{\varpi_1+\varpi_2}$:
\begin{equation*}
\epsilon_1=
\begin{tikzpicture}
\draw (0,1) -- (0,0) node[below,black] {$2$};
\draw (0.5,1) -- (0.5,0) node[below] {$1$};
\end{tikzpicture}
\hspace{0.5em},\hspace{0.5em}
\epsilon_2=
\begin{tikzpicture}
\draw (0,1) -- (0,0) node[below,black] {$1$};
\draw (0.5,1) -- (0.5,0) node[below] {$2$};
\end{tikzpicture}
\hspace{0.5em},\hspace{0.5em}
\begin{tikzpicture}
\draw (0.5,1) -- (0,0) node[below,black] {$1$};
\draw (0,1) -- (0.5,0) node[below] {$2$};
\end{tikzpicture}
\hspace{0.5em},\hspace{0.5em}
\begin{tikzpicture}
\draw (0.5,1) -- (0,0) node[below,black] {$2$};
\draw (0,1) -- (0.5,0) node[below] {$1$};
\end{tikzpicture}
\hspace{0.5em},\hspace{0.5em}
\begin{tikzpicture}
\draw (0,1) -- (0,0) node[below,black] {$1$};
\draw (0.5,1) -- (0.5,0) node[below] {$2$};
\fill (0,0.5) circle (\circwidth);
\end{tikzpicture}
\hspace{0.5em},\hspace{0.5em}
\begin{tikzpicture}
\draw (0,1) -- (0,0) node[below,black] {$2$};
\draw (0.5,1) -- (0.5,0) node[below] {$1$};
\fill (0,0.5) circle (\circwidth);
\end{tikzpicture}
\end{equation*}
The isomorphism of Lemma \ref{lem:isoKLRcycloKLRW} is then easy to check. The algebra $R_{\alpha_1+\alpha_2}^{\varpi_1+\varpi_2}$ has two isoclasses of indecomposable projective modules $Q_i=R_{\alpha_1+\alpha_2}^{\varpi_1+\varpi_2}\epsilon_i$ and the two simple modules $L_i=\tp Q_i$ are $1$ dimensional. One can then check that the algebra is basic and connected. By writing down a basis of the $3$-dimensional projective $P_3$, it is not hard to check that $e_{\cyc}P_3\simeq L_2$ as $R_{\alpha_1+\alpha_2}^{\varpi_1+\varpi_2}$-modules.

The pullbacks of the modules $L_1,L_2$ are simple $R_{\alpha_1+\alpha_2}$-modules satisfying
\begin{equation*}
\dim e((1,2))L_1=0 \hspace{1em}\dim e((2,1))L_1=1 \hspace{1em}\dim e((1,2))L_2=1 \hspace{1em}\dim e((2,1))L_2=0.
\end{equation*}
Following Example \ref{ex:C[N]SL3}, one can compute the pairing $\langle,\rangle:U(\fn)\times\C[N]\to\C$ from \eqref{eq:paringUnCN} restricted to the $(\alpha_1+\alpha_2)$-weight space. One has
\begin{align*}
\langle e_1e_2,ab\rangle &=1 & \langle e_2e_1,ab\rangle &=1\\
\langle e_1e_2,c\rangle &=1 & \langle e_2e_1,c\rangle &=0
\end{align*}
which shows that $[L_1]\mapsto c$ and $[L_2]\mapsto ab-c$ under the map from Theorem \ref{thm:RmuisoCNmu}.
\end{example}
}

\subsection{Characters}
Following \cite[Section~8.2]{baumann2021mirkovic} and \cite[Section~8]{geiss2005semicanonical}, let $\mathfrak{f}$ be the free Lie algebra on the set $\{\hat{e}_i\;;\; i\in I\}$ and consider the $Q_+$-graded associative Hopf algebra $U(\mathfrak{f})$. For $\bfi=(i_1,\dots, i_d)\in \Seq(\nu)$, let $\hat{e}_{\bfi}:=\hat{e}_{i_1}\dots \hat{e}_{i_d}\in U(\mathfrak{f})$. For each $\nu\in Q_+$, the set $\{\hat{e}_{\bfi}\;;\; \bfi\in \Seq(\nu)\}$ is a basis for $U(\mathfrak{f})_{\nu}$. Consider the graded dual $(U(\mathfrak{f}))^\ast$ of $U(\mathfrak{f})$ and for each $\nu\in Q_+$, let $\{\hat{e}_{\bfi}^\ast\;;\; \bfi\in \Seq(\nu)\}\subset (U(\mathfrak{f}))^\ast_{-\nu}$ be the basis dual to $\{\hat{e}_{\bfi}\;;\; \bfi\in \Seq(\nu)\}$. It is easy to check that one can identify the $\C$-algebra $(U(\mathfrak{f}))^\ast$ with the $\C$-vector space $\C\Seq$ endowed with the shuffle product via $\hat{e}_{\bfi}^\ast\mapsto\bfi$.

For each $i\in I$, fix a choice of Chevalley generator $e_i\in \fn$. Then, there is a graded surjective Hopf algebra homomorphism $U(\mathfrak{f})\to U(\fn)$ defined by $\hat{e}_i\mapsto e_i$ which can be dualized to give an injective algebra homomorphism $(U(\fn))^\ast\to (U(\mathfrak{f}))^\ast\simeq \C\Seq$. Recall that there is a $Q_+$-graded $\C$-algebra isomorphism $\C[N]\simeq (U(\fn))^\ast$. Combining these facts with the isomorphism $\gamma$ from \eqref{eq:isoRmuCNmu} yields the \textit{character map}
\begin{equation}\label{eq:characterKLR}
\begin{aligned}
\ch : \textstyle\bigoplus_{\nu\in Q_+} K_0(R_\nu\fmodu)&\to \Z\Seq\\
[M]&\mapsto\textstyle\sum_{\bfi} \big(\dim e(\bfi)M \big)\,\bfi^{\rev}
\end{aligned}
\end{equation}
which is $\Z$-linear and injective by construction. Here, the map $(\trou)^{\rev}:\Seq\to \Seq$ reverses sequences, i.e. $\bfi=(i_1,\dots, i_d)\mapsto\bfi^{\rev}=(i_d,\dots,i_1)$. 

\begin{remark}
The sequence reversal appearing here is needed because of \eqref{eq:pairingOnSequence}. This is a consequence of the fact that the cyclotomic quotients considered are decorated by dots on the right-most strand. 
\end{remark}

Recall the $\DbarPM$ maps defined in Section \ref{subsec:Dbarmap} as well as the rational functions $\DbarPM_{\bfi}$. Moreover, recall the involutive automorphism $\varsigma:\C(\fh)\to \C(\fh)$ defined by $\varsigma(f)=(h\mapsto f(-h))$. Notice that the rational functions $\DbarPM_{\bfi}\in \C(\fh)$ are related by the equalities
\begin{equation}\label{eq:DbarPM_vs_varsigma}
\varsigma(\DbarP_{\bfi})=\DbarM_{\bfi^{\rev}}\quad\text{and}\quad\varsigma(\DbarM_{\bfi})=\DbarP_{\bfi^{\rev}}.
\end{equation}

Consider the map $\Z\Seq\to \C(\fh)$ defined by $\bfi\mapsto \DbarM_{\bfi}$ and extended $\Z$-linearly. By \cite[Proposition~8.4]{baumann2021mirkovic}, this map is a morphism of algebras. Composing with \eqref{eq:characterKLR} provides a \textit{bar-character} map
\begin{align*}
\chbar: \textstyle\bigoplus_{\nu\in Q_+}K_0(R_\nu\fmodu)&\to \C(\fh),\\
[M]&\mapsto\textstyle\sum_{\bfi} \big(\dim e(\bfi)M \big)\,\DbarM_{\bfi^{\rev}}.
\end{align*}
From \eqref{eq:DbarVSDbarbfi}, one has $\DbarM(\gamma([M]))=\chbar([M])$. The map $\chbar$ was first considered in \cite{casbi2021equivariant}.

The objective of the upcoming sections is to relate bar-characters with asymptotic characters for truncated shifted Yangians of Section \ref{sec:Yangian}. The following basic representation theoretic lemma will be useful.		

\begin{lemma}\label{lem:adjonctionIdempotents}
Let $A$ be a finite-dimensional $\C$-algebra and consider $e,e'\in A$ idempotents such that $ee'=e'e=e'$. Then, for $M$ an $A$-module, there are isomorphisms
\begin{equation*}
\Hom_A(Ae',M)\simeq \Hom_{eAe}(eAe',eM).
\end{equation*}
\end{lemma}	

\begin{proof}
As per \cite[Theorem~6.8]{assem2006elements}, the functors $(T_e,\res_e)$ form an adjoint pair, where 
\begin{align*}
T_e:=Ae\otimes_{eAe}(\trou):eAe\modu &\to A\modu & \res_e:=e(\trou):A\modu &\to eAe\modu.
\end{align*}
Moreover, there is an isomorphism of left $A$-modules $T_e(eAe')=Ae\otimes_{eAe}(eAe')\simeq Ae'$. The statement follows by adjunction.
\end{proof}

We will also need the following identity. 

\begin{lemma}\label{lem:DbarPrep}
Let $\beta_1,\dots, \beta_d\in Q_+\setminus\{0\}$. As generating series, there is an equality 
\begin{equation*}
\sum_{\substack{q_1<\dots<q_d<0\\ (q_1,\dots,q_d)\in \Z^d }} e^{q_1\beta_1+\dots +q_d\beta_d}=\dfrac{e^{-(d\beta_1+(d-1)\beta_2+\dots+\beta_d)}}{\prod_{k=1}^d (1-e^{-(\beta_1+\dots+\beta_k)})}.
\end{equation*}
\end{lemma}

\begin{proof}
The formula can be proved by induction on $d$. For $d=0$, the formula holds trivially. 
Now, for the induction step, assume the formula holds for sequences of length $d'<d$. Then, one has 
\begin{align*}
\sum_{q_1<\dots<q_d<0 } e^{q_1\beta_1+\dots +q_d\beta_d}&=\sum_{q_2<\dots<q_d<0} e^{q_2\beta_2+\dots +q_d\beta_d}\cdot\Big(\sum_{q_1<q_2}e^{q_1\beta_2}\Big)\\
&=\sum_{q_2<\dots<q_d<0} e^{q_2\beta_2+\dots +q_d\beta_d} \,\cdot\,\dfrac{e^{(q_2-1)\beta_1}}{1-e^{-\beta_1}}\\
&=\dfrac{e^{-\beta_1}}{1-e^{-\beta_1}}\,\sum_{q_2<\dots<q_d<0} e^{q_2(\beta_1+\beta_2)+q_3\beta_3+\dots +q_d\beta_d}.
\end{align*}
Notice that the right-most term of the previous equation is of the same form as the beginning expression but with $d-1$ elements of $Q_+$. One can apply the induction hypothesis and deduce that 
\begin{align*}
\sum_{q_1<\dots<q_d<0} e^{q_1\beta_1+\dots +q_d\beta_d}&=\dfrac{e^{-\beta_1}}{1-e^{-\beta_1}}\,\cdot\,\dfrac{e^{-\big((d-1)(\beta_1+\beta_2)+(d-2)\beta_3+\dots+\beta_d\big)}}{\prod_{k=2}^{d} (1-e^{-(\beta_1+\dots +\beta_{k})})}
\end{align*}
and the result follows.
\end{proof}

\subsection{Parity KLRW algebras}\label{subsec:parityKLRW}

This section is based on \cite[Section~3.4]{kamnitzer2019category}. Suppose that $\Gamma$ is simply laced and recall notation from Appendix \ref{sec:setofparameters}. Let $\bR$ be an integral set of parameters with $\uvarpi_\bR=(\varpi_{j_1},\dots,\varpi_{j_\ell})$. Then, one can define $\tildeT^{\bR}:=\tildeT^{\uvarpi_\bR}$.\par 

Defining the algebra using an integral set of parameters allows one to associate extra data to the idempotents $e(\bfi,\kappa)$. Every red strand of an idempotent $e(\bfi,\kappa)\in \tildeT^{\bR}$ acquires a labeling by a pair $(i,r)$, where $i\in I$ is a vertex of the Dynkin diagram and $r$ is an integer of the same parity as $i$. The vertex $i\in I$ corresponds to the index of the fundamental weight labeling the red strand. The integer $r$ corresponds the element of the multi-set $R_i$ which produced the fundamental weight $\varpi_i$ in the list $\uvarpi_\bR$. The integer $r$ is called the \textit{longitude} of the red strand.\par 

Using the labeling of the red strands by elements of $I$ as well as the decomposition $I=I_0\sqcup I_1$ into even and odd vertices, one can define the \textit{parity distance} between strands of an idempotent $e(\bfi,\kappa)$. For two successive strands $p$ and $p'$ where $p$ is to the left of $p'$, define 
\begin{equation*}
\delta(p,p')=\begin{cases}
2 & \text{if $p$ and $p'$ have the same parity, $p$ is a black strand and $p'$ is a red strand}\\
1 & \text{if $p$ and $p'$ have different parity}\\
0 & \text{otherwise}
\end{cases}.
\end{equation*}
Now, extend the definition of $\delta$ to any two strands by declaring $\delta(p,p'')=\delta(p,p')+\delta(p',p'')$ if $p$ is to the left of $p'$ which is to the left of $p''$. The quantity $\delta(p,p')$ is called the parity distance between the strands $p$ and $p'$. An idempotent $e(\bfi,\kappa)$ is called a parity idempotent if, for any two red strands $p$ and $p'$ of $e(\bfi,\kappa)$, one has $\delta(p,p')\leq |r-r'|$, where $r$ and $r'$ are the longitudes of the red strands $p$ and $p'$ respectively. Let $e_{\bR}$ denote the sum of all the parity idempotents of $\tildeT^{\bR}$.
{\mathdiagramsetup
\begin{example}\label{ex:parityIdempotents}
Let $\fg=\sl_3$ and let $R_1=\{a\}$, $R_2=\{b\}$ with $a<b$.  One has $\uvarpi_\bR=(\varpi_1,\varpi_2)$ and the red strand labelled $1$ has longitude $a$ while the red strand labelled $2$ has longitude $b$. In the algebra $\tildeT^{\bR}_{\mu=0}$, the idempotent 
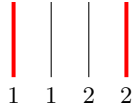
\begin{equation*}
\begin{tikzpicture}
\redstrandLabelled{0}{1}
\redstrandLabelled{1.5}{2}
\blackstrandLabelled{0.5}{1}
\blackstrandLabelled{1}{2}
\end{tikzpicture}
\end{equation*} 
is parity if and only if $3\leq b-a$.
\end{example}
}

\begin{lemma}\label{lem:cycloIsParity}
Idempotents $e(\bfi,\kappa)$ of $\tildeT^\bR$ with $\kappa(1)=m$ are parity idempotents.
\end{lemma}

\begin{proof}
Two successive red strands having different parity must have longitudes of different parity. Hence, their longitudes differ by at least $1$. 
\end{proof}

\begin{definition}
The \textit{parity KLRW algebra} $P^\bR$ is the quotient of the idempotent sandwich algebra $e_{\bR} \tildeT^{\bR} e_{\bR}$ by the ideal generated by parity idempotents $e(\bfi,\kappa)$ such that $\kappa(\ell)<m$. 
\end{definition}

As for the KLRW algebra, put $P^{\bR}_\mu:=e(\lambda-\mu,\tinybullet) P^{\bR}e(\lambda-\mu,\tinybullet)$. Similarly to the cyclotomic KLRW algebras, $P^{\bR}_\mu$ is finite-dimensional.

Given an integral set of parameters $\bS=(S_i)_{i\in I}$ whose cardinality vector is $(n_1,\dots, n_r)$ with $n=\sum n_i$, one can construct an idempotent $e(\bfi,\kappa)\in \tildeT^{\bR}$ as follows. First, order the elements of $\bS$ as $s_1\leq \dots \leq s_n$. Now, consider the sequence $\bfi=(i_1,\dots,i_n)\in \Seq(\sum_{i\in I} n_i \alpha_i)$, where $i_j$ is such that $s_{j}\in S_{i_j}$. From the straight line diagram $e(\bfi)$, declare $s_j$ to be the longitude of the strand $i_j$ in $e(\bfi)$. Then, interlace the red and black strands according to their longitudes, where longitudes should be weakly increasing from left to right and where, if the longitude of a red strand agrees with the longitude of a black strand, the red strand goes to the left of the black strand. This determines a charge $\kappa$. Denote the resulting idempotent by $e(\bS)$.

When $\bS$ is singular, note that, when ordering the elements $s_1\leq \dots \leq s_n$, there exists $j$ such that $s_j=s_{j+1}$. In this case, we have to choose  a total order on $s_1\leq \dots \leq s_n$. We call such a choice a \textit{resolution of multiplicities} of $\bS$. Thus, when $\bS$ is singular, $e(\bS)$ is well-defined up to a resolution of multiplicities. Let $\bS$ be singular and suppose $s_j=s_{j+1}$ with $i_j\neq i_{j+1}$. Let $e_1, e_2\in \tildeT^{\bR}$ be two possible idempotents given by the above construction. Since $s_j=s_{j+1}$, the vertices $i_j$ and $i_{j+1}$ have the same parity and there is no edge between the vertices $i_j$ and $i_{j+1}$ in $\Gamma$. The diagram
\begin{equation*}
\begin{tikzpicture}
\node at (-1.25,0.5) {$\dots$};
\draw plot [smooth, tension=1] coordinates {(-0.75,1) (-0.75,0)} node[below] {$i_{j-1}$};
\draw plot [smooth, tension=1] coordinates {(0,1) (1,0)} node[below] {$i_{j+1}$};
\draw plot [smooth, tension=1] coordinates {(1,1) (0,0)} node[below] {$i_{j}$};
\draw plot [smooth, tension=1] coordinates {(1.75,1) (1.75,0)} node[below] {$i_{j+2}$};
\node at (2.25,0.5) {$\dots$};
\end{tikzpicture}
\end{equation*}
squares to a straight-line diagram by \eqref{eq:KLRpsisquare}. Interlacing it with red strands, it provides an isomorphism $\tildeT^{\bR}e_1\simeq  \tildeT^{\bR}e_2$.

\begin{lemma}[{\cite[Lemma~3.25]{kamnitzer2019category}}]
Given an idempotent $e(\bfi,\kappa)\in \tildeT^{\bR}$, there exists an integral set of parameters $\bS$ such that, upon choosing a resolution of multiplicities, $e(\bfi,\kappa)=e(\bS)$ if and only if $e(\bfi,\kappa)$ is a parity idempotent. 
\end{lemma}

By the previous discussion and the above lemma, for a $P^{\bR}$-module $M$, the quantity $\dim e(\bS)M$ is well-defined and independent of the resolution of multiplicities.

\subsection{An equivalence of categories}\label{subsec:equivOfCats}

Fix $\lambda=\sum_{i\in I} \lambda_i\varpi\in P_+$ along with $\mu\in P$ and set $\nu:=\lambda-\mu=\sum_{i\in I}m_i\alpha_i\in Q_+$. Choose $\bR$ an integral set of parameters of level $\lambda$. Recall from Section \ref{sec:Yangian} that this data allows one to define the truncated shifted Yangian $Y_\mu^\lambda(\bR)$ (for the Lie algebra $\fg^\vee$) together with the category $\cO_\mu^\lambda(\bR)$. 

The following theorem is one of the main result of \cite{kamnitzer2019category}.

\begin{theorem}[{\cite[Theorem~1.2,~5.2]{kamnitzer2019category}}]\label{thm:equivTheta}
There is an equivalence of categories 
\begin{equation*}
\Theta:\cO_\mu^\lambda(\bR) \xrightarrow[]{\sim} P_{\mu}^\bR\modu
\end{equation*}
which sends the simple top of $\Delta(\bS)$ to the simple top of the projective module $P_{\mu}^\bR e(\bS)$. Moreover, for $M\in \ob\big(\cO_\mu^\lambda(\bR)\big)$ and a set of parameters $\bS$, one has
\begin{equation}\label{eq:weightSpace}
\dim W_\bS(M)=\tfrac{1}{\sigma(\bS)} \dim e(\bS) \Theta(M).
\end{equation}
In particular, sets of parameters $\bS$ such that $W_{\bS}(M)\neq 0$ are necessarily integral sets of parameters.
\end{theorem}

\begin{remark}
We warn the reader about a mistake in \cite[Definition~5.20]{kamnitzer2019category} (``$+$" and ``$-$" should be interchanged in the definition of $\cO^{\pm}$, see the proof of \cite[Theorem~5.21]{kamnitzer2019category}). However, this sign issue cancels out miraculously with the sign issue which was pointed out in Remark \ref{rem:erreurWbS}.
\end{remark}

Building on Theorem \ref{thm:equivTheta}, it is possible to rewrite characters of objects of $\cO_\mu^\lambda(\bR)$ in a way which makes the connection with KLR algebras more apparent. First, we recall the following technical computation that was made in the previous sections.

\begin{lemma}\label{lem:weightOfbS}
Let $M$ be an object of $\cO_\mu^\lambda(\bR)$. Under the pullback by $i_\fh$, the weight of vectors in $W_\bS(M)$ is given by
\begin{equation*}
\textstyle\lambda(\bS)=\tfrac{1}{2}\Big(\sum_{i\in I} e_1(S_i)\alpha_i\Big)+\tfrac{1}{2}(\lambda-\mu)-\tfrac{1}{2}\Big(\sum_{i\in I}e_1(R_i)\varpi_i\Big).
\end{equation*}
\end{lemma}

\begin{proof}
This follows from the proof of \eqref{eq:actionAizero}, which in turn relies on \eqref{eq:actionHvsA}. 
\end{proof}

The following example shows how to combine \eqref{eq:weightSpace} and Lemma \ref{lem:weightOfbS} to compute characters.
This will be used in the proof of the main theorem of this section.
{\mathdiagramsetup
\begin{example}\label{ex:calculCaracteresl3}
Continuing Examples \ref{ex:KLRWsl3} and \ref{ex:parityIdempotents}, let $R_1=\{a\}$, $R_2=\{a+3\}$. The algebra $P_{0}^\bR$ is isomorphic to $T_{0}^\bR$ since every idempotent is parity. The representation theory of $T_{0}^\bR$ was laid out in Example \ref{ex:KLRWsl3}. Let $M$ be the object of $\cO_{0}^{\varpi_1+\varpi_2}(\bR)$ corresponding to $P_3$ under the equivalence $\Theta$. Its character can be computed using longitudes of idempotents of $P_{0}^\bR$. A basis of $P_3$ is 
\begin{equation*}
v_1=
\begin{tikzpicture}
\redstrandLabelled{0}{1}
\redstrandLabelled{1.5}{2}
\blackstrandLabelled{0.5}{1}
\blackstrandLabelled{1}{2}
\end{tikzpicture}
\hspace{1em},\hspace{1em}
v_2=
\begin{tikzpicture}
\redstrandXXLabelled{0}{0.5}{1}
\redstrandXXLabelled{1.5}{1.5}{2}
\blackstrandXXLabelled{0.5}{0}{1}
\blackstrandXXLabelled{1}{1}{2}
\end{tikzpicture}
\hspace{1em},\hspace{1em}
v_3=
\begin{tikzpicture}
\redstrandXXLabelled{0}{1}{1}
\redstrandXXLabelled{1.5}{1.5}{2}
\blackstrandXXLabelled{1}{0.5}{2}
\blackstrandXXLabelled{0.5}{0}{1}
\end{tikzpicture}
\end{equation*}
from which it follows that 
\begin{enumerate}[label=\arabic*.]
\item $e(\bS)v_1\neq 0$ if and only if $S_1=\{a\}$ and $S_2=\{a+1\}$, 
\item $e(\bS)v_2\neq 0$ if and only if $\bS=(\{s_1\},\{a+1\})$ with $s_1\leq a-2$ and $\bS$ is integral, 
\item $e(\bS)v_3\neq 0$ if and only if $\bS=(\{s_1\},\{s_2\})$ with $s_1\leq s_2\leq a-1$  and $\bS$ is integral.
\end{enumerate}
For $\bS=(\{s_1\},\{s_2\})$, \eqref{eq:actionAizero} gives 
\begin{equation*}
\lambda(\bS)=\tfrac{1}{2}(s_1\alpha_1+s_2\alpha_2+(\alpha_1+\alpha_2)-a\varpi_1-(a+3)\varpi_2).
\end{equation*}
The above three cases respectively give the weights
\begin{enumerate}[label=\arabic*.]
\item $\lambda(\bS)=0$, 
\item $\lambda(\bS)=-(k+1)\alpha_1$ with $k\in \Z_{\geq 0}$,
\item $\lambda(\bS)=-(k_1+1)\alpha_1-(k_2+1)\alpha_2$  with $k_1,k_2\in \Z_{\geq 0}$, $k_1\geq k_2$.
\end{enumerate}
Consequently, the character of $M$ is 
\begin{align*}
\chi(M)&=\textstyle 1+\Big(\sum_{k\geq 0} e^{-(k+1)\alpha_1}\Big)+\Big(\sum_{k_1\geq k_2\geq 0} e^{-(k_1+1)\alpha_1-(k_2+1)\alpha_2}\Big)\\
&=1+\dfrac{e^{-\alpha_1}}{1-e^{-\alpha_1}}+\dfrac{e^{-\alpha_1-\alpha_2}}{(1-e^{-\alpha_1})(1-e^{-\alpha_1-\alpha_2})}=\dfrac{1}{(1-e^{-\alpha_1})(1-e^{-\alpha_1-\alpha_2})}.
\end{align*}
The reader well acquainted with parabolic category $\cO$ should not be surprised to see the appearance of the character of a parabolic Verma module. It is a consequence of the \textit{quantum Mirkovi\'{c}-Vybornov isomorphism} from \cite{webster2020quantum}.
\end{example}
}
By Lemma \ref{lem:cycloIsParity}, every idempotent appearing in $e_{\cyc}\in T^\bR_\mu$ also appears in $e_{\bR}\in T^\bR_\mu$. It follows that one can consider $e_{\cyc}$ as an idempotent of $P^\bR_\mu$. We write
\begin{equation*}
e_{\cyc}(\trou):P_{\mu}^\bR\modu\to R_{\lambda-\mu}^\lambda\modu
\end{equation*}
for the corresponding idempotent restriction functor. 

Let $\Theta_{\cyc}:\cO_\mu^\lambda(\bR)\to R_{\lambda-\mu}^\lambda\modu$ be the composition of the two functors
\begin{equation}\label{eq:definitionThetaCyc}
\begin{tikzcd}
\cO_\mu^\lambda(\bR)\arrow[r,"\Theta"] & P_{\mu}^\bR\modu \arrow[r,"e_{\cyc}(\trou)"] & R_{\lambda-\mu}^\lambda\modu.
\end{tikzcd}
\end{equation}
Recall the notion of GK-subcategories of order $d$ from Section \ref{subsec:filteredquant}.

\begin{proposition}[{\cite[Proposition~9.19]{kamnitzer2022lie}}]\label{prop:kernelThetacyc}
The kernel of the functor $\Theta_{\cyc}$ is the subcategory $(\cO_\mu^\lambda(\bR))_{<d}$ where $d=\rht(\lambda-\mu)$. Consequently, the induced functor 
\begin{equation*}
\Theta_{\cyc}:(\cO_\mu^\lambda(\bR))_{\tp}\to  R_{\lambda-\mu}^\lambda\modu
\end{equation*}
is an equivalence of categories.
\end{proposition}

The main result of the section is:
\begin{theorem}\label{thm:DbarMvsachi}
Let $M$ be an object of $\cO_{\mu}^\lambda(\bR)$. Then, 
\begin{equation}\label{eq:equationDbarMvsachi}
\varsigma(\achi_d(M))=\chbar\big(\Theta_{\cyc}(M)\big).
\end{equation}
Equivalently, the diagram
\begin{equation*}
\begin{tikzcd}[column sep=1.2em,row sep=0.75em]
K_0(\cO_\mu^\lambda(\bR)) \arrow[r,"\Theta_{\cyc}"]\arrow[dd,"\varsigma\circ\achi_d"'] & K_0(R_{\lambda-\mu}^\lambda\modu)\arrow[d,"\gamma"]\\
& \C[N]_{-(\lambda-\mu)}\arrow[dl,"\DbarM"]\\
\C(\fh) & 
\end{tikzcd}
\end{equation*}
is commutative. 
\end{theorem}

Before attempting to prove the above result, we introduce some notation. From now on, fix $\nu=\lambda-\mu=\sum_{i\in I}m_i\alpha_i$ and write $d=\rht(\nu)$. Recall from Appendix \ref{sec:setofparameters} that $\GTw_\nu$ denotes the collection of integral GT-weights of height $\nu$ and that $\GTwR_\nu$ (respectively $\GTwS_\nu$) denotes the regular (respectively singular) GT-weights of height $\nu$. For each GT-weight, we choose once and for all a resolution of multiplicities as in Section \ref{subsec:parityKLRW}. This data produces a map 
\begin{equation*}
\pi:\Lambda_\nu\to \Seq(\nu)
\end{equation*}
where $\pi(\bS)$ is the sequence obtained from the resolution of multiplicities of $\bS$.

Fix $\lambda\in P_+$ along with an integral set of parameters $\bR=(R_i)_{i\in I}$ of level $\lambda$. From this data, we obtain a subset $\GTw_\nu^\bR$ of GT-weights which produce non-zero idempotents in the algebra $P^\bR_\nu$. The restriction of $\pi$ to this subset is denoted $\pi_\bR$. Given a sequence $\bfi=(i_1,\dots, i_d)\in \Seq(\nu)$, let $\GTw_\nu^\bR(\bfi)=\pi^{-1}_\bR(\bfi)\subset \GTw_\nu^\bR$.

We partition the set $\Lambda_\nu^\bR$ into three disjoint subsets. Doing so, it is possible to rewrite the character $\chi(M)$ of an object $M$ of $\cO_{\mu}^\lambda(\bR)$ as a sum of three expressions. In the proof of Theorem \ref{thm:DbarMvsachi}, we will show that one can use this decomposition to produce \eqref{eq:equationDbarMvsachi}. Note that the techniques used here seem to be related to the ones appearing in the proof of \cite[Theorem~11.4]{baumann2021mirkovic}, which are ultimately related to \cite[Section~3]{knutson1999littelmann}. 

\begin{example}
In Example \ref{ex:calculCaracteresl3}, we saw that 
\begin{equation*}
\chi(M)=\dfrac{e^{-\alpha_1-\alpha_2}}{(1-e^{-\alpha_1})(1-e^{-\alpha_1-\alpha_2})}+1+\dfrac{e^{-\alpha_1}}{1-e^{-\alpha_1}}.
\end{equation*}
Upon computing the asymptotic character, only the first term of the above sum contributes and gives $\achi(M)=\tfrac{1}{\alpha_1(\alpha_1+\alpha_2)}$, which is equal to 
\begin{equation*}
\varsigma(\chbar(\Theta_{\cyc}(M)))=\varsigma(\chbar(L_2))=\varsigma(\Dbar^-_{(1,2)^\varsigma})=\Dbar^+_{(1,2)},
\end{equation*}
according to the computations made in Example \ref{ex:KLRWsl3} along with \eqref{eq:DbarPM_vs_varsigma}.
\end{example}

\subsubsection{Partitioning GT-weights}

Consider the partitioning of $\Lambda_\nu^\bR$ into the following three disjoint subsets:
\begin{enumerate}
\item $\GTw_\nu^{\cyc}:=\GTw_\nu^\bR\cap\GTw_\nu^{\reg}\cap \{\bS\in \GTw_\nu\;;\; e(\bS)e_{\cyc}\neq 0\}$,
\item $\GTw_\nu^\bR\cap\GTw_\nu^{\reg}\cap \{\bS\in \GTw_\nu\;;\; e(\bS)e_{\cyc}=0\}$ and
\item $\GTw_\nu^\bR\cap \GTw_\nu^{\sing}$.
\end{enumerate}
The proof of Theorem \ref{thm:DbarMvsachi} will be split into three lemmas according to the three subsets above.

Fix $\bfi=(i_1,\dots, i_d)\in \Seq(\nu)$. We start by parametrizing GT-weights of 
\begin{equation*}
\GTw_\nu^{\cyc}(\bfi):=\GTw_\nu^{\cyc}\cap \GTw_\nu^{\bR}(\bfi)
\end{equation*}
using a subset of $\Z^d$. For $i\in I$, let 
\begin{equation*}
\pi_i:=\begin{cases}
0 & i\equiv_2 0\\
1 & i\equiv_2 1
\end{cases}
\end{equation*}
and consider the map
\begin{equation}\label{eq:labellingQbfi}
Q_{\bfi}:\{(q_1,\dots, q_d)\in \Z^d\;;\; q_1<q_2<\dots <q_d\}\to \textstyle\prod_{i\in I} \C^{m_i}/\mathfrak{S}_{m_i}
\end{equation}
such that $Q_{\bfi}(q_1,\dots,q_d)=(Q_{\bfi}(q_1,\dots,q_d)_i)_{i\in I}$ is the collection of multisets such that $Q_{\bfi}(q_1,\dots,q_d)_i$ contains the integers $2q_j-\pi_{i_j}$ with $i_j=i$. By construction, $Q_{\bfi}(q_1,\dots,q_d)$ is an integral, regular set of parameters of height $\nu$ which satisfies 
\begin{equation*}
\pi(Q_{\bfi}(q_1,\dots,q_d))=\bfi.
\end{equation*}
Notice that, from the construction of idempotents in terms of GT-weights, $\GTw_\nu^{\cyc}(\bfi)$ is the image under $Q_{\bfi}$ of a subset of the form 
\begin{equation}\label{eq:preimage_GTcyci}
\{(q_1,\dots, q_d)\in \Z^d\;;\; q_1<q_2<\dots <q_d<N\}
\end{equation}
for some integer $N$ which depends on $\bR$. Now, set 
\begin{equation*}
h=\textstyle-\tfrac{1}{2}(\sum_{i\equiv_2 1}m_i\alpha_{i})+\tfrac{1}{2}(\lambda-\mu)-\tfrac{1}{2}(\sum_{i\in I}e_1(R_i)\varpi_i)\in \fh^\ast
\end{equation*}
and notice that $h$ depends only on $\lambda$, $\mu$ and $\bR$. Using Lemma \ref{lem:weightOfbS}, it follows that 
\begin{equation}\label{eq:weightQbfi}
\begin{aligned}
\textstyle\lambda(Q_{\bfi}(q_1,\dots,q_d))&=\textstyle\Big(\sum_{j=1}^d (q_j-\tfrac{1}{2}\pi_{i_j})\alpha_{i_j}\Big)+\tfrac{1}{2}(\lambda-\mu)-\tfrac{1}{2}\Big(\sum_{i\in I}e_1(R_i)\varpi_i\Big)\\
&=\textstyle\Big(\sum_{j=1}^d q_j\alpha_{i_j}\Big)+h.
\end{aligned}
\end{equation}

\begin{lemma}\label{lem:DbariFromChar}
There is an equality 
\begin{equation*}
\lim_{n\to\infty} \dfrac{1}{n^d} \sum_{\bS\in \GTw_\nu^{\cyc}(\bfi)} e^{\lambda(\bS)/n} =\DbarP_{\bfi}
\end{equation*}
as rational functions on $\fh$. 
\end{lemma}

\begin{proof}
By \eqref{eq:preimage_GTcyci}, we have
\begin{align*}
\sum_{\bS\in \GTw_\nu^{\cyc}(\bfi)} e^{\lambda(\bS)}&=\sum_{q_1<\dots<q_d<N} e^{q_1\alpha_{i_1}+\dots+q_d\alpha_{i_d}+h}=e^{\zeta}\sum_{q_1<\dots<q_d<0} e^{q_1\alpha_{i_1}+\dots+q_d\alpha_{i_d}}
\end{align*}
for some $\zeta\in \fh^\ast$. Applying Lemma \ref{lem:DbarPrep} shows that 
\begin{equation*}
\sum_{\bS\in \GTw_\nu^{\cyc}(\bfi)} e^{\lambda(\bS)}=\dfrac{e^{\xi}}{\prod_{k=1}^d (1-e^{-\beta^{\bfi}_k})}
\end{equation*}
for some $\xi\in \fh^\ast$. Computing the desired limit using a Taylor expansion yields
\begin{equation*}
\lim_{n\to \infty}\dfrac{1}{n^d}\dfrac{e^{\xi/n}}{\prod_{k=1}^d (1-e^{-\beta^{\bfi}_k/n})}=\prod_{k=1}^d \dfrac{1}{\beta^{\bfi}_k}=\DbarP_{\bfi}
\end{equation*}
as required.
\end{proof}

Now, consider the collection of GT-weights $E:=\GTw_\nu^\bR\cap\GTw_\nu^{\reg}\cap \{\bS\in \GTw_\nu\;;\; e(\bS)e_{\cyc}=0\}$. For each GT-weight $\bS\in E$, there exists an integer $1\leq k<d$ such that the $k$ first strands of the idempotent $e(\bS)$ are to the left of the leftmost red strand. Diagrammatically, this idempotent looks like
\begin{equation*}
e(\bS)=
\raisebox{-2.5em}{\begin{tikzpicture}
\blackstrandLabelled{0}{i_1}
\diagdots{0.5}
\blackstrandLabelled{1}{i_{k}}
\redstrand{1.5}
\blackstrandLabelled{2}{i_{k+1}}
\diagdots{2.5}
\blackstrandLabelled{3}{i_{d}}
\diagdots{3.5}
\end{tikzpicture}}
\end{equation*}
where the red strand appearing above is the leftmost one. Fix a sequence $\bfi\in \Seq(\nu)$ and set $E(\bfi):=E\cap \GTw_\nu^\bR(\bfi)$. This collection of GT-weights decomposes as 
\begin{equation*}
E(\bfi)=E(\bfi,1)\sqcup\dots\sqcup E(\bfi,d-1)
\end{equation*}
according to how many black strands are left of the leftmost red strand as above. Given $k$ such that $1\leq k<d$, the GT-weights in $E(\bfi,k)$ can be parametrized 
using $Q_\bfi$ by a set
\begin{equation*}
\{(q_1,\dots,q_d)\in \Z^d \;;\;q_1<\dots<q_k<N\leq q_{k+1}<\dots<q_d\text{ and }(q_{k+1},\dots,q_d)\in X\}
\end{equation*}
where $N\in \Z$ is an integer (which depends on $\bR$) and $X\subset \Z^{d-k}$ is a finite subset (which also depends on $\bR)$. 

\begin{lemma}\label{lem:vanishingFprime}
For $k=1,\dots,d$, there is an equality 
\begin{equation*}
\lim_{n\to\infty} \dfrac{1}{n^d} \sum_{\bS\in E(\bfi,k)} e^{\lambda(\bS)/n} =0
\end{equation*}
as rational functions on $\fh$. 
\end{lemma}

\begin{proof}
We use the same idea as in the proof Lemma \ref{lem:DbariFromChar}. First off, we have 
\begin{equation*}
\sum_{\bS\in E(\bfi,k)} e^{\lambda(\bS)}=\sum_{\substack{q_1<\dots<q_k<N\\ (q_{k+1},\dots q_d)\in X}} e^{q_1\alpha_{i_1}+\dots+q_d\alpha_{i_d}+h}=e^{\zeta}\sum_{q_1<\dots<q_k<0} e^{q_1\alpha_{i_1}+\dots+q_k\alpha_{i_k}}
\end{equation*}
for some $\zeta\in \fh^\ast$ from which it follows that
\begin{equation*}
\sum_{\bS\in E(\bfi,k)} e^{\lambda(\bS)}=\dfrac{e^\xi}{\prod_{t=1}^k (1-e^{-\beta^{\bfi}_t})}
\end{equation*}
some $\xi\in \fh^\ast$. Hence, using basic properties of limits, we have
\begin{equation*}
\lim_{n\to\infty} \dfrac{1}{n^d} \sum_{\bS\in E(\bfi,k)} e^{\lambda(\bS)/n}=\lim_{n\to \infty} \dfrac{1}{n^{d-k}} \dfrac{1}{n^k}\dfrac{e^{\xi/n}}{\prod_{t=1}^k (1-e^{-\beta^{\bfi}_t/n})}=\lim_{n\to \infty} \dfrac{1}{n^{d-k}} \prod_{t=1}^k \dfrac{1}{\beta^{\bfi}_t}=0
\end{equation*}
since $d-k>0$.
\end{proof}

Finally, consider the collection of GT-weights $F:=\GTw_\nu^\bR\cap \GTw_\nu^{\sing}$ and let $F(\bfi)=F\cap\GTw_\nu^\bR(\bfi)$. For each $\bS\in F(\bfi)$, there exists $1\leq r<d$ such that the integers appearing in GT-weight $\bS$ can be ordered as $s_1<\dots <s_{r}$ where $s_{j}$ appears with multiplicity $\rho_{j}\in\Z_{\geq 1}$. Given $\rho=(\rho_{1},\dots, \rho_{r})$, let $F(\bfi,\rho)\subset F(\bfi)$ be the subset of GT-weights whose multiplicities are given by $\rho$. If $s_j\in \bS_{i_j}$, then $\pi_\bR(\bS)$ is the sequence given by
\begin{equation*}
\bfi=(i_1^{\rho_1},\dots,i_q^{\rho_r})
\end{equation*}
where the exponent means repetition. Associated to such a sequence, there is a corresponding sequence of elements of $Q_+$ given by $(\rho_1\alpha_{i_1},\dots, \rho_r\alpha_{i_r})$.

As above, the subset $F(\bfi,\rho)$ can be parametrized using integer points. Consider the map 
\begin{equation*}
Q_{\bfi,\rho}:\{(q_1,\dots,q_r)\in \Z^r\;;\; q_1<\dots<q_r\}\to \textstyle\prod_{i\in I} \C^{m_i}/\mathfrak{S}_{m_i}
\end{equation*}
such that $Q_{\bfi}(q_1,\dots,q_r)=(Q_{\bfi}(q_1,\dots,q_r)_i)_{i\in I}$ is the collection of multisets such that $Q_{\bfi}(q_1,\dots,q_r)_i$ contains the integers $2q_j-\pi_{i_j}$ with $i_j=i$ each with multiplicity $\rho_j$. Analogously to \eqref{eq:weightQbfi}, we have
\begin{equation*}
\textstyle\lambda(Q_{\bfi,\rho}(q_1,\dots,q_r))=(\sum_{j=1}^d q_j\rho_j\alpha_{i_j})+h.
\end{equation*}
Notice that for every $N\in \Z$, Lemma \ref{lem:DbarPrep} shows that
\begin{equation*}
\sum_{q_1<\dots<q_r<N} e^{\lambda(Q_{\bfi,\rho}(q_1,\dots,q_r))}=\dfrac{e^\zeta}{\prod_{t=1}^r (1-e^{-(\rho_1\alpha_{i_1}+\dots+\rho_t\alpha_{i_t})})}
\end{equation*}
for some $\zeta\in \fh^\ast$. Since $r<d$, the same proofs as Lemma \ref{lem:DbariFromChar} and \ref{lem:vanishingFprime} yield:

\begin{lemma}\label{lem:vanishingFprimeprime}
There is an equality 
\begin{equation*}
\lim_{n\to\infty} \dfrac{1}{n^d} \sum_{\bS\in F(\bfi,\rho)} e^{\lambda(\bS)/n} =0
\end{equation*}
as rational functions on $\fh$. 
\end{lemma}

\begin{proof}[Proof of Theorem \ref{thm:DbarMvsachi}]
Fix $M$ an object of $\cO_\mu^\lambda(\bR)$ and consider the associated $P^\bR_\mu$-module $N:=\Theta(M)$ along with the $R_\nu^\lambda$-module $L:=e_{\cyc}N$. Using Lemma \ref{lem:weightOfbS}, it is possible to rewrite the character of $M$ as
\begin{equation*}
\chi(M)=\sum_{\bS\in\Lambda_\nu^{\bR}} \dim W_{\bS}(M)\,e^{\lambda(\bS)}=\sum_{\bS\in\Lambda_\nu^{\bR}}\dfrac{1}{\sigma(\bS)} \dim e(\bS)N\,e^{\lambda(\bS)}.
\end{equation*}
Now, splitting the sum according to the previous partitioning of $\Lambda_\nu^{\bR}$, we get
\begin{equation*}
\scalebox{0.9}{$
\displaystyle\chi(M)=\Big(\sum_{\bfi\in \Seq(\nu)} \sum_{\GTw_\nu^{\cyc}(\bfi)} \dim e(\bS)N\,e^{\lambda(\bS)} \Big) +\sum_{\bS\in E} \dim e(\bS)N\,e^{\lambda(\bS)}+\sum_{\bS\in F}\dfrac{1}{\sigma(\bS)} \dim e(\bS)N\,e^{\lambda(\bS)}.
$}
\end{equation*}
For $\bS\in \GTw_\nu^{\cyc}(\bfi)$, applying Lemma \ref{lem:adjonctionIdempotents} (with $A=P^{\bR}_{\mu}$, $e=e_{\cyc}$ and $e'=e(\bS)$) gives
\begin{equation*}
\dim e(\bS)N=\dim e(\bfi)L
\end{equation*}
which shows that $\dim e(\bS)N$ is independent of the choice of $\bS\in \GTw_\nu^{\cyc}(\bfi)$. It follows that the first sum can be rewritten as
\begin{equation}\label{eq:big_sum_cyc}
\sum_{\bfi\in \Seq(\nu)} \sum_{\bS\in \GTw_\nu^{\cyc}(\bfi)} \dim e(\bS)N\,e^{\lambda(\bS)}=\sum_{\bfi\in \Seq(\nu)}  \dim e(\bfi)L\,\big(\sum_{\bS\in \GTw_\nu^{\cyc}(\bfi)} e^{\lambda(\bS)}\big).
\end{equation}

Consider the equivalence relation on $\GTw_\nu$ defined by $\bS_1\sim \bS_2$ if $P^\bR_\mu e(\bS_1)\simeq P^\bR_\mu e(\bS_2)$. Notice that the quantity $\dim e(\bS)N$ does not change over each equivalence class. As $P^{\bR}_\mu$ is a finite-dimensional algebra, the number of equivalence classes under $\sim$ is finite. Suppose that there are $p\geq 1$ such classes. For $s\in \{1,\dots,p\}$, let $e_s:=\dim e(\bS)N$.

Now, for the second sum, recall that the set $E$ can be decomposed as the disjoint union 
\begin{equation*}
E=\sqcup_{\bfi\in \Seq(\nu)}\sqcup_{k=1}^{d-1} E(\bfi,k).
\end{equation*}
For each $\bfi\in \Seq(\nu)$ and $k$, $1\leq k\leq d-1$, the set $E(\bfi,k)$ decomposes further into $p$ disjoint subsets
\begin{equation*}
E(\bfi,k)=\sqcup_{s=1}^p E(\bfi,k)_s
\end{equation*}
according to the equivalence relation $\sim$.
Thus, we have 
\begin{equation}\label{eq:big_sum_E}
\sum_{\bS\in E} \dim e(\bS)N\,e^{\lambda(\bS)}=\sum_{k=1}^{d-1}\sum_{\bfi\in \Seq(\nu)}\sum_{s=1}^p  \sum_{\bS\in E(\bfi,k)_s} e_s e^{\lambda(\bS)}.
\end{equation}
Finally, for the third sum, we can proceed in a similar way and obtain
\begin{equation}\label{eq:big_sum_F}
\sum_{\bS\in F}\dfrac{1}{\sigma(\bS)} \dim e(\bS)N\,e^{\lambda(\bS)}=\sum_{r=1}^{d-1} \sum_{\bfi\in \Seq(\nu)}\sum_{\rho}\sum_{s=1}^p \sum_{\bS\in F(\bfi,\rho)_s} \dfrac{1}{\rho_{i_1}!\dots\rho_{i_r}!} e_s e^{\lambda(\bS)}
\end{equation}
where $F(\bfi,\rho)=\sqcup_{s=1}^p F(\bfi,\rho)_s$ is the decomposition into equivalence classes according to $\sim$. 

Since limits commute with finite sums, applying the preparatory Lemmas \ref{lem:vanishingFprime} and \ref{lem:vanishingFprimeprime} show that \eqref{eq:big_sum_E} and \eqref{eq:big_sum_F} vanish upon computing the asymptotic character. Combining this with \eqref{eq:big_sum_cyc} and the result of Lemma \ref{lem:DbariFromChar}, we get
\begin{equation*}
\achi_d(M)=\sum_{\bfi\in \Seq(\nu)}  \dim e(\bfi) (\Theta_{\cyc}(M)) \DbarP_{\bfi}.
\end{equation*}
Finally, as $\varsigma(\DbarP_{\bfi})=\DbarM_{\bfi^{\rev}}$, applying $\varsigma$ to the above equation yields 
\begin{equation*}
\varsigma(\achi_d(M))=\sum_{\bfi\in \Seq(\nu)}  \dim e(\bfi) (\Theta_{\cyc}(M)) \DbarM_{\bfi^{\rev}}=\chbar(\Theta_{\cyc}(M))
\end{equation*}
as required. 
\end{proof}

\section{Stitching together}\label{sec:stitching}

This section combines results of Sections \ref{sec:Yangian} and \ref{sec:KLR} into the main result of the paper and discusses consequences of it. A few conjectures about the relationship between the tensor structure on the category $\cOsh$ of \cite[Section~3.3]{hernandez2021shifted} and the monoidal structure on $\bigoplus_{\nu\in Q_+} R_{\nu}\fmodu$ are stated. A comment concerning a conjecture of \cite{braverman2016coulomb} is also added.

\subsection{The main theorem}

Let $\lambda\in P_+$, $\mu\in \wt(V(\lambda))$ and choose $\bR$ an integral set of parameters of level $\lambda$. Combining Theorems \ref{thm:epsilonTversusachi} and \ref{thm:DbarMvsachi} yields: 

\begin{theorem}\label{thm:theMainTheorem}
The diagram 
\begin{equation}\label{eq:bigDiag}
\begin{tikzcd}[column sep=1.5em,row sep=1em]
\HH_{\tp}\big((\cWbar{}^\lambda_\mu)_-\big)\arrow[d,"\psi_\lambda"']  & K_0(\cO_\mu^\lambda(\bR)) \arrow[l,"\CC_{\leq \tp}"'] \arrow[r,"\Theta_{\cyc}"]\arrow[dd,"\varsigma\circ\achi_d"'] & K_0(R_{\lambda-\mu}^\lambda\modu)\arrow[d,"\gamma"]\\
\C[N]_{-(\lambda-\mu)}\arrow[dr,"\DbarM"'] & & \C[N]_{-(\lambda-\mu)}\arrow[dl,"\DbarM"]\\
& \C(\fh) & 
\end{tikzcd}
\end{equation}
is commutative.
\end{theorem}

As shown by Lemma \ref{lem:CCpassingToTop} and Proposition \ref{prop:kernelThetacyc}, diagram \eqref{eq:bigDiag} is compatible with passing to $(\cO_\mu^\lambda(\bR))_{\tp}$. We expect that the following strengthening holds: 

\begin{conjecture}\label{conj:changementDeBase}
The diagram 
\begin{equation}\label{eq:bigDiagConjecture}
\begin{tikzcd}[column sep=1.5em]
\HH_{\tp}\big((\cWbar{}^\lambda_\mu)_-\big)\arrow[dr,"\psi_\lambda"'] & K_0(\cO_\mu^\lambda(\bR)_{\tp}) \arrow[l,"\CC_{\tp}"'] \arrow[r,"\Theta_{\cyc}"] & K_0(R_{\lambda-\mu}^\lambda\modu)\arrow[dl,"\gamma"]\\
& \C[N]_{-(\lambda-\mu)} & 
\end{tikzcd}
\end{equation}
is commutative.
\end{conjecture}

\begin{remark}
Since $\DbarM$ is not injective whenever $|\Phi_+|>2\rk \fg$ (see \cite[Remark~8.9]{baumann2021mirkovic}), Conjecture \ref{conj:changementDeBase} does not follows from Theorem \ref{thm:theMainTheorem}.
\end{remark}

Conjecture \ref{conj:changementDeBase} implies that the change-of-basis matrix from the dual canonical basis to the MV basis is provided by $\CC_{\tp}$ which has non-negative integral coefficients. As explained in the introduction, this conjecture is compatible with conjectures coming from cluster algebras.

\subsection{Multiplicative structure and monoidality of \texorpdfstring{$\Theta$}{the equivalence Theta}}

Shifted Yangians admit a family of algebra homomorphism called \textit{shifted coproducts} which were first defined in \cite[Section~4.11]{finkelberg2018comultiplication} and which generalize the usual Hopf algebra coproduct for unshifted Yangians. Given two coweights $\mu_1,\mu_2\in P^\vee$, the corresponding shifted coproduct is an algebra map
\begin{equation}\label{eq:shifted_coproducts}
\Delta_{\mu_1,\mu_2}:Y_{\mu_1+\mu_2}\to Y_{\mu_1}\otimes Y_{\mu_2}.
\end{equation}
These homomorphisms endow the category $\cOsh:=\bigoplus_{\mu \in P^\vee} \cO_\mu$ with a tensor structure, making $K_0(\cOsh)$ into a commutative ring by \cite[Theorem~3.14]{hernandez2021shifted}. Under this tensor structure, the character map $\chi:K_0(\cOsh)\to \charRing$ (see Section \ref{subsec:characters}) becomes an algebra homomorphism (as a corollary of \cite[Theorem~3.14]{hernandez2021shifted}).\par

\begin{lemma}\label{lem:multiplicativeachi}
Let $\lambda_1,\lambda_2\in P^\vee_+$, $\mu_1,\mu_2\in P^\vee$ and for $i=1,2$ set $d_i=\rht(\lambda_i-\mu_i)$ and let $d=d_1+d_2$. Consider $\bR_1$, $\bR_2$ integral sets of parameters of level $\lambda_1$, $\lambda_2$ respectively. Then, if $M_i\in \ob(\cO_{\mu_i}^{\lambda_i}(\bR_i))$ for $i=1,2$, one has
\begin{equation*}
\achi_d(M_1\otimes M_2)=\achi_{d_1}(M_1)\cdot\achi_{d_2}(M_2).
\end{equation*}
\end{lemma}

\begin{proof}
Follows immediately from the fact that $\chi$ respect the multiplicative structure and from properties of limits. 
\end{proof}

Let $R=\bigoplus_{\nu\in Q_+} R_\nu$ be the KLR algebra associated to $\fg$ and consider the category $\bigoplus_{\nu\in Q_+} R_\nu\fmodu $. As stated in Theorem \ref{thm:RmuisoCNmu}, there is an isomorphism of vector spaces
\begin{equation*}
\textstyle \bigoplus_{\nu\in Q_+} \C\otimes_{\Z} K_0 (R_\nu\fmodu )\simeq \C[N].
\end{equation*}
Using the non-unital algebra maps defined in \eqref{eq:iotaKLR}, as shown by \cite{khovanov2009diagrammatic,rouquier20082}, one can endow the category ${\bigoplus_{\nu\in Q_+} R_\nu\fmodu}$ with a monoidal structure often referred to as the \textit{convolution product}. This provides $\bigoplus_{\nu\in Q_+} K_0 (R_\nu\fmodu)$ with a commutative ring structure and upgrades the isomorphism of Theorem \ref{thm:RmuisoCNmu} into an algebra isomorphism. Furthermore, the character map \eqref{eq:characterKLR} also becomes a ring homomorphism and since the $\DbarPM$ satisfy the shuffle identities by \cite[Lemma~8.3]{baumann2021mirkovic}, it follows that the bar-character map also is an algebra homomorphism. 

The following result is a consequence of the previous discussion, of Theorem \ref{thm:DbarMvsachi} and of Lemma \ref{lem:multiplicativeachi}: 

\begin{corollary}\label{cor:chbar_monoidal}
Using the same notation as in Lemma \ref{lem:multiplicativeachi}, one has
\begin{equation*}
\chbar(\Theta_{\cyc}(M_1\otimes M_2))=\chbar(\Theta_{\cyc}(M_1))\cdot\chbar(\Theta_{\cyc}(M_2)).
\end{equation*}
\end{corollary}

For each $\bR$, consider the functor $\Theta_{\cyc,\bR}$ as in \eqref{eq:definitionThetaCyc}. In a follow-up paper \cite{kamnitzer2026kosh}, we show that the shifted coproducts \eqref{eq:shifted_coproducts} descend to truncations, in the sense that there is a map making the following diagram commutative:
\begin{equation*}
\adjustbox{scale=1,center}{
\begin{tikzcd}[row sep =2em]
Y_{\mu_1+\mu_2} \ar[r,"\Delta_{\mu_1,\mu_2}"] \ar[d,two heads] & Y_{\mu_1}\otimes Y_{\mu_2} \ar[d,two heads] \\
Y^{\lambda_1+\lambda_2}_{\mu_1+\mu_2}(\bR_1 \cup \bR_2) \ar[r,dashed] & Y_{\mu_1}^{\lambda_1}(\bR_1) \otimes Y_{\mu_2}^{\lambda_2}(\bR_2)
\end{tikzcd}
}
\end{equation*}
Our proof works under the assumption that the dimensions of the weight spaces $V(\lambda_1)_{\mu_1}$ and $V(\lambda_2)_{\mu_2}$ are both non-zero. Above, the union of sets of parameters is the point-wise union of multisets.

Corollary \ref{cor:chbar_monoidal} provides evidence for the following conjecture: 

\begin{conjecture}\label{con:monoidalThetacyc}
For every pair of integral sets of parameters $\bR_1$, $\bR_2$, the functors 
\begin{equation*}
\Theta_{\cyc,\bR_1}(\trou)\circ\Theta_{\cyc,\bR_2}(\trou)\quad\text{and}\quad\Theta_{\cyc,\bR_1\cup\bR_2}(\trou\otimes\trou)
\end{equation*}
are naturally isomorphic.
\end{conjecture}

In a forthcoming paper \cite{kamnitzer2026kosh}, we show that Conjecture \ref{con:monoidalThetacyc} holds on the level of $K_0$. We expect that there is a generalization of this conjecture which would state that the functor $\Theta$ is monoidal. However, as there is no obvious monoidal structure on the category $\bigoplus_{\mu,\bR}P^\bR_\mu\modu$, we are not able to make a precise statement.

\subsection{On a conjecture of Nakajima}

Let $\fg$ be a symmetric Kac-Moody Lie algebra. In this setting, conjecturally, the replacement for the generalized affine Grassmannian slice $\cWbar^\lambda_\mu$ is the \textit{Coulomb branch} $\mathcal{M}_C(\mathsf{G},\mathsf{N})$ first defined in \cite{nakajima1994instantons} and studied in \cite{braverman2016coulomb}. Note that the data of the groups $\mathsf{G}$ and $\mathsf{N}$ depend on $\lambda$ and $\mu$. For $\fg$ of type $A,D,E$, \cite[Theorem~3.10]{braverman2016coulomb} shows that the spaces $\cWbar{}^\lambda_\mu$ and $\mathcal{M}_C(\mathsf{G},\mathsf{N})$ are isomorphic. As in the case of affine Grassmannian slices, the varieties $\mathcal{M}_C(\mathsf{G},\mathsf{N})$ are endowed with torus actions and it is conjectured \cite[Conjecture~3.25]{braverman2016coulomb} that under a well-chosen cocharacter $2\rho^\vee$, there is an isomorphism
\begin{equation*}
\C\otimes_{\Z}\HH_{\tp}\Big(\big(\mathcal{M}_C(\mathsf{G},\mathsf{N})\big)_-;\Z\Big)\simeq V(\lambda)_\mu
\end{equation*}
where, as previously, $(\trou)_-$ denotes the repelling set. This conjecture is of great interest since it should serve as a replacement for the geometric Satake isomorphism in the Kac-Moody setting. 

Results of \cite{kamnitzer2022lie} provide an extension in the Kac-Moody case of the main theorems of \cite{kamnitzer2019category} and generalize the functor $\Theta$. By carefully checking that the necessary hypotheses hold, it is possible to extend Theorem \ref{thm:theMainTheorem} to produce the commutative diagram 
\begin{equation*}
\begin{tikzcd}[row sep=1.25em, column sep=0.75em]
K_0(\cO_\mu^\lambda(\bR)) \arrow[r]\arrow[d]& K_0(R_{\lambda-\mu}^\lambda\modu)\arrow[d]\\
\HH_{\tp}\Big(\big(\mathcal{M}_C(\mathsf{G},\mathsf{N})\big)_-;\Z\Big) \arrow[r]& \C(\fh)
\end{tikzcd}
\end{equation*} 
where $\cO^\lambda_\mu(\bR)$ is the category $\cO$ for the Coulomb branch algebra $\mathcal{A}(\mathsf{G},\mathsf{N})$ (see \cite[Section~3(iv)]{braverman2016towards} for the definition of $\mathcal{A}(\mathsf{G},\mathsf{N})$). In particular, this commutative diagram implies that
\begin{equation*}
\dim \C\otimes_{\Z}\HH_{\tp}\Big(\big(\mathcal{M}_C(\mathsf{G},\mathsf{N})\big)_-;\Z\Big)\geq \dim \Ima\Big(\C\otimes_{\Z}K_0(R_{\lambda-\mu}^\lambda\modu)\to\C(\fh) \Big).
\end{equation*}
By the main result of \cite{kang2012categorification}, the following corollary is immediate and it provides evidence for \cite[Conjecture~3.25]{braverman2016coulomb}.
\begin{corollary}\label{cor:conjectureNakajima}
If $\mu\in \wt(V(\lambda))$, then $\HH_{\tp}\Big(\big(\mathcal{M}_C(\mathsf{G},\mathsf{N})\big)_-;\Z\Big)\neq \{0\}$.
\end{corollary}

\appendix
\section{Sets of parameters}\label{sec:setofparameters}
Suppose $\fg$ is simply-laced and adopt a parity convention according to Convention \ref{con:parity}.\par 

A \textit{set of parameters} is a collection of multi-sets $\bR=(R_i)_{i\in I}$, where each $R_i$ is a multi-set of elements of $\C$ which has finite cardinality. The \textit{cardinality vector} of a set of parameters is a tuple $|\bR|:=(n_1,\dots,n_r)\in \Z_{\geq 0}^I$ such that $|R_i|=n_i$. If $\lambda=\sum_{i\in I} \lambda_i\varpi_i\in P_+$, a set of parameters $\bR$ has \textit{level} $\lambda$ if $|\bR|=(\lambda_1,\dots, \lambda_r)$. If $\nu=\sum_{i\in I} m_i\alpha_i\in Q_+$, a set of parameters $\bR$ has \textit{height} $\nu$ if $|\bR|=(m_1,\dots,m_r)$. A set of parameters $\bR=(R_i)_{i\in I}$ is called \textit{integral} if, for each $i$, $R_i$ contains only integers having the same parity as $i$.\par 

For $\nu\in Q_+$, let $\GTw_{\nu}$ denote the collection of integral sets of parameters of height $\nu$. Elements of $\Lambda_{\nu}$ are called \textit{Gelfand-Tsetlin weights} (GT-weights) of height $\nu$. Moreover, let $\GTw=\sqcup_{\nu\in Q_+}\GTw_\nu$.

To each integral set of parameters $\bR$, one associates \textit{multiplicity functions} $\rho_i^{\bR}:\Z\to\Z_{\geq 0}$ defined by 
\begin{equation*}
\rho_i^{\bR}(q)=\begin{cases}
\text{multiplicity of }2q\text{ in }R_i&i\text{ even}\\
\text{multiplicity of }2q+1\text{ in }R_i &i\text{ odd}
\end{cases}.
\end{equation*}
Let 
\begin{equation*}
\textstyle\sigma(\bR)=\prod_{i\in I}\prod_{q\in \Z} \rho_i^{\bR}(q)!
\end{equation*}
and call $\bR$ \textit{singular} if $\sigma(\bR)>1$ and call $\bR$ \textit{regular} otherwise. For each set of parameters $\bR$ of cardinality vector $(n_1,\dots,n_r)$, choose $\nu(\bR)\in \prod_{i\in I} \C^{n_i}$ such that $\nu_i$ has entries given by a choice of ordering of the elements of $R_i$. The quantity $\sigma(\bR)$ is the order of the subgroup of $\prod_{i\in I}\mathfrak{S}_{n_i}$ which stabilizes $\nu(\bR)$. 

The set of regular (respectively singular) GT-weights is denoted by $\GTw^{\reg}$ (respectively $\GTw^{\sing}$).

The value of the $k$-th elementary symmetric function $e_k$ on $R_i$ is well defined and it is denoted by $e_k(R_i)\in \C$. As usual, when $k>|R_i|$, we use the convention that $e_k(R_i)=0$.\par 

Let $\bR$ be an integral set of parameters such that $|R_i|=\lambda_i\in\Z_{\geq 0}$. Then, for each $q$, define the list of fundamental weights	
\begin{equation*}
\uvarpi_{\bR}(q)=(\varpi_{i_1}^{\rho_{i_1}^{\bR}(q)},\dots,\varpi_{i_r}^{\rho_{i_r}^{\bR}(q)})
\end{equation*}
where $\varpi_i^n$ means that $\varpi_i$ is repeated $n$ times in the list. Define $\uvarpi_{\bR}=\sqcup_{q\in\Z}\uvarpi_{\bR}(q)$, where $\sqcup$ denotes concatenation of the lists from left to right along the natural order on $\Z$. The list $\uvarpi_{\bR}$ is finite and if $\uvarpi_{\bR}=(\varpi_{j_1},\dots,\varpi_{j_\ell})$, then $\sum_{k=1}^\ell \varpi_{j_k}=\sum_{i\in I} \lambda_i\varpi_i$.

\begin{example}
Let 
\begin{equation*}
\Gamma:1\leftarrow 2\rightarrow 3
\end{equation*}
with $I_0=\{2\}$, $I_1=\{1,3\}$ and choose the total order $I=\{2<1<3\}$. The set of parameters 
\begin{equation*}
\bR=(R_1=\{-1,1,1\},R_2=\{0\},R_3=\emptyset)
\end{equation*}
has level $3\varpi_1+\varpi_2$ and it is singular since $\sigma(\bR)=2$. Moreover, one can check that
\begin{equation*}
\uvarpi_{\bR}=(\varpi_1,\varpi_2,\varpi_1,\varpi_1),
\end{equation*}
which is a list that sums to $3\varpi_1+\varpi_2$.
\end{example}

\section{The product monomial crystal}\label{sec:prodMonomialCrystal}

The monomial crystal $\cB$ was introduced in \cite{nakajima2003tanalogues} and it was shown in \cite{kashiwara2003realizations} that it is a normal $\fg$-crystal. In \cite{kamnitzer2019highest}, it was proven that when $\fg$ is simply-laced, certain subcrystals of $\cB$ parametrize highest weights of modules in category $\cO$ for truncated shifted Yangians. Subsequently, in \cite{kamnitzer2019category}, a connection between KLRW algebras and the product monomial crystal was made. The combinatorial properties of the product monomial crystal have been studied in \cite{gibson2021demazure}. As it will be a useful tool in this paper, we provide the reader with a quick review of the definition of the product monomial crystal in the simply-laced case.

For $\lambda\in P_+$, let $\cB(\lambda)$ be the (unique) normal $\fg$-crystal of highest weight $\lambda$. 

Suppose $\fg$ is simply-laced and adopt a parity convention according to Convention \ref{con:parity}. Consider the variables $\{y_{i,k}\;;\; i\in I,k\in \Z, k\equiv_2 i\}$ and let $\cB$ be the set of all monomials in the variables $y_{i,k}^{\pm1}$. We refer to parameter $k$ in a variable $y_{i,k}$ as its spectral parameter. Let  
\begin{equation*}
z_{i,k}=\dfrac{y_{i,k}y_{i,k+2}}{\prod_{j\sim i} y_{i,k+1}}.
\end{equation*}
Endow $\cB$ with maps $(\wt,\epsilon_i,\phi_i,\tilde{e}_i,\tilde{f}_i)_{i\in I}$ defined as follows. Set $\wt(y_{i,k}):=\varpi_i$ and for $m,m'\in \cB$, define $\wt(mm'):=\wt(m)+\wt(m')$. Fix a monomial $m=\prod_{i\in I,k\in \Z} y_{i,k}^{a_{i,k}}\in \cB$. For each $k\in \Z$ and $i\in I$, define
\begin{align*}
\epsilon_i^{k}(m)&=-\textstyle\sum_{\ell\leq k} a_{i,\ell}, & \phi_i^{k}(m)&=\textstyle\sum_{\ell\geq k } a_{i,\ell},\\
\epsilon_i(m)&=\max \{\epsilon_i^{k}(m) \;;\;k\in \Z\}, & \phi(m)&=\max \{\phi^{k}(m) \;;\;k\in \Z\}.
\end{align*}
Finally, define the following Kashiwara operators: 
\begin{align*}
\tilde{e}_i(m)&=\begin{cases}
0 & {\scriptstyle\text{if }\epsilon_i(m)=0}\\
m \cdot z_{i,k}\phantom{^{-1}} & {\scriptstyle\text{otherwise, where $k$ is the largest integer such that $i\equiv_2 k$ and $\epsilon_i(m)=\epsilon_i^k(m)$}}
\end{cases}\\
\tilde{f}_i(m)&=\begin{cases}
0 & {\scriptstyle\text{if }\phi_i(m)=0}\\
m \cdot z_{i,k-2}^{-1} & {\scriptstyle\text{otherwise, where $k$ is the smallest integer such that $i\equiv_2 k$ and $\phi_i(m)=\phi_i^k(m)$}}
\end{cases}
\end{align*}

\begin{theorem}[{\cite[Proposition~3.1]{kashiwara2003realizations}}]
The set $\cB$ together with maps $(\wt,\epsilon_i,\phi_i,\tilde{e}_i,\tilde{f}_i)_{i\in I}$ defined above is a normal crystal. 
\end{theorem}

For $i\in I$, $k\in \Z$ with $k\equiv_2 i$, let $\cB(\varpi_i,k)$ be the subcrystal of $\cB$ generated by $y_{i,k}$. As an abstract crystal, it is isomorphic to $\cB(\varpi_i)$. Now, for $c\in \C$, let $\cB(\varpi_i,c)$ be the crystal which has the same crystal operators as $\cB(\varpi_i,k)$ and whose elements are monomials of $\cB(\varpi_i,k)$ for which the spectral parameters are translated according to $k\mapsto c-k$. 

For $\lambda=\sum_{i\in I}\lambda_i\varpi_i\in P_+$ and $\bR=(\bR_i)_{i\in I}$ a set of parameters of level $\lambda$, define 
\begin{equation*}
\cB(\lambda,\bR)=\textstyle\prod_{i\in I}\prod_{c\in \bR_i} \cB(\varpi_i,c)
\end{equation*}
where the products are products of monomials of $\cB$. When $\bR$ is an integral set of parameters, $\cB(\lambda,\bR)$ can be endowed with a crystal structure by restricting the crystal operators of $\cB$. This is the \textit{product monomial crystal} associated to $\bR$. By \cite[Theorem~2.2]{kamnitzer2019highest}, it is a normal subcrystal of $\cB$ such that 
\begin{equation*}
\textstyle\cB(\lambda)\subset \cB(\lambda,\bR)\subset \bigotimes_{i\in I}\cB(\varpi_i)^{\otimes \lambda_i}
\end{equation*}
as crystals. 

\begin{example}
Let $\fg=\sl_3$. Then, the crystal generated by $y_{1,c}$ is given by
\begin{equation*}
\cB(\varpi_1,c)=\Big\{\;\; y_{1,c}\xrightarrow[]{\tilde{f}_1}\dfrac{y_{2,c-1}}{y_{1,c-2}}\xrightarrow[]{\tilde{f}_2}\dfrac{1}{y_{2,c-3}} \;\;\Big\}.
\end{equation*}
Consider the set of parameters $\bR=(\bR_1=\{c\},\bR_2=\{c+k\})$ with $k\in 2\Z_{\geq 0}+1$. If $k=1$, then the variables appearing in the crystal $\cB(\varpi_1+\varpi_2,\bR)$ are
\begin{equation*}{\small
\begin{tikzpicture}
\def\ray{2.5}
\node  at (0+30:\ray) {$\dfrac{y_{1,c}^2}{y_{2,c-1}}$};
\node  at (60+30:\ray) {$y_{1,c}y_{2,c+1}$};
\node  at (120+30:\ray) {$\dfrac{y_{2,c-1}y_{2,c+1}}{y_{1,c-2}}$};
\node  at (180+30:\ray) {$\dfrac{y_{2,c-1}}{y_{1,c-2}^2}$};
\node  at (240+30:\ray) {$\dfrac{1}{y_{1,c-2}y_{2,c-3}}$};
\node  at (300+30:\ray) {$\dfrac{y_{1,c}}{y_{2,c-1}y_{2,c-3}}$};
\node  at (0,0) {$\dfrac{y_{1,c}}{y_{1,c-2}},\dfrac{y_{2,c+1}}{y_{2,c-3}}$};
\end{tikzpicture}}
\end{equation*}
and the crystal operators make $\cB(\varpi_1+\varpi_2,\bR)$ isomorphic to $\cB(\varpi_1+\varpi_2)$. On the other hand, when $k\geq 3$, $\cB(\varpi_1+\varpi_2,\bR)$ has $9$ monomials and, as a crystal, it is isomorphic to $\cB(\varpi_1+\varpi_2)\oplus \cB(0)$.
\end{example}

{
\emergencystretch=3em
\printbibliography

@article {anderson2003polytope,
    AUTHOR = {Anderson, Jared E.},
     TITLE = {A polytope calculus for semisimple groups},
   JOURNAL = {Duke Math. J.},
  FJOURNAL = {Duke Mathematical Journal},
    VOLUME = {116},
      YEAR = {2003},
    NUMBER = {3},
     PAGES = {567--588},
      ISSN = {0012-7094},
   MRCLASS = {20G05 (14L99)},
  MRNUMBER = {1958098},
MRREVIEWER = {Timo Neuvonen},
       DOI = {10.1215/S0012-7094-03-11636-1},
       URL = {https://doi.org/10.1215/S0012-7094-03-11636-1},
}

@book{assem2006elements,
    AUTHOR = {Assem, Ibrahim and Simson, Daniel and Skowro\'{n}ski, Andrzej},
     TITLE = {Elements of the representation theory of associative algebras.
              {V}ol. 1},
    SERIES = {London Mathematical Society Student Texts},
    VOLUME = {65},
      NOTE = {Techniques of representation theory},
 PUBLISHER = {Cambridge University Press, Cambridge},
      YEAR = {2006},
     PAGES = {x+458},
      ISBN = {978-0-521-58423-4; 978-0-521-58631-3; 0-521-58631-3},
   MRCLASS = {16G10 (16-02)},
  MRNUMBER = {2197389},
MRREVIEWER = {Peter W. Donovan},
       DOI = {10.1017/CBO9780511614309},
       URL = {https://doi.org/10.1017/CBO9780511614309},
}

@article{baumann2021mirkovic,
    AUTHOR = {Baumann, Pierre and Kamnitzer, Joel and Knutson, Allen},
     TITLE = {The {M}irkovi\'{c}-{V}ilonen basis and\linebreak {D}uistermaat-{H}eckman
              measures},
      NOTE = {With an appendix by Anne Dranowski, Kamnitzer and Calder
              Morton-Ferguson},
   JOURNAL = {Acta Math.},
  FJOURNAL = {Acta Mathematica},
    VOLUME = {227},
      YEAR = {2021},
    NUMBER = {1},
     PAGES = {1--101},
      ISSN = {0001-5962},
   MRCLASS = {22E57 (14D24 14M15 17B37 20G20)},
  MRNUMBER = {4346265},
MRREVIEWER = {Jacinta Torres},
       DOI = {10.4310/ACTA.2021.v227.n1.a1},
       URL = {https://doi.org/10.4310/ACTA.2021.v227.n1.a1},
}

@article {baumann2024bases,
    AUTHOR = {Baumann, Pierre and Gaussent, St\'{e}phane and Littelmann, Peter},
     TITLE = {Bases of tensor products and geometric {S}atake
              correspondence},
   JOURNAL = {J. Eur. Math. Soc. (JEMS)},
  FJOURNAL = {Journal of the European Mathematical Society (JEMS)},
    VOLUME = {26},
      YEAR = {2024},
    NUMBER = {3},
     PAGES = {919--983},
      ISSN = {1435-9855},
   MRCLASS = {22E46 (14M15)},
  MRNUMBER = {4721027},
       DOI = {10.4171/jems/1302},
       URL = {https://doi.org/10.4171/jems/1302},
}

@misc{beilinson1991quantization,
  title={Quantization of Hitchin’s integrable system and Hecke eigensheaves},
  author={Beilinson, Alexander and Drinfeld, Vladimir}
}

@article {berensetein2005cluster,
    AUTHOR = {Berenstein, Arkady and Fomin, Sergey and Zelevinsky, Andrei},
     TITLE = {Cluster algebras. {III}. {U}pper bounds and double {B}ruhat cells},
   JOURNAL = {Duke Math. J.},
  FJOURNAL = {Duke Mathematical Journal},
    VOLUME = {126},
      YEAR = {2005},
    NUMBER = {1},
     PAGES = {1--52},
      ISSN = {0012-7094},
   MRCLASS = {16S99 (05E15 14M17 22E46)},
  MRNUMBER = {2110627},
       DOI = {10.1215/S0012-7094-04-12611-9},
       URL = {https://doi.org/10.1215/S0012-7094-04-12611-9},
}

@article{berenstein2006geometric,
    AUTHOR = {Berenstein, Arkady and Kazhdan, David},
     TITLE = {Geometric and unipotent crystals. {II}. {F}rom unipotent
              bicrystals to crystal bases},
 BOOKTITLE = {Quantum groups},
    SERIES = {Contemp. Math.},
    VOLUME = {433},
     PAGES = {13--88},
 PUBLISHER = {Amer. Math. Soc., Providence, RI},
      YEAR = {2007},
   MRCLASS = {17B37 (14L30 22E47)},
  MRNUMBER = {2349617},
MRREVIEWER = {Jonathan Brundan},
       DOI = {10.1090/conm/433/08321},
       URL = {https://doi.org/10.1090/conm/433/08321},
}

@book{bourbaki1968groupes,
    AUTHOR = {Bourbaki, Nicolas},
     TITLE = {\'{E}l\'{e}ments de math\'{e}matique: groupes et alg\`ebres de {L}ie},
      NOTE = {Chapitres 7 et 8},
 PUBLISHER = {Masson, Paris},
      YEAR = {1975},
     PAGES = {265}
}

@article{braden2012quantizations,
    AUTHOR 		= {Braden, Tom and Proudfoot, Nicholas and Webster, Ben},
	TITLE 		= {Quantizations of conical symplectic resolutions {I}: local and global structure},
	JOURNAL 	= {Ast\'{e}risque},
	FJOURNAL 	= {Ast\'{e}risque},
    NUMBER 		= {384},
	YEAR 		= {2016},
	PAGES		= {1--73},
	ISSN 		= {0303-1179},
	ISBN 		= {978-2-85629-845-9},
	MRCLASS 	= {14F05 (16S38 17B10 53D50)},
	MRNUMBER 	= {3594664},
	MRREVIEWER 	= {Xiaobin Li},
	}

@article{braden2014quantizations,
    AUTHOR = {Braden, Tom and Licata, Anthony and Proudfoot, Nicholas and Webster, Ben},
     TITLE = {Quantizations of conical symplectic resolutions {II}: category {$\mathcal O$} and symplectic duality},
      NOTE = {With an appendix by I. Losev},
   JOURNAL = {Ast\'{e}risque},
  FJOURNAL = {Ast\'{e}risque},
    NUMBER = {384},
      YEAR = {2016},
     PAGES = {75--179},
      ISSN = {0303-1179},
      ISBN = {978-2-85629-845-9},
   MRCLASS = {53D50 (14F05 14J81 16S38 17B10)},
  MRNUMBER = {3594665},
}

@article{braverman2016coulomb,
	AUTHOR = {Braverman, Alexander and Finkelberg, Michael and Nakajima, Hiraku},
     TITLE = {Coulomb branches of {$3d$} {$\mathcal N=4$} quiver gauge theories and slices in the affine {G}rassmannian},
      NOTE = {With two appendices by Braverman, Finkelberg, Joel Kamnitzer, Ryosuke Kodera, Nakajima, Ben Webster and Alex Weekes},
   JOURNAL = {Adv. Theor. Math. Phys.},
  FJOURNAL = {Advances in Theoretical and Mathematical Physics},
    VOLUME = {23},
      YEAR = {2019},
    NUMBER = {1},
     PAGES = {75--166},
      ISSN = {1095-0761},
   MRCLASS = {57R57 (14M15 16G20 17B81 81T13)},
  MRNUMBER = {4020310},
       DOI = {10.4310/ATMP.2019.v23.n1.a3},
       URL = {https://doi.org/10.4310/ATMP.2019.v23.n1.a3},
}

@article{braverman2016towards,
  title={Towards a mathematical definition of Coulomb branches of {$3$}-dimensional {$\mathcal N=4$} gauge theories, II},
  author={Braverman, Alexander and Finkelberg, Michael and Nakajima, Hiraku},
  journal={arXiv preprint arXiv:1601.03586},
  year={2016}
}

@article{brion1997equivariant,
    AUTHOR = {Brion, M.},
     TITLE = {Equivariant {C}how groups for torus actions},
   JOURNAL = {Transform. Groups},
  FJOURNAL = {Transformation Groups},
    VOLUME = {2},
      YEAR = {1997},
    NUMBER = {3},
     PAGES = {225--267},
      ISSN = {1083-4362},
   MRCLASS = {14C15 (14L30 14M15)},
  MRNUMBER = {1466694},
MRREVIEWER = {Dan Edidin},
       DOI = {10.1007/BF01234659},
       URL = {https://doi.org/10.1007/BF01234659},
}

@article{brundan2006shifted,
    AUTHOR = {Brundan, Jonathan and Kleshchev, Alexander},
     TITLE = {Shifted {Y}angians and finite {$W$}-algebras},
   JOURNAL = {Adv. Math.},
  FJOURNAL = {Advances in Mathematics},
    VOLUME = {200},
      YEAR = {2006},
    NUMBER = {1},
     PAGES = {136--195},
      ISSN = {0001-8708},
   MRCLASS = {17B37},
  MRNUMBER = {2199632},
MRREVIEWER = {Chengming Bai},
       DOI = {10.1016/j.aim.2004.11.004},
       URL = {https://doi.org/10.1016/j.aim.2004.11.004},
}

@article {brundan2008highest,
    AUTHOR = {Brundan, Jonathan and Goodwin, Simon M. and Kleshchev,
              Alexander},
     TITLE = {Highest weight theory for finite {$W$}-algebras},
   JOURNAL = {Int. Math. Res. Not. IMRN},
  FJOURNAL = {International Mathematics Research Notices. IMRN},
      YEAR = {2008},
    NUMBER = {15},
     PAGES = {Art. ID rnn051, 53},
      ISSN = {1073-7928,1687-0247},
   MRCLASS = {17B10},
  MRNUMBER = {2438067},
MRREVIEWER = {Serge\ M.\ Skryabin},
       DOI = {10.1093/imrn/rnn051},
       URL = {https://doi.org/10.1093/imrn/rnn051},
}

@book{bump2017crystal,
    AUTHOR = {Bump, Daniel and Schilling, Anne},
     TITLE = {Crystal bases},
      NOTE = {Representations and combinatorics},
 PUBLISHER = {World Scientific Publishing Co. Pte. Ltd., Hackensack, NJ},
      YEAR = {2017},
     PAGES = {xii+279},
      ISBN = {978-981-4733-44-1},
   MRCLASS = {05-01 (05E10 14T05 17B10)},
  MRNUMBER = {3642318},
       DOI = {10.1142/9876},
       URL = {https://doi.org/10.1142/9876},
}

@article {casbi2021equivariant,
    AUTHOR = {Casbi, {\'E}lie},
     TITLE = {Equivariant multiplicities of simply-laced type flag minors},
   JOURNAL = {Represent. Theory},
  FJOURNAL = {Representation Theory. An Electronic Journal of the American
              Mathematical Society},
    VOLUME = {25},
      YEAR = {2021},
     PAGES = {1049--1092},
   MRCLASS = {16G10 (20G05)},
  MRNUMBER = {4353894},
MRREVIEWER = {Jian Min Chen},
       DOI = {10.1090/ert/589},
       URL = {https://doi.org/10.1090/ert/589},
}

@article {casbi2023quantum,
    AUTHOR = {Casbi, {\'E}lie and Li, Jian-Rong},
     TITLE = {Equivariant multiplicities via representations of quantum
              affine algebras},
   JOURNAL = {Selecta Math. (N.S.)},
  FJOURNAL = {Selecta Mathematica. New Series},
    VOLUME = {29},
      YEAR = {2023},
    NUMBER = {1},
     PAGES = {Paper No. 9, 58},
      ISSN = {1022-1824},
   MRCLASS = {17B37 (13F60 16G99 17B67)},
  MRNUMBER = {4514985},
MRREVIEWER = {Hiro-Fumi Yamada},
       DOI = {10.1007/s00029-022-00805-y},
       URL = {https://doi.org/10.1007/s00029-022-00805-y},
}

@article{casbi2024auslander,
  title={Auslander-Reiten combinatorics and {$q$}-characters of representations of affine quantum groups},
  author={Casbi, {\'E}lie and Li, Jian-Rong},
  journal={arXiv preprint arXiv:2410.06046},
  year={2024}
}

@book{chriss1997representation,
	AUTHOR = {Chriss, Neil and Ginzburg, Victor},
     TITLE = {Representation theory and complex geometry},
 PUBLISHER = {Birkh\"{a}user Boston, Inc., Boston, MA},
      YEAR = {1997},
     PAGES = {x+495},
      ISBN = {0-8176-3792-3},
   MRCLASS = {22E47 (14F99 19L47 20G05 22-02 58F05)},
  MRNUMBER = {1433132},
MRREVIEWER = {William M. McGovern},
}

@article {dupont2011generic,
    AUTHOR = {Dupont, G.},
     TITLE = {Generic variables in acyclic cluster algebras},
   JOURNAL = {J. Pure Appl. Algebra},
  FJOURNAL = {Journal of Pure and Applied Algebra},
    VOLUME = {215},
      YEAR = {2011},
    NUMBER = {4},
     PAGES = {628--641},
      ISSN = {0022-4049},
   MRCLASS = {13F60 (16G20)},
  MRNUMBER = {2738377},
MRREVIEWER = {Kyungyong Lee},
       DOI = {10.1016/j.jpaa.2010.06.012},
       URL = {https://doi.org/10.1016/j.jpaa.2010.06.012},
}

@article{drinfeld2013algebraic,
	title	= {On algebraic spaces with an action of {${\mathbb G}_m$}},
  	author	= {Drinfeld, Vladimir},
  	journal	= {arXiv preprint arXiv:1308.2604},
  	year	= {2013}
}

@article{finkelberg2018comultiplication,
	AUTHOR = {Finkelberg, Michael and Kamnitzer, Joel and Pham, Khoa and Rybnikov, Leonid and Weekes, Alex},
     TITLE = {Comultiplication for shifted {Y}angians and quantum open {T}oda lattice},
   JOURNAL = {Adv. Math.},
  FJOURNAL = {Advances in Mathematics},
    VOLUME = {327},
      YEAR = {2018},
     PAGES = {349--389},
      ISSN = {0001-8708},
   MRCLASS = {17B37 (81T25)},
  MRNUMBER = {3761996},
MRREVIEWER = {Huafeng Zhang},
       DOI = {10.1016/j.aim.2017.06.018},
       URL = {https://doi.org/10.1016/j.aim.2017.06.018},
}

@article {geiss2005semicanonical,
    AUTHOR = {Geiss, Christof and Leclerc, Bernard and Schr\"{o}er, Jan},
     TITLE = {Semicanonical bases and preprojective algebras},
   JOURNAL = {Ann. Sci. \'{E}cole Norm. Sup. (4)},
  FJOURNAL = {Annales Scientifiques de l'\'{E}cole Normale Sup\'{e}rieure. Quatri\`eme
              S\'{e}rie},
    VOLUME = {38},
      YEAR = {2005},
    NUMBER = {2},
     PAGES = {193--253},
      ISSN = {0012-9593},
   MRCLASS = {17B37 (16G20)},
  MRNUMBER = {2144987},
MRREVIEWER = {Andrew W. Hubery},
       DOI = {10.1016/j.ansens.2004.12.001},
       URL = {https://doi.org/10.1016/j.ansens.2004.12.001},
}

@article {geiss2012generic,
    AUTHOR = {Geiss, Christof and Leclerc, Bernard and Schr\"{o}er, Jan},
     TITLE = {Generic bases for cluster algebras and the {C}hamber ansatz},
   JOURNAL = {J. Amer. Math. Soc.},
  FJOURNAL = {Journal of the American Mathematical Society},
    VOLUME = {25},
      YEAR = {2012},
    NUMBER = {1},
     PAGES = {21--76},
      ISSN = {0894-0347},
   MRCLASS = {13F60 (14M15 16G20 20G44)},
  MRNUMBER = {2833478},
MRREVIEWER = {Kyungyong Lee},
       DOI = {10.1090/S0894-0347-2011-00715-7},
       URL = {https://doi.org/10.1090/S0894-0347-2011-00715-7},
}

@article{gerasimov2005class,
    AUTHOR = {Gerasimov, A. and Kharchev, S. and Lebedev, D. and Oblezin, S.},
     TITLE = {On a class of representations of the {Y}angian and moduli space of monopoles},
   JOURNAL = {Comm. Math. Phys.},
  FJOURNAL = {Communications in Mathematical Physics},
    VOLUME = {260},
      YEAR = {2005},
    NUMBER = {3},
     PAGES = {511--525},
      ISSN = {0010-3616},
   MRCLASS = {53C07 (17B37 53D17 53D30)},
  MRNUMBER = {2182434},
MRREVIEWER = {Olivier G. Schiffmann},
       DOI = {10.1007/s00220-005-1417-3},
       URL = {https://doi.org/10.1007/s00220-005-1417-3},
}

@article {gibson2021demazure,
    AUTHOR = {Gibson, Joel},
     TITLE = {A {D}emazure character formula for the product monomial crystal},
   JOURNAL = {Algebr. Comb.},
  FJOURNAL = {Algebraic Combinatorics},
    VOLUME = {4},
      YEAR = {2021},
    NUMBER = {2},
     PAGES = {301--327},
   MRCLASS = {17B37 (05E10)},
  MRNUMBER = {4244375},
MRREVIEWER = {Jacinta Torres},
       DOI = {10.5802/alco.156},
       URL = {https://doi.org/10.5802/alco.156},
}

@article{ginzburg1995perverse,
  title={Perverse sheaves on a loop group and Langlands' duality},
  author={Ginzburg, Victor},
  journal={arXiv preprint arXiv:alg-geom/9511007},
  year={1995}
}

@article{ginzburg2003primitive,
    AUTHOR = {Ginzburg, Victor},
     TITLE = {On primitive ideals},
   JOURNAL = {Selecta Math. (N.S.)},
  FJOURNAL = {Selecta Mathematica. New Series},
    VOLUME = {9},
      YEAR = {2003},
    NUMBER = {3},
     PAGES = {379--407},
      ISSN = {1022-1824},
   MRCLASS = {16S10 (16S30)},
  MRNUMBER = {2006573},
MRREVIEWER = {Anthony Joseph},
       DOI = {10.1007/s00029-003-0338-2},
       URL = {https://doi.org/10.1007/s00029-003-0338-2},
}

@article{ginzburg2003category,
    AUTHOR = {Ginzburg, Victor and Guay, Nicolas and Opdam, Eric and Rouquier, Rapha\"{e}l},
     TITLE = {On the category {$\mathscr O$} for rational {C}herednik algebras},
   JOURNAL = {Invent. Math.},
  FJOURNAL = {Inventiones Mathematicae},
    VOLUME = {154},
      YEAR = {2003},
    NUMBER = {3},
     PAGES = {617--651},
      ISSN = {0020-9910},
   MRCLASS = {20C08 (20F55)},
  MRNUMBER = {2018786},
MRREVIEWER = {Iain G. Gordon},
       DOI = {10.1007/s00222-003-0313-8},
       URL = {https://doi.org/10.1007/s00222-003-0313-8},
}

@article {gross2018canonical,
    AUTHOR = {Gross, Mark and Hacking, Paul and Keel, Sean and Kontsevich, Maxim},
     TITLE = {Canonical bases for cluster algebras},
   JOURNAL = {J. Amer. Math. Soc.},
  FJOURNAL = {Journal of the American Mathematical Society},
    VOLUME = {31},
      YEAR = {2018},
    NUMBER = {2},
     PAGES = {497--608},
      ISSN = {0894-0347},
   MRCLASS = {13F60 (14J33)},
  MRNUMBER = {3758151},
MRREVIEWER = {Ralf Schiffler},
       DOI = {10.1090/jams/890},
       URL = {https://doi.org/10.1090/jams/890},
}

@article {hernandez2016cluster,
    AUTHOR = {Hernandez, David and Leclerc, Bernard},
     TITLE = {A cluster algebra approach to {$q$}-characters of
              {K}irillov-{R}eshetikhin modules},
   JOURNAL = {J. Eur. Math. Soc. (JEMS)},
  FJOURNAL = {Journal of the European Mathematical Society (JEMS)},
    VOLUME = {18},
      YEAR = {2016},
    NUMBER = {5},
     PAGES = {1113--1159},
      ISSN = {1435-9855},
   MRCLASS = {17B67 (13F60)},
  MRNUMBER = {3500832},
MRREVIEWER = {Xueqing Chen},
       DOI = {10.4171/JEMS/609},
       URL = {https://doi.org/10.4171/JEMS/609},
}

@article{hernandez2021shifted,
    AUTHOR = {Hernandez, David and Zhang, Huafeng},
     TITLE = {Shifted {Y}angians and polynomial {$R$}-matrices},
   JOURNAL = {Publ. Res. Inst. Math. Sci.},
  FJOURNAL = {Publications of the Research Institute for Mathematical Sciences},
    VOLUME = {60},
      YEAR = {2024},
    NUMBER = {1},
     PAGES = {1--69},
      ISSN = {0034-5318},
   MRCLASS = {20G42 (16T25 81R50)},
  MRNUMBER = {4803333},
       DOI = {10.4171/prims/60-1-1},
       URL = {https://doi.org/10.4171/prims/60-1-1},
}

@book{humphreys2008representations,
    AUTHOR = {Humphreys, James E.},
     TITLE = {Representations of semisimple {L}ie algebras in the {BGG} category {$\mathscr{O}$}},
    SERIES = {Graduate Studies in Mathematics},
    VOLUME = {94},
 PUBLISHER = {American Mathematical Society, Providence, RI},
      YEAR = {2008},
     PAGES = {xvi+289},
      ISBN = {978-0-8218-4678-0},
   MRCLASS = {17B10},
  MRNUMBER = {2428237},
MRREVIEWER = {Serge M. Skryabin},
       DOI = {10.1090/gsm/094},
       URL = {https://doi.org/10.1090/gsm/094},
}

@article {kamnitzer2026kosh,
	title={Category {$\mathcal{O}$} for truncated shifted Yangians and the bi-infinite Bott-Samelson variety},
  	author={Kamnitzer, Joel and Labelle, Antoine and Leroux-Lapierre, Alexis and Pinet, Th{\'e}o and Weekes, Alex},
  	journal={arXiv preprint arXiv:2607.04480},
  	year={2026}
}

@article {kamnitzer2010mirkovic,
    AUTHOR = {Kamnitzer, Joel},
     TITLE = {Mirkovi\'{c}-{V}ilonen cycles and polytopes},
   JOURNAL = {Ann. of Math. (2)},
  FJOURNAL = {Annals of Mathematics. Second Series},
    VOLUME = {171},
      YEAR = {2010},
    NUMBER = {1},
     PAGES = {245--294},
      ISSN = {0003-486X},
   MRCLASS = {20G05 (14M15)},
  MRNUMBER = {2630039},
MRREVIEWER = {Peter Fiebig},
       DOI = {10.4007/annals.2010.171.245},
       URL = {https://doi.org/10.4007/annals.2010.171.245},
}

@article{kamnitzer2014yangians,
    AUTHOR = {Kamnitzer, Joel and Webster, Ben and Weekes, Alex and Yacobi, Oded},
     TITLE = {Yangians and quantizations of slices in the affine {G}rassmannian},
   JOURNAL = {Algebra Number Theory},
  FJOURNAL = {Algebra \& Number Theory},
    VOLUME = {8},
      YEAR = {2014},
    NUMBER = {4},
     PAGES = {857--893},
      ISSN = {1937-0652},
   MRCLASS = {17B37 (14D24 14M15 20G15 53D55)},
  MRNUMBER = {3248988},
MRREVIEWER = {Christian Ohn},
       DOI = {10.2140/ant.2014.8.857},
       URL = {https://doi.org/10.2140/ant.2014.8.857},
}

@article{kamnitzer2019highest,
    AUTHOR = {Kamnitzer, Joel and Tingley, Peter and Webster, Ben and Weekes, Alex and Yacobi, Oded},
     TITLE = {Highest weights for truncated shifted {Y}angians and product monomial crystals},
   JOURNAL = {J. Comb. Algebra},
  FJOURNAL = {Journal of Combinatorial Algebra},
    VOLUME = {3},
      YEAR = {2019},
    NUMBER = {3},
     PAGES = {237--303},
      ISSN = {2415-6302},
   MRCLASS = {16G20 (05E10 16T99 17B10)},
  MRNUMBER = {4011667},
MRREVIEWER = {Aleksandr Panov},
       DOI = {10.4171/JCA/32},
       URL = {https://doi.org/10.4171/JCA/32},
}

@article{kamnitzer2019category,
    AUTHOR = {Kamnitzer, Joel and Tingley, Peter and Webster, Ben and  Weekes, Alex and Yacobi, Oded},
     TITLE = {On category {$\mathcal O$} for affine {G}rassmannian slices and categorified tensor products},
   JOURNAL = {Proc. Lond. Math. Soc. (3)},
  FJOURNAL = {Proceedings of the London Mathematical Society. Third Series},
    VOLUME = {119},
      YEAR = {2019},
    NUMBER = {5},
     PAGES = {1179--1233},
      ISSN = {0024-6115},
   MRCLASS = {14M15 (17B37 20G42)},
  MRNUMBER = {3968721},
MRREVIEWER = {Huafeng Zhang},
       DOI = {10.1112/plms.12254},
       URL = {https://doi.org/10.1112/plms.12254},
}

@article{kamnitzer2022lie,
    AUTHOR = {Kamnitzer, Joel and Webster, Ben and Weekes, Alex and Yacobi, Oded},
     TITLE = {Lie algebra actions on module categories for truncated shifted yangians},
   JOURNAL = {Forum Math. Sigma},
  FJOURNAL = {Forum of Mathematics. Sigma},
    VOLUME = {12},
      YEAR = {2024},
     PAGES = {Paper No. e18, 69},
   MRCLASS = {20G05 (14M15 17B67)},
  MRNUMBER = {4699878},
MRREVIEWER = {Linliang Song},
       DOI = {10.1017/fms.2024.3},
       URL = {https://doi.org/10.1017/fms.2024.3},
}

@article{kamnitzer2022hamiltonian,
    AUTHOR = {Kamnitzer, Joel and Pham, Khoa and Weekes, Alex},
     TITLE = {Hamiltonian reduction for affine {G}rassmannian slices and truncated shifted {Y}angians},
   JOURNAL = {Adv. Math.},
  FJOURNAL = {Advances in Mathematics},
    VOLUME = {399},
      YEAR = {2022},
     PAGES = {Paper No. 108281, 52},
      ISSN = {0001-8708},
   MRCLASS = {14M15 (17B10)},
  MRNUMBER = {4385132},
MRREVIEWER = {Felipe Zald\'{\i}var},
       DOI = {10.1016/j.aim.2022.108281},
       URL = {https://doi.org/10.1016/j.aim.2022.108281},
}

@incollection{kamnitzer2022perfect,
    AUTHOR = {Kamnitzer, Joel},
     TITLE = {Perfect bases in representation theory: three mountains and their springs},
 BOOKTITLE = {I{CM}---{I}nternational {C}ongress of {M}athematicians. {V}ol. 4. {S}ections 5--8},
     PAGES = {2976--2996},
 PUBLISHER = {EMS Press, Berlin},
      YEAR = {[2023] \copyright 2023},
   MRCLASS = {22E46 (05E10 13F60 14M15 16G20)},
  MRNUMBER = {4680349},
}

@article{kang2012categorification,
    AUTHOR = {Kang, Seok-Jin and Kashiwara, Masaki},
     TITLE = {Categorification of highest weight modules via\linebreak {K}hovanov-{L}auda-{R}ouquier algebras},
   JOURNAL = {Invent. Math.},
  FJOURNAL = {Inventiones Mathematicae},
    VOLUME = {190},
      YEAR = {2012},
    NUMBER = {3},
     PAGES = {699--742},
      ISSN = {0020-9910},
   MRCLASS = {17B67 (17B37)},
  MRNUMBER = {2995184},
MRREVIEWER = {Volodymyr Mazorchuk},
       DOI = {10.1007/s00222-012-0388-1},
       URL = {https://doi.org/10.1007/s00222-012-0388-1},
}

@incollection {kashiwara2003realizations,
    AUTHOR = {Kashiwara, Masaki},
     TITLE = {Realizations of crystals},
 BOOKTITLE = {Combinatorial and geometric representation theory ({S}eoul, 2001)},
    SERIES = {Contemp. Math.},
    VOLUME = {325},
     PAGES = {133--139},
 PUBLISHER = {Amer. Math. Soc., Providence, RI},
      YEAR = {2003},
   MRCLASS = {17B37 (17B10)},
  MRNUMBER = {1988989},
MRREVIEWER = {Jae-Hoon Kwon},
       DOI = {10.1090/conm/325/05668},
       URL = {https://doi.org/10.1090/conm/325/05668},
}

@article{khovanov2009diagrammatic,
    AUTHOR = {Khovanov, Mikhail and Lauda, Aaron D.},
     TITLE = {A diagrammatic approach to categorification of quantum groups {I}},
   JOURNAL = {Represent. Theory},
  FJOURNAL = {Representation Theory. An Electronic Journal of the American
              Mathematical Society},
    VOLUME = {13},
      YEAR = {2009},
     PAGES = {309--347},
   MRCLASS = {17B37},
  MRNUMBER = {2525917},
MRREVIEWER = {Fan Xu},
       DOI = {10.1090/S1088-4165-09-00346-X},
       URL = {https://doi.org/10.1090/S1088-4165-09-00346-X},
}

@article{khovanov2011diagrammatic,
    AUTHOR = {Khovanov, Mikhail and Lauda, Aaron D.},
     TITLE = {A diagrammatic approach to categorification of quantum groups {II}},
   JOURNAL = {Trans. Amer. Math. Soc.},
  FJOURNAL = {Transactions of the American Mathematical Society},
    VOLUME = {363},
      YEAR = {2011},
    NUMBER = {5},
     PAGES = {2685--2700},
      ISSN = {0002-9947},
   MRCLASS = {17B37 (16T20)},
  MRNUMBER = {2763732},
MRREVIEWER = {Volodymyr Mazorchuk},
       DOI = {10.1090/S0002-9947-2010-05210-9},
       URL = {https://doi.org/10.1090/S0002-9947-2010-05210-9},
}

@article {kleshchev2010homogeneous,
    AUTHOR = {Kleshchev, Alexander and Ram, Arun},
     TITLE = {Homogeneous representations of {K}hovanov-{L}auda algebras},
   JOURNAL = {J. Eur. Math. Soc. (JEMS)},
  FJOURNAL = {Journal of the European Mathematical Society (JEMS)},
    VOLUME = {12},
      YEAR = {2010},
    NUMBER = {5},
     PAGES = {1293--1306},
      ISSN = {1435-9855},
   MRCLASS = {20C08 (17B67)},
  MRNUMBER = {2677617},
MRREVIEWER = {Andrew Mathas},
       DOI = {10.4171/JEMS/230},
       URL = {https://doi.org/10.4171/JEMS/230},
}

@article {kleshchev2011representations,
    AUTHOR = {Kleshchev, Alexander and Ram, Arun},
     TITLE = {Representations of {K}hovanov-{L}auda-{R}ouquier algebras and
              combinatorics of {L}yndon words},
   JOURNAL = {Math. Ann.},
  FJOURNAL = {Mathematische Annalen},
    VOLUME = {349},
      YEAR = {2011},
    NUMBER = {4},
     PAGES = {943--975},
      ISSN = {0025-5831},
   MRCLASS = {16S99 (16G10 20C08)},
  MRNUMBER = {2777040},
MRREVIEWER = {Selvaraj Chelliah},
       DOI = {10.1007/s00208-010-0543-1},
       URL = {https://doi.org/10.1007/s00208-010-0543-1},
}

@article {knutson1999littelmann,
    AUTHOR = {Knutson, Allen},
     TITLE = {A {L}ittelmann-type formula for {D}uistermaat-{H}eckman
              measures},
   JOURNAL = {Invent. Math.},
  FJOURNAL = {Inventiones Mathematicae},
    VOLUME = {135},
      YEAR = {1999},
    NUMBER = {1},
     PAGES = {185--200},
      ISSN = {0020-9910},
   MRCLASS = {53D20 (37J15)},
  MRNUMBER = {1664699},
MRREVIEWER = {Christopher T. Woodward},
       DOI = {10.1007/s002220050283},
       URL = {https://doi.org/10.1007/s002220050283},
}

@book{krause2000growth,    
    AUTHOR = {Krause, G\"{u}nter R. and Lenagan, Thomas H.},
     TITLE = {Growth of algebras and {G}elfand-{K}irillov dimension},
    SERIES = {Graduate Studies in Mathematics},
    VOLUME = {22},
   EDITION = {Revised},
 PUBLISHER = {American Mathematical Society, Providence, RI},
      YEAR = {2000},
     PAGES = {x+212},
      ISBN = {0-8218-0859-1},
   MRCLASS = {16P90},
  MRNUMBER = {1721834},
MRREVIEWER = {Martha K. Smith},
       DOI = {10.1090/gsm/022},
       URL = {https://doi.org/10.1090/gsm/022},
}

@article{krylov2018integrable,
    AUTHOR = {Krylov, V. V.},
     TITLE = {Integrable crystals and a restriction on a {L}evi subgroup via
              generalized slices in the affine {G}rassmannian},
   JOURNAL = {Funktsional. Anal. i Prilozhen.},
  FJOURNAL = {Funktsional\cprime ny\u{\i} Analiz i ego Prilozheniya},
    VOLUME = {52},
      YEAR = {2018},
    NUMBER = {2},
     PAGES = {40--65},
      ISSN = {0374-1990},
   MRCLASS = {14D24 (20G05)},
  MRNUMBER = {3799411},
MRREVIEWER = {Dmitry V. Artamonov},
       DOI = {10.1007/s10688-018-0217-4},
       URL = {https://doi.org/10.1007/s10688-018-0217-4},
}

@article{krylov2021almost,
    AUTHOR = {Krylov, Vasily and Perunov, Ivan},
     TITLE = {Almost dominant generalized slices and convolution diagrams
              over them},
   JOURNAL = {Adv. Math.},
  FJOURNAL = {Advances in Mathematics},
    VOLUME = {392},
      YEAR = {2021},
     PAGES = {Paper No. 108034, 45},
      ISSN = {0001-8708},
   MRCLASS = {14M15 (20G05)},
  MRNUMBER = {4316674},
MRREVIEWER = {Ryan David Kinser},
       DOI = {10.1016/j.aim.2021.108034},
       URL = {https://doi.org/10.1016/j.aim.2021.108034},
}

@article{kashiwara2010deformation,
    AUTHOR = {Kashiwara, Masaki and Schapira, Pierre},
     TITLE = {Deformation quantization modules},
   JOURNAL = {Ast\'{e}risque},
  FJOURNAL = {Ast\'{e}risque},
    NUMBER = {345},
      YEAR = {2012},
     PAGES = {xii+147},
      ISSN = {0303-1179},
      ISBN = {978-2-85629-345-4},
   MRCLASS = {53D55 (14F43 17B63 53D17)},
  MRNUMBER = {3012169},
MRREVIEWER = {Andrea D'Agnolo},
}

@article{kashiwara1993global,
    AUTHOR = {Kashiwara, Masaki},
     TITLE = {Global crystal bases of quantum groups},
   JOURNAL = {Duke Math. J.},
  FJOURNAL = {Duke Mathematical Journal},
    VOLUME = {69},
      YEAR = {1993},
    NUMBER = {2},
     PAGES = {455--485},
      ISSN = {0012-7094},
   MRCLASS = {17B37 (17B10 81R50)},
  MRNUMBER = {1203234},
MRREVIEWER = {Kailash C. Misra},
       DOI = {10.1215/S0012-7094-93-06920-7},
       URL = {https://doi.org/10.1215/S0012-7094-93-06920-7},
}

@article{lauda2011crystals,
    AUTHOR = {Lauda, Aaron D. and Vazirani, Monica},
     TITLE = {Crystals from categorified quantum groups},
   JOURNAL = {Adv. Math.},
  FJOURNAL = {Advances in Mathematics},
    VOLUME = {228},
      YEAR = {2011},
    NUMBER = {2},
     PAGES = {803--861},
      ISSN = {0001-8708},
   MRCLASS = {17B37},
  MRNUMBER = {2822211},
MRREVIEWER = {Peter W. Tingley},
       DOI = {10.1016/j.aim.2011.06.009},
       URL = {https://doi.org/10.1016/j.aim.2011.06.009},
}

@article{losev2017categories,
    AUTHOR = {Losev, Ivan},
     TITLE = {On categories {$\mathcal O$} for quantized symplectic resolutions},
   JOURNAL = {Compos. Math.},
  FJOURNAL = {Compositio Mathematica},
    VOLUME = {153},
      YEAR = {2017},
    NUMBER = {12},
     PAGES = {2445--2481},
      ISSN = {0010-437X},
   MRCLASS = {16E35 (16G99 18F99 32C99 53D50)},
  MRNUMBER = {3705295},
MRREVIEWER = {Vida Milani},
       DOI = {10.1112/S0010437X17007382},
       URL = {https://doi.org/10.1112/S0010437X17007382},
}

@article{losev2021localization,
  title={Localization theorems for quantized symplectic resolutions},
  author={Losev, Ivan},
  journal={arXiv:2103.11193},
  year={2021}
}

@article{lusztig1990canonicalI,
    AUTHOR = {Lusztig, G.},
     TITLE = {Canonical bases arising from quantized enveloping algebras},
   JOURNAL = {J. Amer. Math. Soc.},
  FJOURNAL = {Journal of the American Mathematical Society},
    VOLUME = {3},
      YEAR = {1990},
    NUMBER = {2},
     PAGES = {447--498},
      ISSN = {0894-0347},
   MRCLASS = {17B35 (16A64)},
  MRNUMBER = {1035415},
MRREVIEWER = {Ya. S. So\u{\i}bel\cprime man},
       DOI = {10.2307/1990961},
       URL = {https://doi.org/10.2307/1990961},
}

@article{lusztig1990canonicalII,
    AUTHOR = {Lusztig, G.},
     TITLE = {Canonical bases arising from quantized enveloping algebras.
              {II}},
      NOTE = {Common trends in mathematics and quantum field theories
              (Kyoto, 1990)},
   JOURNAL = {Progr. Theoret. Phys. Suppl.},
  FJOURNAL = {Progress of Theoretical Physics. Supplement},
    NUMBER = {102},
      YEAR = {1990},
     PAGES = {175--201 (1991)},
      ISSN = {0375-9687},
   MRCLASS = {17B35},
  MRNUMBER = {1182165},
MRREVIEWER = {Ya. S. So\u{\i}bel\cprime man},
       DOI = {10.1143/PTPS.102.175},
       URL = {https://doi.org/10.1143/PTPS.102.175},
}

@article {lusztig2000semicanonical,
    AUTHOR = {Lusztig, G.},
     TITLE = {Semicanonical bases arising from enveloping algebras},
   JOURNAL = {Adv. Math.},
  FJOURNAL = {Advances in Mathematics},
    VOLUME = {151},
      YEAR = {2000},
    NUMBER = {2},
     PAGES = {129--139},
      ISSN = {0001-8708},
   MRCLASS = {17B67 (16S30)},
  MRNUMBER = {1758244},
MRREVIEWER = {William M. McGovern},
       DOI = {10.1006/aima.1999.1873},
       URL = {https://doi.org/10.1006/aima.1999.1873},
}

@article{lusztig1983singularities,
    AUTHOR = {Lusztig, George},
     TITLE = {Singularities, character formulas, and a {$q$}-analog of weight multiplicities},
 BOOKTITLE = {Analysis and topology on singular spaces, {II}, {III} ({L}uminy, 1981)},
    SERIES = {Ast\'{e}risque},
    VOLUME = {101},
     PAGES = {208--229},
 PUBLISHER = {Soc. Math. France, Paris},
      YEAR = {1983},
   MRCLASS = {17B10 (05A30 20G05 22E47)},
  MRNUMBER = {737932},
MRREVIEWER = {James E. Humphreys},
}

@article {mandel2017theta,
    AUTHOR = {Mandel, Travis},
     TITLE = {Theta bases are atomic},
   JOURNAL = {Compos. Math.},
  FJOURNAL = {Compositio Mathematica},
    VOLUME = {153},
      YEAR = {2017},
    NUMBER = {6},
     PAGES = {1217--1219},
      ISSN = {0010-437X},
   MRCLASS = {13F60},
  MRNUMBER = {3705255},
MRREVIEWER = {Man-Wai Cheung},
       DOI = {10.1112/S0010437X17007060},
       URL = {https://doi.org/10.1112/S0010437X17007060},
}

@article{mcconnell1989gelfand,
    AUTHOR = {McConnell, J. C. and Stafford, J. T.},
     TITLE = {Gelfand-{K}irillov dimension and associated graded modules},
   JOURNAL = {J. Algebra},
  FJOURNAL = {Journal of Algebra},
    VOLUME = {125},
      YEAR = {1989},
    NUMBER = {1},
     PAGES = {197--214},
      ISSN = {0021-8693},
   MRCLASS = {16A03},
  MRNUMBER = {1012671},
MRREVIEWER = {Arno van den Essen},
       DOI = {10.1016/0021-8693(89)90301-3},
       URL = {https://doi.org/10.1016/0021-8693(89)90301-3},
}

@article {mcnamara2015finite,
    AUTHOR = {McNamara, Peter J.},
     TITLE = {Finite dimensional representations of {K}hovanov-{L}auda-{R}ouquier algebras {I}: {F}inite type},
   JOURNAL = {J. Reine Angew. Math.},
  FJOURNAL = {Journal f\"{u}r die Reine und Angewandte Mathematik. [Crelle's Journal]},
    VOLUME = {707},
      YEAR = {2015},
     PAGES = {103--124},
      ISSN = {0075-4102},
   MRCLASS = {17B10 (17B22 17B37)},
  MRNUMBER = {3403455},
MRREVIEWER = {Volodymyr Mazorchuk},
       DOI = {10.1515/crelle-2013-0075},
       URL = {https://doi.org/10.1515/crelle-2013-0075},
}

@article {mcnamara2019folding,
    AUTHOR = {McNamara, Peter J.},
     TITLE = {Folding {KLR} algebras},
   JOURNAL = {J. Lond. Math. Soc. (2)},
  FJOURNAL = {Journal of the London Mathematical Society. Second Series},
    VOLUME = {100},
      YEAR = {2019},
    NUMBER = {2},
     PAGES = {447--469},
      ISSN = {0024-6107},
   MRCLASS = {17B37},
  MRNUMBER = {4017150},
MRREVIEWER = {Changchang Xi},
       DOI = {10.1112/jlms.12218},
       URL = {https://doi.org/10.1112/jlms.12218},
}

@book{miller2005combinatorial,
    AUTHOR = {Miller, Ezra and Sturmfels, Bernd},
     TITLE = {Combinatorial commutative algebra},
    SERIES = {Graduate Texts in Mathematics},
    VOLUME = {227},
 PUBLISHER = {Springer-Verlag, New York},
      YEAR = {2005},
     PAGES = {xiv+417},
      ISBN = {0-387-22356-8},
   MRCLASS = {13-01 (05-01 05E99 13D02 14M15 14M25)},
  MRNUMBER = {2110098},
MRREVIEWER = {Joseph Gubeladze},
}

@article{mirkovic2007geometric,
    AUTHOR = {Mirkovi\'{c}, I. and Vilonen, K.},
     TITLE = {Geometric {L}anglands duality and representations of algebraic
              groups over commutative rings},
   JOURNAL = {Ann. of Math. (2)},
  FJOURNAL = {Annals of Mathematics. Second Series},
    VOLUME = {166},
      YEAR = {2007},
    NUMBER = {1},
     PAGES = {95--143},
      ISSN = {0003-486X},
   MRCLASS = {22E55 (11R39 20G05)},
  MRNUMBER = {2342692},
MRREVIEWER = {Peter Fiebig},
       DOI = {10.4007/annals.2007.166.95},
       URL = {https://doi.org/10.4007/annals.2007.166.95},
}

@article {muthiah2021weyl,
    AUTHOR = {Muthiah, Dinakar},
     TITLE = {Weyl group action on weight zero {M}irkovi\'{c}-{V}ilonen basis
              and equivariant multiplicities},
   JOURNAL = {Adv. Math.},
  FJOURNAL = {Advances in Mathematics},
    VOLUME = {385},
      YEAR = {2021},
     PAGES = {Paper No. 107793, 40},
      ISSN = {0001-8708},
   MRCLASS = {22E57 (14M15 17B10 20G20)},
  MRNUMBER = {4261167},
MRREVIEWER = {Spencer Leslie},
       DOI = {10.1016/j.aim.2021.107793},
       URL = {https://doi.org/10.1016/j.aim.2021.107793},
}

@article {nakada2008colored,
    AUTHOR = {Nakada, Kento},
     TITLE = {Colored hook formula for a generalized {Y}oung diagram},
   JOURNAL = {Osaka J. Math.},
  FJOURNAL = {Osaka Journal of Mathematics},
    VOLUME = {45},
      YEAR = {2008},
    NUMBER = {4},
     PAGES = {1085--1120},
      ISSN = {0030-6126},
   MRCLASS = {17B67 (05E10 17B22)},
  MRNUMBER = {2493972},
MRREVIEWER = {Hiro-Fumi Yamada},
       URL = {http://projecteuclid.org/euclid.ojm/1227708835},
}

@article{nakajima2003tanalogues,
    AUTHOR = {Nakajima, Hiraku},
     TITLE = {{$t$}-analogs of {$q$}-characters of quantum affine algebras
              of type {$A_n,D_n$}},
 BOOKTITLE = {Combinatorial and geometric representation theory ({S}eoul, 2001)},
    SERIES = {Contemp. Math.},
    VOLUME = {325},
     PAGES = {141--160},
 PUBLISHER = {Amer. Math. Soc., Providence, RI},
      YEAR = {2003},
   MRCLASS = {17B37 (81R50 82B23)},
  MRNUMBER = {1988990},
MRREVIEWER = {\c{S}erban Raianu},
       DOI = {10.1090/conm/325/05669},
       URL = {https://doi.org/10.1090/conm/325/05669},
}

@article {nakajima1994instantons,
    AUTHOR = {Nakajima, Hiraku},
     TITLE = {Instantons on {ALE} spaces, quiver varieties, and {K}ac-{M}oody algebras},
   JOURNAL = {Duke Math. J.},
  FJOURNAL = {Duke Mathematical Journal},
    VOLUME = {76},
      YEAR = {1994},
    NUMBER = {2},
     PAGES = {365--416},
      ISSN = {0012-7094},
   MRCLASS = {53C25 (17B67 58D27 58E15)},
  MRNUMBER = {1302318},
MRREVIEWER = {Andrew Dancer},
       DOI = {10.1215/S0012-7094-94-07613-8},
       URL = {https://doi.org/10.1215/S0012-7094-94-07613-8},
}

@article {nakajima2023coulomb,
    AUTHOR = {Nakajima, Hiraku and Weekes, Alex},
     TITLE = {Coulomb branches of quiver gauge theories with symmetrizers},
   JOURNAL = {J. Eur. Math. Soc. (JEMS)},
  FJOURNAL = {Journal of the European Mathematical Society (JEMS)},
    VOLUME = {25},
      YEAR = {2023},
    NUMBER = {1},
     PAGES = {203--230},
      ISSN = {1435-9855},
   MRCLASS = {14D21 (16G20 17B37 81T13)},
  MRNUMBER = {4556783},
MRREVIEWER = {\^{A}ngela Mestre},
       DOI = {10.4171/jems/1176},
       URL = {https://doi.org/10.4171/jems/1176},
}

@article {qin2017triangular,
    AUTHOR = {Qin, Fan},
     TITLE = {Triangular bases in quantum cluster algebras and monoidal categorification conjectures},
   JOURNAL = {Duke Math. J.},
  FJOURNAL = {Duke Mathematical Journal},
    VOLUME = {166},
      YEAR = {2017},
    NUMBER = {12},
     PAGES = {2337--2442},
      ISSN = {0012-7094},
   MRCLASS = {13F60 (81R50)},
  MRNUMBER = {3694569},
MRREVIEWER = {Alfredo N\'{a}jera Ch\'{a}vez},
       DOI = {10.1215/00127094-2017-0006},
       URL = {https://doi.org/10.1215/00127094-2017-0006},
}

@article {qin2024bases,
    AUTHOR = {Qin, Fan},
     TITLE = {Bases for upper cluster algebras and tropical points},
   JOURNAL = {J. Eur. Math. Soc. (JEMS)},
  FJOURNAL = {Journal of the European Mathematical Society (JEMS)},
    VOLUME = {26},
      YEAR = {2024},
    NUMBER = {4},
     PAGES = {1255--1312},
      ISSN = {1435-9855},
   MRCLASS = {13F60},
  MRNUMBER = {4721032},
       DOI = {10.4171/jems/1308},
       URL = {https://doi.org/10.4171/jems/1308},
}

@incollection {qin2023cluster,
    AUTHOR = {Qin, Fan},
     TITLE = {Cluster algebras and their bases},
 BOOKTITLE = {Representations of algebras and related structures},
    SERIES = {EMS Ser. Congr. Rep.},
     PAGES = {335--369},
 PUBLISHER = {EMS Press, Berlin},
      YEAR = {2023},
   MRCLASS = {13F60},
  MRNUMBER = {4693645},
}

@article{rouquier20082,
	title	={2-kac-moody algebras},
  	author	={Rouquier, Rapha{\"e}l},
  	journal	={arXiv preprint arXiv:0812.5023},
  	year	={2008}
}

@article{soergel1986equivalences,
    AUTHOR = {Soergel, Wolfgang},
     TITLE = {\'{E}quivalences de certaines cat\'{e}gories de {${\mathfrak g}$}-modules},
   JOURNAL = {C. R. Acad. Sci. Paris S\'{e}r. I Math.},
  FJOURNAL = {Comptes Rendus des S\'{e}ances de l'Acad\'{e}mie des Sciences. S\'{e}rie I. Math\'{e}matique},
    VOLUME = {303},
      YEAR = {1986},
    NUMBER = {15},
     PAGES = {725--728},
      ISSN = {0249-6291},
   MRCLASS = {17B10 (22E47)},
  MRNUMBER = {872544},
MRREVIEWER = {James E. Humphreys},
}

@article {stembridge2001minuscule,
    AUTHOR = {Stembridge, John R.},
     TITLE = {Minuscule elements of {W}eyl groups},
   JOURNAL = {J. Algebra},
  FJOURNAL = {Journal of Algebra},
    VOLUME = {235},
      YEAR = {2001},
    NUMBER = {2},
     PAGES = {722--743},
      ISSN = {0021-8693,1090-266X},
   MRCLASS = {20F55 (17B67)},
  MRNUMBER = {1805477},
MRREVIEWER = {G\"{o}tz\ Pfeiffer},
       DOI = {10.1006/jabr.2000.8488},
       URL = {https://doi.org/10.1006/jabr.2000.8488},
}

@article{tauvel1982dimension,
    AUTHOR = {Tauvel, Patrice},
     TITLE = {Sur la dimension de {G}elfand-{K}irillov},
   JOURNAL = {Comm. Algebra},
  FJOURNAL = {Communications in Algebra},
    VOLUME = {10},
      YEAR = {1982},
    NUMBER = {9},
     PAGES = {939--963},
      ISSN = {0092-7872},
   MRCLASS = {16A55 (16A64 17B35)},
  MRNUMBER = {654050},
MRREVIEWER = {Martha K. Smith},
       DOI = {10.1080/00927878208822758},
       URL = {https://doi.org/10.1080/00927878208822758},
}

@article{varagnolo2011canonical,
    AUTHOR = {Varagnolo, M. and Vasserot, E.},
     TITLE = {Canonical bases and {KLR}-algebras},
   JOURNAL = {J. Reine Angew. Math.},
  FJOURNAL = {Journal f\"{u}r die Reine und Angewandte Mathematik. [Crelle's Journal]},
    VOLUME = {659},
      YEAR = {2011},
     PAGES = {67--100},
      ISSN = {0075-4102},
   MRCLASS = {17B37 (16T20)},
  MRNUMBER = {2837011},
MRREVIEWER = {Nicolas Jacon},
       DOI = {10.1515/CRELLE.2011.068},
       URL = {https://doi.org/10.1515/CRELLE.2011.068},
}

@article{varagnolo2025representations,
  title={Representations of shifted affine quantum groups and Coulomb branches},
  author={Varagnolo, Michela and Vasserot, Eric},
  journal={arXiv preprint arXiv:2503.06262},
  year={2025}
}

@article{webster2017knot,
    AUTHOR = {Webster, Ben},
     TITLE = {Knot invariants and higher representation theory},
   JOURNAL = {Mem. Amer. Math. Soc.},
  FJOURNAL = {Memoirs of the American Mathematical Society},
    VOLUME = {250},
      YEAR = {2017},
    NUMBER = {1191},
     PAGES = {v+141},
      ISSN = {0065-9266},
      ISBN = {978-1-4704-2650-7; 978-1-4704-4206-4},
   MRCLASS = {57M27 (17B10 18D05 57M25)},
  MRNUMBER = {3709726},
MRREVIEWER = {Stefan K. Friedl},
       DOI = {10.1090/memo/1191},
       URL = {https://doi.org/10.1090/memo/1191},
}

@article{webster2020quantum,
    AUTHOR = {Webster, Ben and Weekes, Alex and Yacobi, Oded},
     TITLE = {A quantum {M}irkovi\'{c}-{V}ybornov isomorphism},
   JOURNAL = {Represent. Theory},
  FJOURNAL = {Representation Theory. An Electronic Journal of the American Mathematical Society},
    VOLUME = {24},
      YEAR = {2020},
     PAGES = {38--84},
   MRCLASS = {17B37 (16S80 20C99)},
  MRNUMBER = {4052554},
MRREVIEWER = {Aleksandr Panov},
       DOI = {10.1090/ert/536},
       URL = {https://doi.org/10.1090/ert/536},
}

@phdthesis{weekes2016highest,
    AUTHOR = {Weekes, Alex},
     TITLE = {Highest {W}eights for {T}runcated {S}hifted {Y}angians},
      NOTE = {Thesis (Ph.D.)--University of Toronto (Canada)},
      YEAR = {2016},
     PAGES = {102},
      ISBN = {978-1369-85412-1},
   MRCLASS = {Thesis},
  MRNUMBER = {3697595},
}

@incollection {zhu2017introduction,
    AUTHOR = {Zhu, Xinwen},
     TITLE = {An introduction to affine {G}rassmannians and the geometric {S}atake equivalence},
 BOOKTITLE = {Geometry of moduli spaces and representation theory},
    SERIES = {IAS/Park City Math. Ser.},
    VOLUME = {24},
     PAGES = {59--154},
 PUBLISHER = {Amer. Math. Soc., Providence, RI},
      YEAR = {2017},
   MRCLASS = {14M15 (14D24 20F65 22E57)},
  MRNUMBER = {3752460},
MRREVIEWER = {Felipe Zald\'{\i}var},
}
}

\end{document}